\documentclass{siamart1116}

\usepackage{lipsum}
\usepackage{amssymb,amsfonts,amsopn,amsmath}
\usepackage[utf8]{inputenc}
\usepackage[english]{babel}

\usepackage[T1]{fontenc}
\usepackage{ae,aecompl}

\usepackage{bm}

\usepackage{scrtime}

\usepackage{setspace}
\usepackage{url}
\usepackage{float}
\usepackage{placeins}

\usepackage{multirow}

\usepackage{graphicx}
\usepackage{ulem}  % for the command \sout{}

\usepackage{epstopdf}
\usepackage{algorithmic}
\ifpdf
  \DeclareGraphicsExtensions{.eps,.pdf,.png,.jpg}
\else
  \DeclareGraphicsExtensions{.eps}
\fi

\numberwithin{theorem}{section}

\newcommand{\TheTitle}{APPROXIMATIONS AND INTERPOLATIONS} 
\newcommand{\TheAuthors}{K. ZHANG, E. CROOKS, and A. ORLANDO}

\headers{\TheTitle}{\TheAuthors}

\title{{Compensated Convexity Methods\\ for Approximations and Interpolations\\ 
of Sampled Functions in Euclidean Spaces: \\ 
Applications to Contour Lines, Sparse Data and Inpainting}\thanks{Submitted to the editors DATE.
}}
\author{
	Kewei Zhang\thanks{School  of Mathematical Sciences,
		University of Nottingham, University Park, Nottingham, NG7 2RD, UK 
		(\email{kewei.zhang@nottingham.ac.uk}).} 
	\and
	Elaine Crooks\thanks{Department of Mathematics, Swansea University,
		Singleton Park, Swansea, SA2 8PP, UK
		(\email{e.c.m.crooks@swansea.ac.uk}).}
	\and
	Antonio Orlando\thanks{CONICET, Departamento de Bioingenier\'{i}a, FACET,
		Universidad Nacional de Tucum\'an, Argentina
		(\email{aorlando@herrera.unt.edu.ar}).}
	}

\ifpdf
\hypersetup{
  pdftitle={\TheTitle},
  pdfauthor={\TheAuthors}
}
\fi

\numberwithin{equation}{section}
\definecolor{purple}{RGB}{160,32,40}

\usepackage{ulem}  % for the command \sout{}

\newtheorem{teo}{Theorem}[section]
\newtheorem{nota}[teo]{Remark}
\newtheorem{ex}[teo]{Example}
\newtheorem{coro}[teo]{Corollary}
\newtheorem{defi}[teo]{Definition}
\newtheorem{prop}[teo]{Proposition}
\newtheorem{lema}[teo]{Lemma}

\newcommand{\R}{\ensuremath{\mathbb{R}} }

\newcommand{\dist}{\mathrm{dist}}

\newcommand{\sign}{\mathrm{sign}}

\DeclareMathOperator{\co}{\mathsf{co}}
\DeclareMathOperator{\Span}{\mathsf{span}}
\DeclareMathOperator{\diam}{\mathsf{diam}}

\DeclareMathOperator{\ddiv}{\mathsf{div}}

\begin{document}

%%%%%%%%%%%%%%%%%%%%%%%%%%%%%%%

\maketitle

%%%%%%%%%%%%%%%%%%%%%%%%%%%%%%%

\begin{abstract}  
	This paper is concerned with applications of the theory of approximation and interpolation based on 
	compensated convex transforms developed in \cite{ZCO16}.  
	We apply our methods to 
	$(i)$ surface reconstruction starting from the knowledge of finitely many level sets (or `contour lines'); 
	$(ii)$ scattered data approximation;  
	$(iii)$ image inpainting. 
	For $(i)$ and $(ii)$ our methods give interpolations. For the case of finite sets (scattered data),
	in particular, our approximations  provide a natural triangulation and piecewise affine interpolation.  
	Prototype examples of explicitly calculated approximations and inpainting results are
	presented for both finite and compact sets.  We also show numerical experiments for
	applications of our methods to high density salt \& pepper noise reduction in image 
	processing, for image inpainting and for approximation and interpolations of continuous 
	functions sampled on finitely many level sets and on scattered points. 
\end{abstract}

% REQUIRED
\begin{keywords}
  compensated convex transforms, scattered data, contour lines, interpolation, 
  approximation, inpainting, Hausdorff stability,
  maximum principle, convex density radius, image inpainting, high density salt \& pepper noise reduction
\end{keywords}

% REQUIRED
\begin{AMS}
  90C25, 90C26, 49J52, 52A41, 65K10
\end{AMS}

%%%%%%%%%%%%%%%%%%%%%%%%%%%%%%%%%%%%%%%%%%%%%%%%%%%%%%%%%%%%%%%%%%%%%%%%%%%%%%%%%%%%%%%%%%%%%%%%%%%%%
%%%%%%%%%%%%%%%%%%%%%%%%%%%%%%%%%%%%%%%%%%%%%%%%%%%%%%%%%%%%%%%%%%%%%%%%%%%%%%%%%%%%%%%%%%%%%%%%%%%%%
%%%%%%%%%%%%%%%%%%%%%%%%%%%%%%%%%%%%%%%%%%%%%%%%%%%%%%%%%%%%%%%%%%%%%%%%%%%%%%%%%%%%%%%%%%%%%%%%%%%%%

\setcounter{equation}{0}
\section{Introduction}

This paper is concerned with the application of  the compensated-convexity based theory 
for approximation and interpolation of sampled functions
that was presented in our previous article \cite{ZCO16} to 
surface reconstruction based on knowledge from finitely many level sets, scattered data approximation, 
and image inpainting.  

\medskip\noindent 
In general, approximation theory is concerned with the problem of finding in the set of
simple known functions one that is close in some sense to a more
complicated otherwise unknown function. The variational theory is developed by specifying a 
priori the class of the approximating functions and the criteria that allow selecting an element of such class. 
In the implementation of the theory, the approximating functions generally depend on unknowns parameters that control 
their form, so that the problem boils down to selecting the parameters that allow meeting 
the chosen criteria. Such criteria are usually related to the error between the approximating functions and 
what is known about the function to be approximated and might contain some regularizing term that determines 
the regularity of the approximating function and makes the whole problem well posed.

\medskip\noindent 
Different classes of approximating functions, such as,  for instance, algebraic polynomials \cite{Tre13}, 
trigonometric polynomials \cite{Tre00,Tre13}, radial basis functions \cite{Wen05,Buh04,Fas07}, 
continuous piecewise polynomials \cite{Pow81}, have been considered, and while their definition is usually motivated 
by good approximating properties for a given field of application, on the other hand the specific nature 
of a class of functions also represents a restriction that limits their general application. 

\medskip \noindent 
Total variation-type models \cite{ROF92,BKP10}, \cite[Ch. 6]{CS05} and  
geometric partial differential equations \cite{CMS98},\cite[Ch. 1]{Wei98},\cite[Ch. 8]{Sap01} have
also been used as interpolation models. Their use has been principally motivated by applications in the field 
of image processing and geoscience. We mention in particular the applications to salt \& pepper noise
reduction \cite{CHN05}, image inpainting (by using TV-inpainting models \cite{BFW16,Get12},\cite[Ch. 6]{CS05}, 
Curvature Diffusion Driven inpainting model \cite{CS01}, geometric PDE based inpainting model \cite{BSCB00} 
or other PDE-based models discussed in the monograph \cite{Sch15})
and image interpolation \cite{BBBW09,CMS98,Get11}, among others. For the applications
to geoscience, and in particular to the construction of digital elevation models,
PDE based interpolation models, such as the one considered in 
\cite{ACGR02}, where the interpolant is sought as 
the absolutely minimizing Lipschitz extension \cite{Aro67,Jen93} of the known values, 
have also been proposed and shown to be 
competitive against the classical interpolation methods such as the geodesic distance 
transformation method \cite{Soi91}, the thin plate model \cite{Duc76,Fra82} and the kriging method \cite{Ste99}.

\medskip
\noindent 
As for these latter methods, although there is a well-developed mathematical
theory on the existence and uniqueness of weak solutions of variational models \cite{ABCM01a,BCN02,GMS79}, 
and of the viscosity solution \cite{Aro67,Jen93} of
the PDE based interpolation model used in \cite{CMS98}, 
the quantitative effectiveness of such methods is mostly assessed on the basis of numerical experiments.

\medskip
\noindent 
The new approximation and interpolation theory introduced in \cite{ZCO16} is based, on the other hand, on 
the theory of compensated convex transforms \cite{Zha08, ZOC15a, ZOC15b, ZCO15} and can be applied to general bounded
real-valued functions sampled from either a compact set $K \subset \mathbb{R}^n$ or the complement 
$K = \mathbb{R}^n \setminus \Omega$ of a bounded open set $\Omega$. 
The methods presented in \cite{ZCO16} centre on the so-called average approximation that is recalled 
in Definition \ref{Def.LwrUpAv} below.  
Importantly, \cite{ZCO16} establishes error estimates for the approximation of bounded uniformly continuous functions,
or Lipschitz functions, and of $C^{1,1}$-functions, and proves rigorously that the approximation methods
are stable with respect to the Hausdorff distance between samples.

\medskip
\noindent 
Here we apply the average approximation method developed in \cite{ZCO16} to three important problems:
level set and scattered data approximation and interpolation, 
for which the sample set $K \subset \mathbb{R}^n$ is compact, and the inpainting problem in image processing,  
where the aim is to reconstruct an image in a damaged region based on 
the image values in the undamaged part and  the sample set $K= \mathbb{R}^n \setminus \Omega$ is the complement 
of a bounded open set $\Omega$ representing the damaged area of the image. 
We will also present  a series of prototype examples of explicitly calculated approximations that build insight 
into the behaviour of the average approximation introduced in \cite{ZCO16}, as well as a selection of  
illustrative numerical experiments.

\medskip
\noindent 
Before outlining the rest of the paper, we first recall the definitions of compensated convex transforms \cite{Zha08} 
and average approximation \cite{ZCO16}.  
Suppose $f:\mathbb{R}^n\to \mathbb{R}$ is bounded.
The quadratic lower and upper compensated convex transform \cite{Zha08}
(lower and upper transforms for short)  are defined for each $\lambda>0$ by
\begin{equation}\label{Eq.Def.UpLwTr}
\begin{split}
	&C^l_\lambda(f)(x)=\co[\lambda|\cdot|^2+f](x)-\lambda|x|^2,\\
	\text{resp.}\quad&
	C^u_\lambda(f)(x)=\lambda|x|^2-\co[\lambda|\cdot|^2-f](x),\qquad x\in\mathbb{R}^n,
\end{split}
\end{equation}
where $|x|$ is the standard Euclidean norm of $x\in \mathbb{R}^n$ and $\co[g]$ denotes the convex
envelope \cite{H-UL01,Roc70} of a function $g:\mathbb{R}^n\to \mathbb{R}$ that is bounded below.

\medskip
\noindent Let $K\subset \R^n$  be a non-empty closed set. Given a function $f:\mathbb{R}^n\to\mathbb{R}$, 
we denote by $f_K:\R^n\supset K \to \mathbb{R}$
the restriction of $f$ to $K$, which can be thought of as a sampling of the original function $f$, 
which we would like to approximate, on the convex hull of the set $K$.

 \medskip
\noindent 
Suppose that for some constant $A_0>0$, $|f_K(x)|\leq A_0$ for  all $x\in K$.
Then given $M>0$, we define two bounded functions that extend $f_K$ to $\R^n\setminus K$, namely
\begin{equation}\label{Eq.Def.ExtFnct}
\begin{array}{lll}
	\displaystyle f^{-M}_K(x) & \displaystyle = f(x)\chi_K(x)-M\chi_{\R^n\setminus K}
				& \displaystyle =
		\left\{\begin{array}{ll}
				f_K(x),	& x\in K,\\[1.5ex]
				-M,	& x\in\mathbb{R}^n\setminus K\,;
			\end{array}\right.\\[2.5ex]
	\displaystyle f^{M}_K(x) & \displaystyle = f(x)\chi_K(x)+M\chi_{\R^n\setminus K}
				& \displaystyle =
		\left\{\begin{array}{ll}
			f_K(x),	& x\in K,\\[1.5ex]
			M,	& x\in\mathbb{R}^n\setminus K\,,
		\end{array}\right.
\end{array}
\end{equation}
where  $\chi_G$ denotes the characteristic function of a set $G$.

\begin{defi}\label{Def.LwrUpAv}
	The {\bf average compensated convex approximation}  with scale $\lambda>0$ and module $M>0$
of the sampled function $f_K:K \to \R$ is defined by
\begin{equation}\label{Eq.Def.AvAprx}
	A^{M}_{\lambda}(f_K)(x)= \frac{1}{2}
		\left(C^l_\lambda(f^{M}_K)(x)+C^u_\lambda(f^{-M}_K)(x)\right),
		\quad x\in \mathbb{R}^n\,.
\end{equation}
\end{defi}

\bigskip
\noindent In addition,  we can also set $M=+\infty$ in place of \eqref{Eq.Def.ExtFnct} and consider the following functions, 
commonly used in convex analysis, 
\begin{equation}
	f^{-\infty}_K(x)=\left\{\begin{array}{l} f(x),\quad x\in K,\\
	-\infty,\quad x\in\mathbb{R}^n\setminus K;\end{array}\right.\qquad
	f^{+\infty}_K(x)=\left\{\begin{array}{l} f(x),\quad x\in K,\\
	+\infty,\quad x\in\mathbb{R}^n\setminus K.\end{array}\right.
\end{equation}
and define the corresponding average approximation 
approximation, 
\begin{equation}\label{Eq.Def.Ainfty}
	A_{\lambda}^{\infty}(f_K)(x) : = \frac{1}{2}
		\left(C^l_\lambda(f^{+\infty}_K)(x)+C^u_\lambda(f^{-\infty}_K)(x)\right),
		\quad x\in \mathbb{R}^n\,.
\end{equation}
By doing so, we can establish better approximation results than those obtained 
using $f^{-M}_K$ and $f^M_K$, but $A_{\lambda}^{\infty}(f_K)$ 
is not Hausdorff stable with respect to sample sets, in contrast to the basic average 
approximation $A^M_{\lambda}(f_K)$ (see \cite[Thm. 4.12]{ZCO16}).

\medskip
\noindent  
The plan of the rest of the paper is as follows. Section \ref{sect-not} introduces notation and recalls key definitions and  
results from our article \cite{ZCO16}, including error estimates for the average approximation $A^{M}_{\lambda}(f_K)$ of bounded and uniformly
continuous, Lipschitz, and $C^{1,1}$ functions. In Section \ref{Sec.LSAprx}, we consider level set interpolation and approximation, 
for which  $f$ is continuous and  $K$ consists of finitely many compact level sets. We give conditions so that $A^M_\lambda(f_{K})$ 
is an interpolation between level sets and also establish a maximum principle.
Section \ref{Sec.AprxSctD} treats the case of scattered data, when $K$  is finite.
In this case,  we show that when  $\lambda>0$ is sufficiently large and when $M>>\lambda$, $A^M_\lambda(f_{K})$ is a piecewise affine 
interpolation of $f_K$ in the convex hull of $K$. Moreover, if $K$ is regular
in the sense of the Delaunay triangulation, we show that   $A^M_\lambda(f_{K})$ agrees with the piecewise
interpolation given by the Delaunay method. In the irregular case that the Delaunay sphere $S_r$ contains more than $n+1$ points
in $\mathbb{R}^n$, $A^M_\lambda(f_{K})$ is the average of the maximum and minimum piecewise affine interpolation over 
the convex hull of $K\cap S_r$. Section \ref{Sec.Inp} presents error estimates for our average approximation in the context 
of the inpainting problem, and compares and contrasts these estimates with the error analysis in \cite{CK06}. 
We also give a simple one-dimensional example to illustrate the effect of the upper and lower compensated 
convex transforms $C^u_\lambda(f)$, $C^l_\lambda(f)$ and the average approximation $A^{M}_{\lambda}(f_K)$  on a jump function, 
to provide insight into how jump discontinuities behave under our approach. 

\medskip
\noindent 
Section \ref{Sec.ProtEx}  contains explicitly calculated prototype examples  in $\mathbb{R}^2$, including both 
examples where the sample set $K$ is finite, and also examples where $K$ is not finite.   We present  graphs of 
our calculated average approximation for two irregular Delaunay cells,  for $4$ and for $8$ points  on the unit circle. 
We also present  prototype examples of contour line approximations, as well as prototypes for inpainting of functions 
that  show that singularities such as ridges and jumps can be preserved subject to compensated convex approximations to 
the original function when the singular parts are close to each other.
Section \ref{Sec.NumEx} discusses several numerical experiments for level set and point clouds 
reconstructions of functions and images, 
for image inpainting, and for restoration of images with heavy salt \& pepper
noise. Though such experiments are carried out only on a proof-of-concept level, we briefly 
report on the comparison of our method with some state-of-art methods.
In Section \ref{Sec.Proofs}  we conclude the paper with  proofs of our main theorems stated in Sections \ref{Sec.LSAprx},  
\ref{Sec.AprxSctD} and \ref{Sec.Inp}.

%%%%%%%%%%%%%%%%%%%%%%%%%%%%%%%%%%%%%%%%%%%%%%%%%%%%%%%%%%%%%%%%%%%%%%%%%%%%%%%%%%%%%%%%%%%%%%%%%%%%
%%%%%%%%%%%%%%%%%%%%%%%%%%%%%%%%%%%%%%%%%%%%%%%%%%%%%%%%%%%%%%%%%%%%%%%%%%%%%%%%%%%%%%%%%%%%%%%%%%%%
%%%%%%%%%%%%%%%%%%%%%%%%%%%%%%%%%%%%%%%%%%%%%%%%%%%%%%%%%%%%%%%%%%%%%%%%%%%%%%%%%%%%%%%%%%%%%%%%%%%%

\setcounter{equation}{0}
\section{Notation and Preliminaries}
\label{sect-not}

\noindent Throughout the paper, we adopt the following notation and recall those results from \cite{ZCO16} that will be used here for our proofs,
 to make the development as self-contained as possible. For the necessary background in convex analysis, 
we refer to the monographs \cite{Roc70,H-UL01}.

\medskip
\noindent For a given set $E\subset \R^n$, with $\R^n$ a $n-$dimensional Euclidean space, $\bar{E}$,
$\partial{E}$, $\mathring{E}$, $E^c$ and $\co[E]$ stand for the closure, the boundary, the interior, the complement and the 
convex hull of $E$, i.e. the smallest convex set which contains $E$, respectively. For a convex set $E\subset \R^n$,
we define  the dimension of $E$, $\dim(E)$, as the dimension of the intersection of all affine manifolds that contain 
$E$, where by affine manifold we mean a translated subspace, i.e. a set $N$ of the form $N=x+S$
with $x\in\mathbb{R}^n$ and $S$ a subspace of $\mathbb{R}^n$. We then define $\dim(N)=\dim(S)$.
We use the term of convex body to denote a compact convex set with non-empty interior. 
The convex hull of a finite set of points is called a polytope and with the notation $\#(E)$ we denote the cardinality 
of the finite set $E$. If $E=\{x_1,\ldots,x_{k+1}\}$
and $\dim(E)=k$, then $\co[E]$ is called a $k-$dimensional simplex and the points 
$x_1,\ldots,x_{k+1}$ are called vertices. A zero-dimensional simplex is a point; a one-dimensional simplex is a line segment; 
a two-dimensional simplex is a triangle; a three-dimensional simplex is a tetrahedron. The condition that 
$\dim(E)=k$ is equivalent to require that the vectors $x_2-x_1,\ldots, x_{k+1}-x_1$ are linearly independent. 

\medskip
\noindent The open ball centered at $x\in\mathbb{R}^n$ and of radius $r>0$ is denoted by 
$B(x;r)=\{y\in\mathbb{R}^n: |y-x|<r\}$ where $|\cdot|$ stands for the
Euclidean norm in $\mathbb{R}^n$, thus $|x-y|$ is the distance between the points $x,\,y\in\R^n$. 
The diameter of the set $E\subset\mathbb{R}^n$, $\diam(E)$,
is then defined as $\diam(E)=\sup_{x,y\in E}|x-y|$.

\medskip
\noindent In this paper, we will assume, unless otherwise specified, that  $K\subset \R^n$ is either a compact set or the complement of a bounded open set, that is, 
$K=\Omega^c$ where $\Omega\subset \mathbb{R}^n$ is a bounded open set.
A function $g:\co[K]\subset \R^n\to \mathbb{R}$ is  said to be an 
{\it interpolation} of  $f_K$ if $g=f$ in $K$, while for $\lambda>0$, a family of
functions $g_\lambda:\co[K]\subset \R^n\to \mathbb{R}$ is said to {\it approximate} $f$ 
if $\displaystyle \lim_{\lambda\to+\infty} g_\lambda=f$
uniformly in $K$.

\medskip
\noindent The error estimates obtained in \cite{ZCO16} are expressed in terms of the modulus of continuity of the 
underlying function $f$ to be approximated and of the convex density radius of $K$. 
For the convenience of the reader, these definitions are recalled here.
The modulus of continuity of a bounded and uniformly continuous functions $f$ is defined as follows \cite{DL93,Hen62}.

\begin{defi}\label{Sec2.Def.MoC}
	Let $f:\R^n \to \R$ be a bounded and uniformly continuous function in $\R^n$. Then,
	\begin{equation}\label{Sec2.Eq.Def.MoC}
		\omega_f:t\in[0,\,\infty)\to
		\omega_{f}(t)=\sup\Big\{|f(x)-f(y)|:\,x,y\in\R^n\text{ and }|x-y|\leq t\Big\}
	\end{equation}
	is called the modulus of continuity of $f$.
\end{defi}

\noindent We also recall that the modulus of continuity of $f$ has the following properties \cite[page 19-21]{Hen62}.

\begin{prop}\label{Sec2.Pro.MoC}
Let $f:\R^n \to \R$ be a bounded and uniformly continuous function in $\R^n$.
Then the modulus of continuity $\omega_f$ of $f$ satisfies the following properties:
\begin{equation}\label{Sec2.Eq.Pro.MoC}
	\begin{array}{ll}
		$(i)$	&\omega_f(t)\to \omega_f(0)=0,\text{ as }t\to 0;	\\[1.5ex]
		$(ii)$	&\omega_f \text{ is non-negative and non-decreasing continuous function on }[0,\infty);	\\[1.5ex]
		$(iii)$	&\omega_f \text{ is subadditive: }\omega_f(t_1+t_2)\leq\omega_f(t_1)+\omega_f(t_2)
			\text{ for all }t_1,\,t_2\geq 0\,.
	\end{array}
\end{equation}
\end{prop}

\noindent Any function $\omega$ defined on $[0,\,\infty)$ and satisfying \eqref{Sec2.Eq.Pro.MoC}{\it (i), (ii), (iii)} is 
called {\it a
modulus of continuity}. A modulus of continuity
$\omega$ can be bounded from above by an affine function (see \cite[Lemma 6.1]{DL93}),
that is, there exist  constants $a>0$ and $b\geq 0$ such that
\begin{equation}\label{Sec2.Eq.Pro.MoC.0}
	\omega(t)\leq at+b\quad(\text{for all }t\geq 0).
\end{equation}
As a result, given $\omega_f$,  one can define the least concave majorant of $\omega_f$,
which we denote by $\omega$, which is also
a modulus of continuity with the property (see \cite{DL93})
\begin{equation}\label{Sec2.Eq.Pro.MoC.01}
	\frac{1}{2}\omega(t)\leq \omega_f(t) \leq \omega(t)\quad\text{for all }t\in[0,\,\infty)\,.
\end{equation}

\noindent The convex density radius of  a point $x\in\co[K]$ with 
respect to the set $K$ and the convex density radius of $K$ in $\co[K]$
are the geometrical 
quantities that describe the set $K$ with respect to its convex hull and are
such properties which enter
the error estimates for our approximation operators. 
We recall next their definition from \cite{ZCO16}.

\begin{defi}\label{Def.DnstRadius}		
Suppose $K\subset \mathbb{R}^n$ is a non-empty and closed
set, and denote by $\dist(x;\,K)$ the Euclidean distance of $x$ to $K$.
For $x\in \co[K]$, consider the balls $B(x;\,r)$ such that $x\in \co[\bar{B}(x;\,r)\cap K]$. 
The convex density radius of $x$ with respect to $K$ is defined as follows
\begin{equation}\label{Eq.Def.rc}
	r_c(x)=\inf\{r\geq 0 \text{ such that } x\in \co[\bar{B}(x;\,r)\cap K]\}\,,
\end{equation}
whereas the convex density radius of $K$ in $\co[K]$ is defined by
\begin{equation}\label{Eq.Def.rcSet}
	r_c(K)=\sup\{r_c(x),\; x\in \co[K]\}\,.
\end{equation}
\end{defi}

\noindent 
Here it is also  useful to introduce the following, more geometric quantities.
Let $Q \subset \R^n$  be a bounded set, and given $x\in Q$
and  $\nu\in\mathbb{R}^n$ with $|\nu|=1$, define the  quantity
\begin{equation*}
	d_\nu(x)=d^{+}_{\nu}(x)+d^{-}_{\nu}(x)\,, 
\end{equation*}	
where
\begin{equation*}
	d^{+}_{\nu}(x)=\sup\Big\{ t>0:\; x+s\nu\in Q \text{ for } 0\leq s\leq t\Big\}\,\text{ and }\, 
	d^{-}_{\nu}(x)=\sup\Big\{ t>0:\; x-s\nu\in Q \text{ for } 0\leq s\leq t\Big\}\,.
\end{equation*}
It is then easy to see that $d_{\nu}(x)$ is the length of the line segment with direction $\nu$ passing through $x$ 
and intersecting $\partial Q$ at two points on each side. We also define
\begin{equation}\label{Eq.Def.di}
	d(x)=\inf \Big\{d_{\nu}(x),\; \nu\in\mathbb{R}^n,\;|\nu|=1\Big\}\,\,
\end{equation}
and the thickness of the set $Q\subset \R^n$ as
\begin{equation}\label{Eq.Def.di}
	D_Q=\sup \Big\{d(x),\;x\in Q \Big\}\,.
\end{equation}

\begin{nota}
\begin{itemize}
\item[(a)]
	Given a non-empty bounded open set $Q= \Omega \subset \R^n$, 
	by comparing definition \eqref{Eq.Def.rc} of $r_c(x)$ and 
	\eqref{Eq.Def.di} of $d(x)$, it is straightforward to verify that 
	\begin{equation}\label{Eq.Ineq01}
		r_c(x)\leq d(x) 
	\end{equation}
	for $x\in \Omega$. 
	\item[(b)] If the interior 
	$\mathring{Q}=\varnothing$, such as in the case of a discrete set,
	then its thickness $D_Q$ is zero.
	\end{itemize}
\end{nota}

\noindent We recall next the error estimates for our average approximation operators developed in \cite{ZCO16}
and refer to \cite{ZCO16} for proofs and details. For the case of $K$ compact and $M=+\infty$,
we have the following. 

\begin{teo}\label{Thm.AprxUF.Cmpct}(See \cite[Theorem 3.6]{ZCO16})	
	Suppose $f:\mathbb{R}^n\to \mathbb{R}$ is bounded and uniformly continuous, 
	satisfying $|f(x)|\leq A_0$ for some constant $A_0>0$ and all 
	$x\in \R^n$, and let $K\subset \mathbb{R}^n$ be a non-empty compact set. 
	Denote by $\omega$ the least concave majorant of the modulus of continuity
	of $f$. Let $a\geq 0$, $b\geq 0$ be such that $\omega(t)\leq at+b$ for $t\geq 0$. 
	Then for all $\lambda >0$ and $x\in \co[K]$,
	\begin{equation}
		|A^\infty_\lambda(f_K)(x)-f(x)|\leq 
			\omega\left(r_c(x)+\frac{a}{\lambda}+\sqrt{\frac{2b}{\lambda}}\right)\,.
	\end{equation}
	where $r_c(x)\geq 0$ is the convex density radius of $x$ with respect to $K$.  
	If we further assume that  $f$ is a globally Lipschitz function with
				Lipschitz constant $L>0$, then for all $\lambda>0$ and   $x\in \co[K]$,
			\begin{equation}
				|A^\infty_\lambda(f_K)(x)-f(x)|\leq Lr_c(x)+\frac{L^2}{\lambda}\,.
			\end{equation}
\end{teo}

\noindent Section \ref{Sec.AprxSctD} will discuss an application of  Theorem \ref{Thm.AprxUF.Cmpct} to the case of scattered data approximation.
We will apply Theorem \ref{Thm.AprxUF.Cmpct} also to the case of salt-and-pepper noise removal, 
where $K$ is the compact set given by the part of the image which is noise free. Section \ref{Sec.NumEx}
contains a numerical experiment showing such an application.  

\medskip
\noindent A similar statement to Theorem \ref{Thm.AprxUF.Cmpct} is obtained with  $M$ finite in the case that 
$K=\Omega^c$, where $\Omega\subset \mathbb{R}^n$ is a non-empty bounded open set. 
In this case, clearly $\co[K]=\R^n$ and the error estimate of the average approximation 
$A^M_\lambda(f_{K})$ is  as follows.

\begin{teo}\label{Thm.AprxCnt}(See \cite[Theorem 3.7]{ZCO16}) 
	Suppose $f:\mathbb{R}^n\to \mathbb{R}$ is bounded and uniformly continuous, satisfying $|f(x)|\leq A_0$ for
	some constant $A_0>0$ and all $x\in \R^n$. 
	Let $\Omega\subset \mathbb{R}^n$ be a bounded open set, $d_\Omega$ the diameter of $\Omega$ and $K=\Omega^c$.
	Denote by $\omega$ the least concave majorant of the modulus of continuity of $f$ and 
	let $a\geq 0$, $b\geq 0$ be such that $\omega(t)\leq at+b$ for $t\geq 0$. 
	Then for $\lambda >0$, $M> A_0 + \lambda d_{\Omega}^2$, and all $x\in \mathbb{R}^n$, 
		\begin{equation}\label{Eq.Thm.AprxCnt.1} 
			|A^M_\lambda(f_{K})(x)-f(x)|\leq \omega\left(r_c(x)+\frac{a}{\lambda}+
			\sqrt{\frac{2b}{\lambda}}\right)\,,
		\end{equation} 
		where $r_c(x)\geq 0$ is the convex density radius of $x$ with respect to $K$.
	If we further assume that $f$ is a globally Lipschitz function with Lipschitz constant $L>0$,
	 then for $\lambda>0$,  $M>A_0+\lambda d_\Omega^2$ and all $x\in \mathbb{R}^n$, we have
	\begin{equation}
		|A^M_\lambda(f_{K})(x)-f(x)|\leq Lr_c(x)+\frac{L^2}{\lambda}\,.
	\end{equation}
\end{teo}

\noindent Under an additional restriction on $f$ and on $K$, it is possible to extend the results of Theorem \ref{Thm.AprxCnt}
 to the case when $K$ is a compact set and thus to obtain  error estimates independent of $M$. 
More precisely, the following result refers to the case where we are given the values of the function $f$
on the union of a compact set and the complement of a bounded open set. This extension allows
the application of Theorem \ref{Thm.AprxCnt} to the problem of inpainting, for instance.

\begin{coro}\label{Cor.AprxCnt}(See \cite[Corollary 3.9]{ZCO16}) 
	Suppose $f:\mathbb{R}^n\to \mathbb{R}$ is bounded and uniformly continuous satisfying $|f(x)|\leq A_0$ 
	for some constant $A_0>0$ and all $x\in\mathbb{R}^n$. Assume that
	$f(x)=c_0$ for $|x|\geq r>0$, where $c_0\in \mathbb{R}$ and $r>0$ are constants. 
	Let $K\subset \mathbb{R}^n$ be a non-empty compact set satisfying $K\subset \bar{B}(0;\, r)$.  
	For $R>r$, define $K_R:= K\cup B^c(0;\,R)$. Denote by $\omega$  the least concave majorant of the modulus 
	of continuity of $f$.
		Let $a\geq 0$, $b\geq 0$ be such that $\omega(t)\leq at+b$ for $t\geq 0$.
		Then for all $\lambda>0$, $M>A_0+\lambda (R+r)^2$
		and all $x\in \co[K]$, 
		\begin{equation} 
			|A^M_\lambda(f_{K_R})(x)-f(x)|\leq 
			\omega\left(r_c(x)+\frac{a}{\lambda}+\sqrt{\frac{2b}{\lambda}}\right)\,,
		\end{equation}
	where $r_c(x)\geq 0$ is the convex density radius of $x$ with respect to $K$.
	If we further assume that  $f$ is a globally Lipschitz function with Lipschitz constant $L>0$,
		then for $\lambda>0$,  $M>A_0+\lambda (R+r)^2$ and all $x\in \co[K]$, we have
		\begin{equation}
			|A^M_\lambda(f_{K_R})(x)-f(x)|\leq Lr_c(x)+\frac{L^2}{\lambda}\,.
		\end{equation}
	If we further assume that  $f$ is a $C^{1,1}$ function such that
		$|Df(x)-Df(y)|\leq L|x-y|$ for all $x,\,y\in \mathbb{R}^n$ and  $L>0$ is a constant,
		then for  $\lambda>L$,
		$M>A_0+\lambda (R+r)^2$ and all $x\in \co[K]$, we have
		\begin{equation}
			|A^M_\lambda(f_{K_R})(x)-f(x)|\leq
				\frac{L}{4}\left(\frac{\lambda+L/2}{\lambda-L/2}+1\right)r_c^2(x)\,.
		\end{equation}
		Furthermore, in  case $(iii)$,  $A^M_\lambda(f_{K_R})$ is an interpolation of $f_K$ in $\R^n$.
\end{coro}

\noindent The conditions of Corollary \ref{Cor.AprxCnt} can be realized, for instance, in the case 
we can define $f$ to be zero outside a large ball containing $K$. 

\medskip
\noindent Theorem \ref{Thm.AprxCnt} and Corollary \ref{Cor.AprxCnt} will be applied to 
the case of $(i)$ surface reconstructions from a finitely many level sets representation
and $(ii)$ inpainting of damaged images, where $\Omega$ is the domain to be inpainted and $K=\Omega^c$. 
We will discuss such applications in Section \ref{Sec.LSAprx} and Section \ref{Sec.Inp}, respectively,
whereas Section \ref{Sec.NumEx} contains some numerical experiments of
both applications.

We conclude this section by giving the following property which will be useful 
in Section \ref{Sec.AprxSctD} that deals with scattered data approximations.

\begin{prop}\label{Prp.Rest} (The restriction property) 
Let $m\geq 1$, $n\geq 1$. Suppose $f:\mathbb{R}^n\to \mathbb{R}$ is bounded, satisfying
$|f(x)|\leq M$ for some $M>0$  and for all $x\in \mathbb{R}^n$. 
Let $g^{\pm M}:\mathbb{R}^n\times \mathbb{R}^m\to \mathbb{R}$ be defined, respectively, as follows
\begin{equation*}
\begin{split}
	& g^{M}(x,y)=\left\{\begin{array}{ll}
			\displaystyle	f(x), & \displaystyle x\in\mathbb{R}^n,\; y=0\in\mathbb{R}^m,\\[1.5ex]
			\displaystyle	M,    & \displaystyle x\in \mathbb{R}^n,\; y\in \mathbb{R}^m,\; y\neq 0;
			\end{array}
			\right.	
	\\[1.5ex]
	& g^{-M}(x,y)=\left\{\begin{array}{ll}
			\displaystyle	f(x), & \displaystyle x\in\mathbb{R}^n,\; y=0\in\mathbb{R}^m,\\[1.5ex]
			\displaystyle	-M,   & \displaystyle x\in \mathbb{R}^n,\; y\in \mathbb{R}^m,\; y\neq 0\,.
			\end{array}
			\right.
\end{split}
\end{equation*}
Then
\begin{equation*}
	C^l_\lambda(g^{M})(x,0)=C^l_\lambda(f)(x)\quad\text{and}\quad 
	C^u_\lambda(g^{-M})(x,0)=C^u_\lambda(f)(x)\quad (\text{for }x\in\mathbb{R}^n)\,.
\end{equation*}
\end{prop}

In the case the sampled set $K$ is compact, the restriction property and Corollary \ref{Cor.AprxCnt} imply that 
if $K$ is contained in a $k$-dimensional plane $E\subset \mathbb{R}^n$, we can then calculate 
the average approximation operator $A^M_\lambda(f_K(x))$ for $x\in \co[K]\subset E$ 
by restricting our calculations in $E$.

%%%%%%%%%%%%%%%%%%%%%%%%%%%%%%%%%%%%%%%%%%%%%%%%%%%%%%%%%%%%%%%%%%%%%%%%%%%%%%%%%%%%%%%%%%%%%%%%%%%%
%%%%%%%%%%%%%%%%%%%%%%%%%%%%%%%%%%%%%%%%%%%%%%%%%%%%%%%%%%%%%%%%%%%%%%%%%%%%%%%%%%%%%%%%%%%%%%%%%%%%
%%%%%%%%%%%%%%%%%%%%%%%%%%%%%%%%%%%%%%%%%%%%%%%%%%%%%%%%%%%%%%%%%%%%%%%%%%%%%%%%%%%%%%%%%%%%%%%%%%%%

\setcounter{equation}{0}
\section{Level Set Approximations}\label{Sec.LSAprx}

We consider the case where the sampled set is given by the union of finitely many compact level sets, that is,
we know the values of a continuous function $f$ only on finitely many compact contour lines, and we want to study 
the structure of $A^M_\lambda(f_K)$.
We will establish a result which gives a natural bound on the value of $A^M_\lambda(f_K)$, ensuring that, for 
$\lambda>0$ sufficiently large,  the value of 
$A^M_\lambda(f_K)$ at points between level sets is between the values of the corresponding level sets, and present an error estimate for $A^M_\lambda(f_K)$.

\medskip
\noindent Let $f:\mathbb{R}^n\to \mathbb{R}$ be a continuous function and $a\in \mathbb{R}$.
Denote by $\Gamma_a=\{x\in \mathbb{R}^n,\; f(x)=a\}$ the level set of $f$ of level $a$
and by $V_a:=\{x\in \mathbb{R}^n,\; f(x)\leq a\}$ the sublevel set of $f$ of level $a$. 

\medskip
\noindent We then have the following result.

\begin{teo} \label{Thm.LSaprx}
Suppose $f:\mathbb{R}^n\to \mathbb{R}$ is continuous and that for $a_0<a_1<\cdots<a_m$, $m\in \mathbb{N}$,
the level sets $\Gamma_{a_i}=\{x\in \mathbb{R}^n,\; f(x)=a_i\}$ are compact
for $i=0,1,\ldots,m$. Denote by
\begin{equation*}
	\delta_0=\min\Big\{\dist(\Gamma_{a_i},\, \Gamma_{a_j}),\; 0\leq i,\,j\leq m,\; i\neq j\Big\}>0\,,
\end{equation*}
the minimum Euclidean distance between two different level sets.
Define $K=\cup_{i=0}^m\Gamma_{a_i}$ and denote by $d_K$ the diameter of $K$. 
If $\lambda>(a_m-a_0)/\delta_0^2$  and $M>\lambda d_K^2+\max_K|f|$, then 
\begin{itemize}
	\item[(i)] $A^M_\lambda(f_K)$ is an interpolation of $f$ from $K$ to $\co[K]$, that is, for $x_0\in \Gamma_{a_i}$, $i=0,1,\ldots,m$,
		\begin{equation}
		A^M_\lambda(f_K)(x_0)=a_i.\end{equation}
			\item[(ii)]  
		For each $x_0$ satisfying $a_i \leq f(x_0) \leq a_{i+1}$ for some $0\leq i\leq m-1$, 
			\begin{equation} 
			\label{ineq}
			a_i\leq A^M_\lambda(f_K)(x_0)\leq a_{i+1}\,.
		\end{equation}
	\item[(iii)] $A^M_\lambda(f_K)(x_0)=a_0$ for $x_0\in V_{a_0}$.
\end{itemize}
\end{teo}

\begin{nota}
\begin{itemize}
	\item[(a)] A sufficient condition for the level set $\Gamma_a$ to be compact is that 
		$f$ is continuous and either $\lim_{|x|\to\infty}f(x)=+\infty$ or $\lim_{|x|\to\infty}f(x)=-\infty$.
	\item[(b)] It might happen that there is an open subset of $\{x\in \mathbb{R}^n,\; a_i\leq f(x)\leq a_{i+1}\}$ 
			on which $A^M_\lambda(f_K)(x)=a_i$ or $A^M_\lambda(f_K)(x)=a_{i+1}$.
		Therefore Theorem \ref{Thm.LSaprx} gives a weak maximum principle.
	\item[(c)] In $\mathbb{R}^2$, it is not difficult to see that if two neighbouring level sets are parallel lines, 
		then our interpolation gives a plane passing through these two lines.
		However, if the function under consideration is not continuous, different level-sets can `intersect' each other.
		In general, it is not clear what the natural level-set approximations for functions with jump discontinuity
		will be like. In Section \ref{Sec.ProtEx} we will present a prototype example of two level lines which 
		are not parallel 
		to each other and work out an analytical expression of the interpolation operator $A^M_\lambda(f_K)$ for such a case.
\end{itemize}
\end{nota}

\medskip
\noindent We next give  an error estimate for our level set average approximation $A^M_\lambda(f_K)$, which is obtained by applying 
Corollary \ref{Cor.AprxCnt} \cite[Corollary 3.9]{ZCO16}. 
We first introduce some further definitions that are needed for the application of this result.
Under the assumptions of Theorem \ref{Thm.LSaprx}, for $i=0,1,\dots,m-1$, define the  open set
\begin{equation}\label{Eq.Def.omi}
	\Omega_i=\{x\in\mathbb{R}^n,\, a_i< f(x)<a_{i+1}\}\,,
\end{equation}
and then for $x\in\Omega_i$, define $d_i(x)$ using \eqref{Eq.Def.di} with $Q=\Omega_i$.  
	Suppose that $V_{a_m}$ is compact,  let  $R>0$  be such 
that $V_{a_m}\subset B(0;\,R)$, and set $V^m_R=V_{a_m}\cup B^c(0;\,R)$. Then define the 
auxiliary function
\begin{equation*}
	\tilde{f}_{V^m_R}(x)=\left\{\begin{array}{ll} 
				\displaystyle f(x),  & x\in V_{a_m},\\[1.5ex]
				\displaystyle a_m+1, & x\in B^c(0;\,R)\,.
			\end{array}\right.
\end{equation*}

\noindent We consider the following two cases.
\begin{itemize}
	\item[$(i)$] If $f$ is continuous, $\tilde{f}_{V^m_R}$ is bounded and uniformly continuous in $V^m_R$. 
		Therefore, by the Tietze extension theorem \cite[pag. 149]{Dug70}, 
		$\tilde{f}_{V^m_R}$ can be extended to $\mathbb{R}^n$ 
		as a bounded uniformly continuous function. We denote this extension by $\tilde{f}$ 
		and by $\tilde{A}_0>0$ an upper bound of $|\tilde{f}|$. 
		Clearly, $\tilde{f}(x)=f(x)$ for $x\in V_{a_m}$. 
		Furthermore, we denote 
		by $\tilde{\omega}(t)$ the least concave majorant of the modulus of continuity of $\tilde{f}$,
		which is itself a modulus of continuity, thus satisfies the properties \eqref{Sec2.Eq.Pro.MoC}, and in particular,
	 can be bounded from above by an affine function, that is, there exist some
		constants $\tilde a\geq 0$ and $\tilde b\geq 0$ such that 
		$\tilde \omega(t)\leq \tilde at+\tilde b$ for all $t\geq 0$.
	\item[$(ii)$] If $f$ is Lipschitz continuous with Lipschitz modulus $L>0$, then  
		$\tilde{f}_{V^m_R}$ is bounded and Lipschitz continuous in $V^m_R$ with 
		 a possibly different Lipschitz modulus $\tilde{L}$ such that 
		\begin{equation}\label{Eq.TilLip}	
			\tilde{L}\leq\max\Big\{L,\, \max_{V_{a_m}}|f|+|a_m+1|\Big\}\,.
		\end{equation}
		By Kirszbraun's theorem \cite[pag. 202]{Fed69}, $\tilde{f}_{V^m_R}$ can then be 
		extended to $\mathbb{R}^n$ 
		as a bounded Lipschitz function. Again we denote this extension by $\tilde{f}$ 
		and assume that $|\tilde{f}(x)|\leq \tilde{A}_0$ for all $x\in\mathbb{R}^n$.
\end{itemize}

\medskip
\noindent With the  notation above, we have the following error estimates for $A^M_\lambda(f_K)$. 

\begin{prop} \label{Prop.AprxLS}
	Suppose $f:\mathbb{R}^n\to \mathbb{R}$ is  continuous and that for $a_0<a_1<\cdots <a_m$, the sublevel sets 
	$V_{a_0}\subset V_{a_1}\subset\cdots\subset V_{a_m}$ are non-empty and compact. 
	Let  $\Gamma_{a_i}$  be  the level set of $f$ of level $a_i$, $K=\cup_{i=0}^m\Gamma_{a_i}$, and $d_i(x)$, $\Omega_i$  
	be defined by \eqref{Eq.Def.di}, \eqref{Eq.Def.omi}, respectively, for $i=0,1,\dots,m-1$. 
	Denote by $\tilde f$ the function defined in $(i)$ above, and by $\tilde{A}_0$ an upper bound 
	of $|\tilde f|$. If $\lambda> a_m-a_0 + 1$ and $M> \tilde{A}_0+\lambda (2R+1)^2$, then 
	for all $x\in\Omega_i$, $i=0,\,\ldots,m-1$, we have 
	\begin{equation}\label{Eq.Prop.AprxLS.1} 
		|A^M_\lambda(f_{K})(x)-f(x)|\leq \tilde{\omega}\left(d_i(x)+\frac{\tilde{a}}{\lambda}+
		\sqrt{\frac{2\tilde{b}}{\lambda}}\right)\,,
	\end{equation}
	where $\tilde{\omega}$ is the least concave majorant of the modulus of continuity of $\tilde{f}$.
	If we further assume that $f$ is a globally Lipschitz function of Lipschitz constant $L>0$, 
	$\lambda>a_m-a_0 + 1$ and $M>\tilde{A}_0+\lambda (2R+1)^2$, then for all $x\in\Omega_i$, $i=0,\,\ldots,m-1$, we have 
	\begin{equation} 
		|A^M_\lambda(f_{K})(x)-f(x)|\leq \tilde{L} d_i(x)+\frac{\tilde{L}^2}{\lambda}\,,
	\end{equation}
	where  $\tilde{L}$ is defined by \eqref{Eq.TilLip}.
\end{prop}

%%%%%%%%%%%%%%%%%%%%%%%%%%%%%%%%%%%%%%%%%%%%%%%%%%%%%%%%%%%%%%%%%%%%%%%%%%%%%%%%%%%%%%%%%%%
%%%%%%%%%%%%%%%%%%%%%%%%%%%%%%%%%%%%%%%%%%%%%%%%%%%%%%%%%%%%%%%%%%%%%%%%%%%%%%%%%%%%%%%%%%%
%%%%%%%%%%%%%%%%%%%%%%%%%%%%%%%%%%%%%%%%%%%%%%%%%%%%%%%%%%%%%%%%%%%%%%%%%%%%%%%%%%%%%%%%%%%

\setcounter{equation}{0}
\section{Scattered Data Approximations}\label{Sec.AprxSctD}

We now turn  our attention to the so-called case of `scattered data' approximation \cite{Wen05} 
corresponding to a discrete sampled set $K$.  
Since for any function
$f:\mathbb{R}^n\to\mathbb{R}$, the restriction $f_K$ of $f$ to a finite set $K$ is always
a Lipschitz function, the following result provides a sufficient condition for our
upper and lower transforms to be interpolations in this case.

\begin{teo}\label{Thm.ULAprx} %theorem 5.3
Suppose $K=\{x_1,x_2,\dots,x_m\} \subset \mathbb{R}^n$ is a finite set
with distinct points and assume $f:K\subset \mathbb{R}^n\to\mathbb{R}$ is a function. 
Assume $-M<f(x_j)<M$ for $j=1,\ldots,m$ and let $L>0$ be the Lipschitz
constant of $f:K\subset \mathbb{R}^n\to \mathbb{R}$. Define
$\alpha=\min\{|x_i-x_j|,\; x_i,\, x_j\in K,\, i\neq j\}>0$. Then for
$\lambda\geq L/\alpha$,
\begin{equation*}
	C^u_\lambda(f^{-M}_K)(x_j)=f(x_j)\quad\text{and}\quad 
	C^l_\lambda(f^M_K)(x_j)=f(x_j)\quad \text{for }x_j\in K\,.
\end{equation*}
\end{teo}

\noindent Let $K\subset \mathbb{R}^n$ be a finite set. Without loss of generality, from now on, 
we assume that
$\dim(\co[K])=n$, that is, that $\co[K]\subset \mathbb{R}^n$ is a convex body. In the case
$\dim(\co[K])=k<n$, we can simply translate $K$ so that $0\in K$, and let $E_k=\Span(\co[K])$ where 
$\Span(\co[K])$ is the 
$k$-dimensional subspace spanned by $\co[K]$. In this case, $E_k\subset\mathbb{R}^n$ is a supporting plane of 
$\co[K]$ and we only need to work in $E_k$ given that in our interpolation problem 
we are only interested in values of our approximation in $\co[K]$.
We can therefore reduce our approximation/interpolation problem to $E_k$ 
by applying Proposition \ref{Prp.Rest}.

\medskip
\noindent In order to describe our approximation/interpolation results,  we first need  to introduce
notions related to the Voronoi diagram and Delaunay triangulation for a 
finite set $K$ \cite{BKOO93,OBSC00,Ede87}.

\medskip
\noindent Let $K=\{x_1,\ldots,x_m\}$ be a finite set of distinct points of $\R^n$,
and denote $m=\#(K)$. We define $\mathcal{V}(K)$, the Voronoi diagram of $K$,
to be the partition of $\R^n$ into $m$ cells, one for each point of $K$, with the property that
a point $x\in \R^n$ belongs to the cell corresponding to the point $x_i\in K$ if $|x-x_i|<|x-x_j|$ for each
$x_j\in K$ with $j\not=i$. We then denote by $M(K)$
the Voronoi edges of the Voronoi diagram $\mathcal{V}(K)$ of $K$,  meaning the set of the edges of $\mathcal{V}(K)$
where a point $y\in M(K)$ if 
there are at least two different points $x_i,\, x_j\in K$ such that $\dist(y,\,K)=|y-x_i|=|y-x_j|>0$. Then there are finitely many points   $y_1,\dots, y_l\in
M(K)$, called Voronoi vertices and whose set is denoted by $V(K)$, with the property 
that there are corresponding radii $r_1,\dots, r_l>0$, 
such that for each $y_i \in V(K)$, there are $m_i\geq n+1$ points $x^i_1,\dots,x_{m_i}^i\in K$ 
such that $\dist(y_i,\, K)= |y_i-x^i_j|=r_i$ so that the open ball $B(y_i;\,r_i)$ 
does not intersect $K$ and $\bar{B}(y_i;\,r_i)\cap K=\{x_1^i,\ldots,x^i_{m_i}\}$. If we write  $K_i=\{x^i_1,\ldots, x^i_{m_i}\}$ for each $i \in \{1, \ldots, l\}$, 
we also have that $\dim(\co[K_i])=n$, $\cup_{j=1}^l\co[K_j]=\co[K]$, and if $i\neq j$, either $\dim(\co[K_i]\cap \co[K_j])<n$  or
$\co[K_i]\cap \co[K_j] = \emptyset$ \cite{OBSC00}.

\medskip
\noindent 
For each $i=1,\ldots, l$, $\co[K_i]$ is referred to as a Delaunay cell with generator $K_i$, centre $y_i$
and radius $r_i$ and the ball $B(y_i;\,r_i)$ is called the associated open ball of 
the Delaunay cell $\co[K_i]$. We have $K_i=K\cap \partial B(y_i;\,r_i)$ while $K\cap B(y_i;\,r_i)=\varnothing$. 
A Delaunay cell is then said regular if it is an $n$-dimensional simplex (so in particular,  a triangle if $n=2$ and a tetrahedron if $n=3$). If each Delaunay cell $\co[K_i]$ in $\co[K]$ is regular, the set $\{ \co[K_1], \co[K_2], \ldots, \co[K_l]\}$ is said to be the regular Delaunay
triangulation of $\co[K]$.

\medskip
\noindent In the following, we consider two different situations. 
\begin{itemize}
	\item[$(i)$] 
			Each Delaunay cell $\co[K_i]$ is an $n$-dimensional simplex, that is, $\co[K]$ has a regular Delaunay triangulation;
	\item[$(ii)$] For some or for all $K_i$'s, $\dim(\co[K_i])=\dim(\co[K])=n$ and 
		$\#(K_i)>n+1$, that is, the Delaunay cell is a convex polytope that is not an $n$-dimensional simplex.
\end{itemize}

\noindent We will show that if $(i)$ holds, that is, if we have a regular Delaunay triangulation of $\co[K]$, 
then our average approximation $A^M_\lambda(f_K)$ defines the usual piecewise affine interpolation
based on this Delaunay triangulation \cite[page. 191]{OBSC00} when $\lambda>0$ and $M>>\lambda$ are sufficiently large. 
If $(ii)$ occurs, our average approximation $A^M_\lambda(f_K)$ will be
the average of the minimum and maximum piecewise affine interpolations of $f_K$ in the cell. 

\begin{nota}
A remarkable difference between our average approximation $A^M_\lambda(f_K)$ and the usual 
design of piecewise affine constructions is that we do not need to know or compute 
the Delaunay cells in advance. Our method simply directly generates the  piecewise affine function.
\end{nota}

\noindent Before we state our first structural theorem on the effect of the upper, 
lower and average approximations over a regular cell, we need the following lemma.

\begin{lema}\label{Lem.BndAf}	%Lemma 5.4
Let $B(x^\ast;\,r)\subset  \mathbb{R}^n$ be the open ball centred at $x^\ast$
with radius $r>0$ and $S=\{x_1,x_2,\dots,x_{m}\}\subset \partial B(x^\ast;\,r)$ be a finite set with 
distinct points and with $\#(S)=m\geq n+1$. Assume $\co[S]\subset \bar B(x^\ast;\,r)$ to be the convex hull of $S$ satisfying 
$\dim(\co[S])=n$. Suppose $f_S:S\to\mathbb{R}$ is a real-valued function with
Lipschitz constant $L>0$. If there is an affine function $\ell_s:\mathbb{R}^n\to \mathbb{R}$ 
such that $\ell_s(x_i)=f_S(x_i)$ for all $x_i\in S$, then there is a constant $C_s>0$ such that the gradient
of $\ell$ satisfies $|D\ell_s(x)|\leq C_sL$.
\end{lema}

\begin{nota}
	In Lemma \ref{Lem.BndAf}, if $m=n+1$, then $\co[S]$ is a $n$-dimensional simplex and there is an
	affine function $\ell_s$ such that $\ell_s(x)=f_S(x)$ for $x\in S$. However if $m>n+1$, 
	in general one can not find an affine function satisfying $\ell_s(x)=f_S(x)$ for $x\in S$.
	We will deal with such a case together with a more general one in Lemma \ref{Lem.AfAprx} and in 
	Theorem \ref{Thm.AfAprx}.
\end{nota}

\noindent We now calculate the transforms $C^u_\lambda(f^{-M}_K)$,
$C^l_\lambda(f^{M}_K)$ and $A^M_\lambda(f_K)$ in a regular Delaunay cell
$\co[S]$ satisfying $m=\#(S)=n+1$ and $\dim(\co[S])=n$. For each regular cell $\co[S]$,  define
\begin{equation*}
	\sigma_s=\min\Big\{|x_j-x_s|-r_s,\; x_j\in K\setminus S\Big\}>0
\end{equation*}
where $x_s, r_s$ are the centre and radius respectively of the associated Delaunay ball $B(x_s; r_s)$ of $\co[S]$, and let $C_s$ be the constant given by Lemma \ref{Lem.BndAf} for the affine function $\ell_s$ 
associated with $\{(x,f_S(x)), \, x\in S\}$. We then have the following result. 

\begin{teo}\label{Thm.AprxRC} 
Let $K=\{x_i\}_{i=1}^m\subset\mathbb{R}^n$ be a finite set with distinct points and
let $f_K:K\to\mathbb{R}$ be a function with Lipschitz constant $L>0$ and bound $A_0>0$, that is,
$|f_K(x)|\leq A_0$ for $x\in K$.
Suppose $S=\{x_1,x_2,\dots, x_{l+1}\}\subset K$ satisfies that $\co[S]$ is  a regular  Delaunay cell
with associated Delaunay ball $B(x_s;\,r_s)$. Let $\ell_s:\mathbb{R}^n\to \mathbb{R}$
be the affine function given by Lemma \ref{Lem.BndAf} for $S$ and $f_K$ restricted on $S$.
Then, for every $x\in \co[S]$,
\begin{equation}\label{Thm.AprxRC.1}
\begin{array}{c}
	\displaystyle C^u_\lambda(f^{-M}_K)(x)=\lambda|x-x_s|^2-\lambda r_s^2+\ell_s(x),\qquad
	\displaystyle C^l_\lambda(f^{M}_K)(x)=\lambda r_s^2-\lambda|x-x_s|^2+\ell_s(x),\\[1.5ex]
	\displaystyle A^M_\lambda(f_K)(x)=\frac{C^u_\lambda(f^{-M}_K)(x)+C^l_\lambda(f^{M}_K)(x)}{2}=\ell_s(x)\,,
	\end{array}
\end{equation}
whenever
\begin{equation}\label{Thm.AprxRC.2}
	\lambda>\frac{2A_0}{\sigma_s(2r_s+\sigma_s)}+\frac{C_sL}{\sigma_s}
\end{equation}
and
\begin{equation}\label{Thm.AprxRC.3}
	M>\lambda r_s^2+C_sLr_s+A_0+\frac{C_s^2L^2}{4\lambda}.
\end{equation}
\end{teo}

\medskip
\begin{nota}\label{Nota.AprxRCS} 
If we replace our functions $f_K^{-M}$ and $f_K^{M}$ by $f_K^{-\infty}$ and 
$f_K^{\infty}$, respectively, defined by
\begin{equation*}
	f_K^{-\infty}(x)=\left\{\begin{array}{ll} 
			\displaystyle f_K(x),	&\displaystyle 	\text{if } x\in\in K\,,\\[1.5ex]
			\displaystyle -\infty,	&\displaystyle 	\text{if } x\in \mathbb{R}^n\setminus K
		\end{array}\right.\qquad\text{and}\qquad
	f_K^{\infty}(x)=\left\{\begin{array}{ll} 
			\displaystyle f_K(x),	&\displaystyle \text{if } x\in K\,,\\[1.5ex]
			\displaystyle +\infty,  &\displaystyle \text{if } x\in \mathbb{R}^n\setminus K\,,
	\end{array}\right.
\end{equation*}
then  Condition \eqref{Thm.AprxRC.2} alone is  sufficient to obtain \eqref{Thm.AprxRC.1}.
Although by setting $M=+\infty$ we have a mathematically simpler statement,  
the resulting approximations would not, however, meet the Hausdorff stability property (see \cite[Thm. 4.12]{ZCO16}  for a Hausdorff stability theorem for $A^M_{\lambda}(f_K)$).
\end{nota}

\noindent If we further assume that for the given finite set $K$ there is a regular Delaunay triangulation of $\co[K]$, 
which thus consists of $n$-dimensional simplices,
we can then easily give  global explicit descriptions of 
$C^u_\lambda(f^{-M}_K)$ and $C^l_\lambda(f^{M}_K)$, and hence of $A^M_\lambda(f_K)$ in each $n$-dimensional 
Delaunay simplex. This, however, requires $\lambda>0$ and $M>0$ to be sufficiently large.

\begin{coro}\label{Cor.AprxRCS} 
Let $K\subset \mathbb{R}^n$ be a finite set with distinct points such that it admits a regular Delaunay triangulation
$\mathcal{D}(K)$ of $\co[K]$ thus comprising of the $n$-dimensional simplices $\co[S_1], \ldots,\co[S_l]$
where $V(K)$ the set of vertices of the Voronoi diagram $\mathcal{V}(K)$ of $K$ with $\#(V(K))=l$. For each Delaunay cell $S_i$  for $i=1,\ldots,l$, 
consider its associated open ball 
$B(y_i;\,r_i)$ such that $B(y_i;\,r_i)\cap K=\varnothing$
and $K\cap\bar{B}(y_i;\,r_i)=S_i$ for $i=1,\ldots,l$. Define $\sigma_i=\min\{|x-y_i|-r_i,\; x\in K\setminus S_i\}$.

\medskip
\noindent Let $f_K:K\subset \R^n\to \mathbb{R}$ be a function with Lipschitz
constant $L>0$ satisfying, for some $A_0>0$, $|f_K(x)|\leq A_0$ for all $x\in K$.
Let $\ell_i$ be the affine function defined in Lemma \ref{Lem.BndAf} for $S_i$, such that $\ell_i(x)=f_K(x)$ for
$x\in S_i$ and $|D\ell_i(x)|\leq C_iL$ for some constant $C_i>0$, $i=1,\ldots,l$. 
Then in each simplex $\co[S_i]$, $i=1,\ldots,l$, and for every $x\in \co[S_i]$, we have
\begin{equation}\label{Eq.Cor.AprxRCS.1}
\begin{array}{c}
	\displaystyle C^u_\lambda(f^{-M}_K)(x)= \lambda|x-x_i|^2-\lambda r_i^2+\ell_i(x),\qquad
	\displaystyle C^l_\lambda(f^{M}_K)(x)=\lambda r_i^2-\lambda|x-x_i|^2+\ell_i(x)\,,\\[1.5ex]
	\displaystyle A^M_\lambda(f_K)(x)=\frac{C^u_\lambda(f^{-M}_K)(x)+C^l_\lambda(f^{M}_K)(x)}{2}=\ell_i(x)\,,
\end{array}
\end{equation}
whenever
\begin{equation}\label{Eq.Cor.AprxRCS.2}
	\lambda>\max_{1\leq i\leq m}\left(\frac{2A_0}{\sigma_i(2r_i+\sigma_i)}+\frac{C_iL}{\sigma_i}\right)
\end{equation}
and
\begin{equation*}
	M>\max_{1\leq i\leq m}\left(\lambda r_i^2+C_iLr_i+A_0+\frac{C_i^2L^2}{4\lambda}\right)\,.
\end{equation*}
\end{coro}

\begin{nota}\label{Nota.AprxRCSInf} 
A similar observation to Remark \ref{Nota.AprxRCS} for Theorem \ref{Thm.AprxRC} can be made for Corollary \ref{Cor.AprxRCS}. Under the assumptions of Corollary \ref{Cor.AprxRCS}, condition \eqref{Eq.Cor.AprxRCS.2}
is sufficient to ensure that  \eqref{Eq.Cor.AprxRCS.1} holds with $f_K^{-\infty}$,  $f_K^{\infty}$ and 
$A_{\lambda}^{\infty}(f_K)$, respectively, for $i=1,\ldots,l$ and for every $x\in \co[S_i]$. \end{nota}

\medskip
\noindent Let $S=\{x_1,\ldots,x_m\}\subset \mathbb{R}^n$. 
Next we study the structure of our upper, lower transforms and average approximations when the $n$-dimensional 
Delaunay cell $\co[S]$ is not a simplex, that is, $\#(S)=m>n+1$. In this case, we say that the 
$n$-dimensional Delaunay cell $\co[S]$ is not regular.
Without loss of generality we may assume that  there is an  open ball $B(0;\,r)$ centred at $0$ 
with radius $r>0$, such that $S\subset \partial B(0;\,r)$. 
Let $f_S:S\to \mathbb{R}$ be a given function, and write $f_S(x_i)=v_i$, $i=1,\ldots,m$.
Let $\Gamma_s=\{(x_i,v_i),\, i=1,\ldots,m\}$ be the graph of $f_S$ in $S\times \mathbb{R}$, 
we may assume that the convex envelope $\co[\Gamma_s]\subset \mathbb{R}^n\times \mathbb{R}$ 
of $\Gamma_s$  is an $n+1$-dimensional convex polytope, otherwise there will be a 
single affine function as in Lemma \ref{Lem.BndAf} 
satisfying $\ell_s(x_i)=v_i$ and we are back to the situation of Theorem \ref{Thm.AprxRC}.

\medskip
\noindent Let $D=\co[S]\subset \mathbb{R}^n$ and $\Gamma=\partial \co[\Gamma_s]$ be the boundary of
the convex polytope $\co[\Gamma_s]$. We have the following result. 

\begin{lema}\label{Lem.AfAprx}
Let $S$, $f_S$ and $\Gamma_s$ be as defined above. Then
\begin{itemize}
	\item[(i)] There are two continuous piecewise affine functions 
			$p_+(x)$ and $p_-(x)$ in $D=\co[S]$ defined by
		\begin{equation*} 
		\begin{split} 
			p_+(x)	& =\max\{v,\;(x,v)\in \co[\Gamma_s]\}\\[1.5ex]
				& =\max\left\{\sum^m_{i=1}\lambda_iv_i,\; x_i\in S,\;\lambda_i\geq 0,\; 
				i=1,\ldots,m,\; \sum^m_{i=1}\lambda_i=1,\; \sum^m_{i=1}\lambda_ix_i=x\right\}\,,\\[1.5ex]
			p_-(x)	& =\min\{v,\; (x,v)\in \co[\Gamma_s]\}\\[1.5ex]
				& =\min\left\{\sum^m_{i=1}\lambda_iv_i,\; x_i\in S,\;\lambda_i\geq 0,\; 
				i=1,\ldots,m,\; \sum^m_{i=1}\lambda_i=1,\; \sum^m_{i=1}\lambda_ix_i=x\right\}\,,
		\end{split}
		\end{equation*}
		where $p_+$ and $p_-$  are piecewise affine concave and convex functions in $D$ respectively;
	\item[(ii)] For every $x\in \mathring{D}$, the interior of $D$,  $p_-(x)< p_+(x)$.
	\item[(iii)] The convex polytope $D\subset\mathbb{R}^n$ has two decompositions
		$D=\cup^k_{i=1}D^+_i$ and $D=\cup^l_{j=1}D^-_j$ such that
		$D^+_k$ and $D^-_j$ are closed convex
	 $n$-dimensional polytopes, 
		$\mathring{D}^+_i\cap \mathring{D}^+_j=\varnothing$
		and $\mathring{D}^-_i\cap \mathring{D}^-_j=\varnothing$ for $1\leq i\neq j\leq l$. 
		On each $D^+_k$ (respectively, $D^-_j$),  $p_+(x)$ (respectively, $p_-(x)$) 
		is an affine function, that is, $p_+(x):=\ell_k^+(x)=a^+_k\cdot x+b_k^+$, $x\in D^+_k$ 
		(respectively,  $p_-(x):=\ell_j^-(x)=a_j^-\cdot x+b_j^-$, $x\in D_j^-$).  
		Furthermore, the affine function $\ell_k^+(x)$ (respectively, $\ell_j^-(x)$) 
		defined in $\mathbb{R}^n$ as above, satisfies $\ell_k^+(x)\geq p_+(x)$ 
		(respectively, $\ell_j^-(x)\leq p_-(x)$) for $x\in D$.
	\item[(iv)] Let $S^+_k\subset D^+_k$ be the set of all vertices of $D^+_k$ for 
		$k=1,\ldots,m$, then $S^+_k\subset S$,  
		and $\cup^m_{k=1}S^+_k=S$. On each $S^+_k$, $p_+(x)=f_S(x)$.
	\item[(v)] Let $S^-_j\subset D^-_j$ be the set of all vertices of $D^-_j$ 
		for $j=1,\ldots,l$, then $S^-_j\subset S$,
		and $\cup^l_{j=1}S^-_j=S$. On each $S^-_k$, $p_-(x)=f_S(x)$.
	\end{itemize}
\end{lema}

\begin{nota}
In Lemma \ref{Lem.AfAprx}, the piecewise affine functions $p_+$ and $p_-$ 
are replacements of $\ell_s$ in Theorem \ref{Thm.AprxRC}.
For the average approximation, the average $\frac{p_+ \, +\, p_-}{2}$ of the piecewise affine functions  $p_+$ and $p_-$
 gives the new interpolation formula in $D=\co[S]$, 
replacing the affine function $\ell_S$. This means that our interpolation $A^M_{\lambda}(f_K)$ might 
introduce extra nodes  in $\co[S]$ in a unique way, in the sense that $D$ is the union of $q$ $n$-dimensional 
convex polytopes $D_i^{av}$, $i \in \{1, \ldots, q \}$,  such that
$\frac{p_+\, +\, p_-}{2}$ is affine on each $D_i^{av}$ but not all vertices of $D_i^{av}$ are contained in $S$.
\end{nota}

\noindent The following is a generalisation of Theorem \ref{Thm.AprxRC}.

\begin{teo}\label{Thm.AfAprx}
Let $K=\{x_i\}_{i=1}^m\subset\mathbb{R}^n$ be a finite set with distinct points and let $f_K:K\to\mathbb{R}$ 
be a function with Lipschitz constant $L>0$ and bound $A_0>0$, that is,
$|f(x)|\leq A_0$ for $x\in K$.
Suppose $S=\{x_1,x_1,\ldots, x_{m}\}\subset K$ generates a 
Delaunay cell $ \co[S]$ satisfying $\dim(\co[S])=n$ and  
$\dim(\co[\Gamma_s])=n+1$, where $\Gamma_s=\{(x,\, f_K(x)),\; x\in S\}$ 
is the graph of $f_K$ restricted to $S$. Let
$B(y_s;\,r_s)$ be the associated open ball of the cell $ \co[S]$. 
Let $p_+: \co[S]\to \mathbb{R}$ be the piecewise affine concave 
function and $p_-: \co[S]\to \mathbb{R}$ be the piecewise affine 
convex function defined in Lemma \ref{Lem.AfAprx}, and let $ \co[S]=\cup^m_{k=1}D^+_k$ and $ \co[S]=\cup^l_{j=1}D^-_j$ be the 
decompositions of $ \co[S]$ given by Lemma \ref{Lem.AfAprx}.
Let
\begin{equation*}
	C^+_sL=\max_{1\leq k\leq m} 
	C^+_kL,\quad  C^-_sL=\max_{1\leq j\leq l}C^-_jL, \quad C_sL=\max\{C^+_sL,\, C^-_sL\}\,,
\end{equation*}
where $C^+_kL$ and $C^-_jL$ are the positive upper bounds given by Lemma \ref{Lem.BndAf} for 
$|Dp_+(x)|$ and $|Dp_-(x)|$, respectively, on $D^+_k$ and $D^-_j$.
Let $\sigma_s=\min\{|x-x_s|-r_s,\; x\in K\setminus S\}>0$. Then for every $x\in  \co[S]$,
\begin{equation}
\begin{array}{c}
	\displaystyle C^u_\lambda(f^{-M}_K)(x)=\lambda|x-x_s|^2-\lambda r_s^2+p_+(x)\,,\qquad
	\displaystyle C^l_\lambda(f^{M}_K)(x)=\lambda r_s^2-\lambda|x-x_s|^2+p_-(x)\,,\\[1.5ex]
	\displaystyle A^M_\lambda(f_K)(x)=\frac{p_+(x)+p_-(x)}{2}\,,
\end{array}
\end{equation}
whenever
\begin{equation}
	\lambda>\frac{2A_0}{\sigma_s(2r_s+\sigma_s)}+\frac{C_sL}{\sigma_s}
\end{equation}
and
\begin{equation}
	M>\lambda r_s^2+C_sLr_s+A_0+\frac{C_s^2L^2}{4\lambda}\,.
\end{equation}
\end{teo}

\begin{nota}
Under the assumptions of Lemma \ref{Lem.AfAprx} and Theorem \ref{Thm.AfAprx}, 
we see that $p_+(x)$ and $p_-(x)$ 
are the maximal and minimal piecewise affine interpolations over $ \co[S]$. 
It is well-known \cite{OBSC00} that in this irregular case, 
there still exist Delaunay triangulations of $\co[S]$ consisting of $n$-dimensional simplices, but
the triangulation is not unique. The average approximation 
\begin{equation*}
	A^M_\lambda(f_K)(x)=\frac{p_+(x)+p_-(x)}{2}
\end{equation*}
given by Theorem \ref{Thm.AfAprx} is exactly the average of the maximal and minimal 
interpolation in a Delaunay cell.
\end{nota}

%%%%%%%%%%%%%%%%%%%%%%%%%%%%%%%%%%%%%%%%%%%%%%%%%%%%%%%%%%%%%%%%%%%%%%%%%%%%%%%%%%%%%%%%%%%%%%%%%%%%
%%%%%%%%%%%%%%%%%%%%%%%%%%%%%%%%%%%%%%%%%%%%%%%%%%%%%%%%%%%%%%%%%%%%%%%%%%%%%%%%%%%%%%%%%%%%%%%%%%%%
%%%%%%%%%%%%%%%%%%%%%%%%%%%%%%%%%%%%%%%%%%%%%%%%%%%%%%%%%%%%%%%%%%%%%%%%%%%%%%%%%%%%%%%%%%%%%%%%%%%%

\setcounter{equation}{0}
\section{Inpainting revisited}\label{Sec.Inp}
Consider now  inpainting of damaged areas of an image. This is the problem where 
we are given an image that is damaged in some parts and we want to reconstruct the values in the damaged part
on the basis of the known values of the image. 
To specify the setting of the problem, let $\Lambda\subset\R^n$ be a convex compact set 
representing the domain of the image $f$ which, without loss 
of generality, we assume to be a grayscale image, and is thus represented by a function $f:\Lambda \subset \R^n\to \R$.
We assume that $f$ is bounded and uniformly continuous. 
See below, in Remark \ref{Rmk.ErEst} and  the comments on Example \ref{Ex.AprxJump}, for a discussion of this assumption in the case of an image.

\medskip
\noindent Denote  by $\Omega\subset \Lambda$ an open set representing the damaged areas of the image and let $K=\Lambda\setminus \Omega$.
We have then $\Omega\subset \co[K]$. 

\medskip\noindent On the basis of the values of $f$ in $K$, we reconstruct the values of 
$f$ in $\Omega$ by using the average approximation $A_{\lambda}^M(f_K)$. 
 In this section, we want
to assess the error of this approximation. 

\medskip
\noindent The next result, which follows from an application of Corollary \ref{Cor.AprxCnt},
is the main error estimate for our inpainting method. 

\begin{prop}\label{Pro.EstEr.Inp}
Let $\Lambda \subset\R^n$ be a convex compact set and $\Omega\subset \Lambda$ a non-empty open set. Assume $f:\Lambda\subset \R^n\to \R$
be bounded and uniformly continuous, such that for $A_0>0$ we have that  $|f(x)|\leq A_0$ for all $x\in K=\Lambda\setminus \Omega$.
Let $\tilde{f}$ be a bounded and uniformly continuous extension of $f$ to $\R^n$, derived by the Tietze extension theorem,
with $\tilde{f}(x)=c_0$ outside an open ball $B(0; r)$ with $r>0$ and such that $K\subset B(0; r)$. For $R>r$, 
define $K_R=K\cup B^c(0;R)$ and let $f_{K_R}(x)=f_K(x)$ for $x\in K$ and $f_{K_R}(x)=c_0$ for $x\in B^c(0; R)$. 
Denote by $\omega$  the least concave majorant of the modulus of continuity of $\tilde{f}$.
Let $a\geq 0$, $b\geq 0$ be such that $\omega(t)\leq at+b$ for $t\geq 0$.
Then for all $\lambda>0$, $M>A_0+\lambda (R+r)^2$ and all $x\in \co[K]$, we have
\begin{equation}\label{Eq.ErEstImpFcnt} 
	|A^M_\lambda(f_{K})(x)-\tilde{f}(x)|\leq 
	\omega\left(r_c(x)+\frac{a}{\lambda}+\sqrt{\frac{2b}{\lambda}}\right)\,,
\end{equation}
where $r_c(x)\geq 0$ is the convex density radius of $x$ with respect to $K$. 

\medskip
\noindent 
If we further assume that  $f$ is a globally Lipschitz function with Lipschitz constant $L>0$,
		then for $\lambda>0$,  $M>A_0+\lambda (R+r)^2$ and all $x\in \co[K]$, we have
		\begin{equation}
			|A^M_\lambda(f_{K_R})(x)-f(x)|\leq Lr_c(x)+\frac{L^2}{\lambda}\,.
		\end{equation}

\medskip
\noindent If we further assume that  $\tilde{f}$ is a $C^{1,1}$ function such that
$|D\tilde{f}(x)-D\tilde{f}(y)|\leq L|x-y|$ for all $x,\,y\in \mathbb{R}^n$ with $L>0$ the Lipschitz 
constant of $D\tilde{f}$, then for  $\lambda>L$, $M>A_0+\lambda (R+r)^2$ and all $x\in \co[K]$, we have
\begin{equation}\label{Eq.ErEstImpFC11}
	|A^M_\lambda(f_{K_R})(x)-\tilde{f}(x)|\leq
	\frac{L}{4}\left(\frac{\lambda+L/2}{\lambda-L/2}+1\right)r_c^2(x)\,.
\end{equation}
Furthermore, in  this case,  $A^M_\lambda(f_{K_R})$ is an interpolation of $f_K$ in $\R^n$.
\end{prop}

\begin{nota}\label{Rmk.ErEst}
\begin{itemize}
	\item[$(i)$] Using \eqref{Eq.Ineq01}, it follows that the estimates \eqref{Eq.ErEstImpFcnt} and
		\eqref{Eq.ErEstImpFC11} hold with $r_c(x)$ replaced by $d(x)$. Although the resulting estimates are less
		sharp, they have a  clearer meaning in light of the geometric interpretation of the gap $d(x)$.
	\item[$(ii)$] While the assumption of boundedness of the image $f$ is a plausible one, 
		the assumption on the continuity of $f$ seems to be less reasonable for applications to images which might have sharp changes in grayscale intensity.
		However, Example \ref{Ex.AprxJump} at the end of this section,
		illustrates the fact  that our average 
		approximation operator well approximates jump discontinuities.
\end{itemize}
\end{nota}

\medskip
\noindent It is interesting to compare our error estimates \eqref{Eq.ErEstImpFcnt} and \eqref{Eq.ErEstImpFC11} with the error analysis for 
image inpainting discussed in \cite{CK06}. 
Let $\Omega\subset\mathbb{R}^2$ be a smooth domain, which is the damaged area of the image to be reconstructed, 
and let $u$ be a $C^2$ function in a larger domain containing $\bar\Omega$. 
Let $u_0=u$ on $\partial\Omega$ and consider the solution $v$ of the 
boundary value problem $\Delta v(x)=0$ with $v=u_0$ on $\partial \Omega$. The function $v$ is the reconstruction of $u$
within $\Omega$. The error estimate obtained in \cite{CK06} is then given by
\begin{equation}
\label{chanest}
	|v(x)-u(x)|\leq \frac{T\beta^2}{4},\quad x\in\Omega\,,
\end{equation}
where $T=\max\{|\Delta u(x)|,\; x\in\bar\Omega\}$ and $\beta$ is the shorter semi-axis of
any ellipse covering $\Omega$. \cite{CK06} also contains variations of estimate \eqref{chanest}  by
deforming (if possible) a general long thin domain into one for which $\beta$ 
is reasonably small.

\medskip
\noindent Note that in light of Remark \ref{Rmk.ErEst}$(i)$, the error bound \eqref{Eq.ErEstImpFC11} 
		depends explicitly on $d(x)$ and the Lipschitz constant $L$ of the gradient $D\tilde{f}$, which is 
		comparable with the bound $T$ for the Laplacian of $u$. 
Moreover, our assumptions on the smoothness of the domain $\Omega$ and the underlying function are weaker 
		than those considered in \cite{CK06}. In fact, we do not require any smoothness of the boundary 
		$\partial\Omega$. 
	Our estimate is particularly sharp for more general thin domains given its dependance on $d(x)$.
		As remarked in \cite{CK06}, the short semi-axis $\beta^2$ used in the error estimate for harmonic 
		inpainting cannot be replaced by $d^2(x)$ which better accounts for the geometric structure of the 
		damaged area to be inpainted.
Due to the Hausdorff stability property of the average approximation (see \cite[Theorem 4.12]{ZCO16}), 
if $\Omega_\epsilon$ is another domain whose Hausdorff distance to
$\Omega$ is small, we can also obtain   similar results to  estimate \eqref{Eq.ErEstImpFC11} for such domains.

\medskip
\noindent Reference \cite{CK06} contains also error estimates for the TV inpainting model using the energy $\int_\Omega|v(x)|dx$ 
under the Dirichlet condition $v|_{\partial\Omega}=u_0$. However, it is not clear how such estimates can be made rigorous. 
Comparing with Proposition \ref{Pro.EstEr.Inp} where we assumed the underlying 
function to be bounded and uniformly continuous, the TV model, in contrast, allows the
function to have jumps, thus the TV inpainting model tries to preserve such jump discontinuities. 
However, such a model cannot be Hausdorff stable. Also, in order to establish 
the existence of solutions for this model, we note that the boundary condition has to be relaxed. 
Even for the more regular minimal graph energy $\int_{\Omega}\sqrt{1+|Dv(x)|^2}dx$, existence of 
solutions for the  Dirichlet problem may not be guaranteed \cite{Giu84}.
On the other hand, the average approximation always exists and is unique. 
See Example \ref{Ex6.5} in Section \ref{Sec.ProtEx}  for an illustration of this.

\medskip
\noindent Compared with our model for inpainting, we also note that for the relaxed Dirichlet problem of the minimal graph 
or of the TV model, as the boundary value of the solution does not have to agree with the original boundary value, 
extra jumps can be introduced along the boundary. By comparison, since our average 
approximation is continuous, it will not introduce such a jump discontinuity at the boundary.  

\medskip
\noindent One of the motivations for using  TV related models \cite{CS05} for the inpainting problem is that functions 
of bounded variations can have jump discontinuities \cite{AFP00}. Some authors argue that continuous functions cannot
be used to model digital image related functions as functions representing images may have jumps \cite{CS05}.  
However, from the human vision perspective, it is hard to distinguish between a jump discontinuity, where values
change abruptly, and a continuous function with sharp changes within a very small transition layer. 
The following is a simple one-dimensional example showing the effects of our upper, lower and average
compensated convex transforms on a jump function. More explicitly calculated prototype examples 
of inpainting by using our method over jump discontinuity and continuous edges are given in 
Section \ref{Sec.ProtEx}.

%%%%%%%%%%%%%%%%%%%%%%%%%%%%%%%%%%%%%%%%%%%%%%%%%%%%%%%%%%%%%%%%%%%%%%%%%%%%%%%%%%%%

\begin{ex}\label{Ex.AprxJump} 
	Let $f(x)=\sign(x)$ be the sign function defined by $\sign(x)=1$ if $x>0$, $\sign(x)=-1$ if $x<0$. 
	For $\lambda>0$, we have
\begin{equation}\label{Eq.AprxJump}
	\begin{split} 
		C^l_\lambda(f)(x)& = \left\{\begin{array}{ll} 
			\displaystyle	-1, & \displaystyle x\leq 0,\\[1.5ex]
			\displaystyle	1-\lambda(x-\sqrt{2/\lambda})^2, & \displaystyle 0\leq x\leq \sqrt{2/\lambda},\\[1.5ex]
			\displaystyle	1, & \displaystyle x\geq \sqrt{2/\lambda};
			\end{array}\right.\\[1.5ex]
		C^u_\lambda(f)(x)& = \left\{\begin{array}{ll} 
	\displaystyle	-1,				& \displaystyle x\leq -\sqrt{2/\lambda}\,,\\[1.5ex]
	\displaystyle \lambda(x+\sqrt{2/\lambda})^2-1,	& \displaystyle -\sqrt{2/\lambda}\leq x\leq 0\,,\\[1.5ex]
	\displaystyle 1,				& \displaystyle x\geq 0;
			\end{array}\right.\\[1.5ex]
	\displaystyle	\frac{1}{2}(C^l_\lambda(f)(x)+C^u_\lambda(f)(x)) & =\left\{\begin{array}{ll}
	\displaystyle	-1,					   & \displaystyle x\leq -\sqrt{2/\lambda},\\[1.5ex]
	\displaystyle	\frac{\lambda}{2}(x+\sqrt{2/\lambda})^2-1, & \displaystyle -\sqrt{2/\lambda}\leq x\leq 0,\\[1.5ex]
	\displaystyle	1-\frac{\lambda}{2}(x-\sqrt{2/\lambda})^2, & \displaystyle 0\leq x\leq \sqrt{2/\lambda},\\[1.5ex]
	\displaystyle	1,& \displaystyle x\geq \sqrt{2/\lambda};
			\end{array}\right.
	\end{split}
\end{equation}

Figure \ref{Fig.1} displays the graphs of these transforms with $\lambda=100$ which give 
very good approximations of the jump function with the square of the $L^2$-error equal to 
$2\sqrt{2}/(5\sqrt{\lambda})$ for the average approximation 
and equal to $\sqrt{2}/(5\sqrt{\lambda}$ for the lower and upper transform. Therefore these transforms 
can be used quite well to replace the jump discontinuity. For further prototype examples of inpainting with jump discontinuity, 
see Section \ref{Sec.ProtEx}.

\begin{figure}[htp]
  \centerline{
  $\begin{array}{ccc}
    \includegraphics[height=4cm]{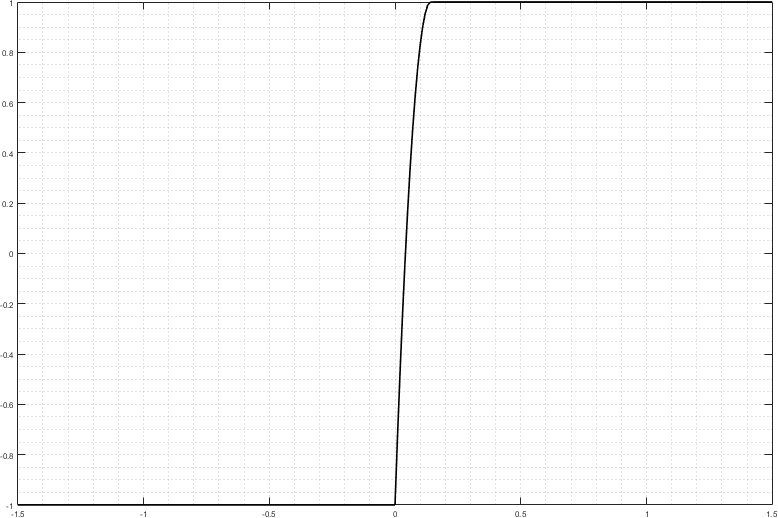}&
    \includegraphics[height=4cm]{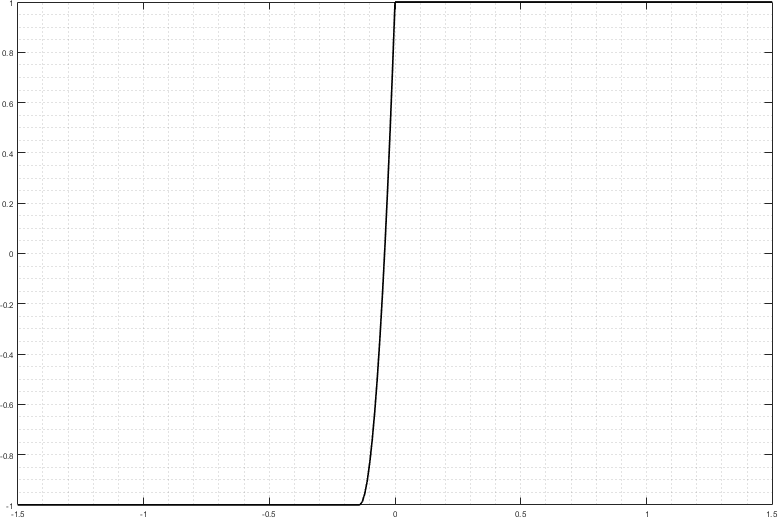}&
    \includegraphics[height=4cm]{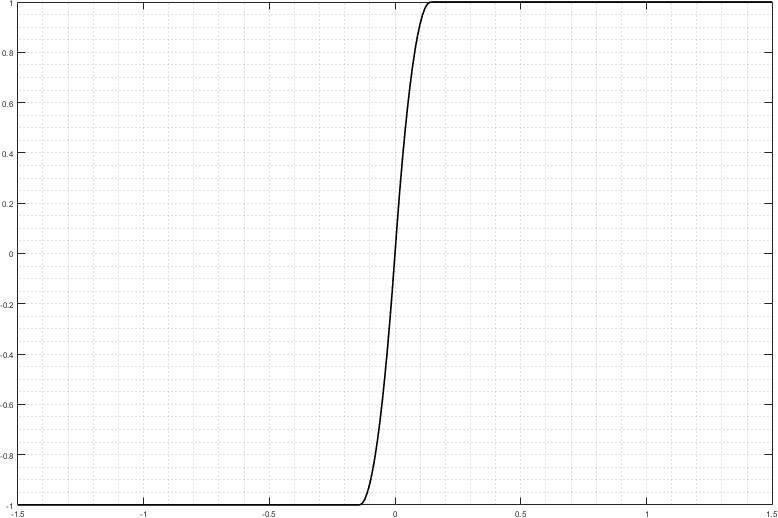}\\
    $(a)$ & $(b)$ & $(c)$ 
\end{array}$}
   \caption{\label{Fig.1}
	$(a)$ Lower transform of the sign function for $\lambda=100$. 
	$(b)$ Upper transform of the sign function for $\lambda=100$. 
	$(c)$ Average approximation of the sign function for $\lambda=100$.}
\end{figure}

\end{ex} 

%%%%%%%%%%%%%%%%%%%%%%%%%%%%%%%%%%%%%%%%%%%%%%%%%%%%%%%%%%%%%%%%%%%%%%
\medskip
\noindent We conclude this section by presenting a result on inpainting in bounded convex domains 
which we state only for continuous functions defined on the closure of the domain. 
For Lipschitz and $C^{1,1}$ functions,  similar results can  be established.

\begin{coro}\label{Cor.Inp}
Suppose $\Omega\subset \mathbb{R}^n$ is a non-empty, bounded, open and convex set and
$U\subset \bar{U}\subset\Omega$ is an open subset whose closure $\bar{U}$ is contained in $\Omega$. 
Suppose $f:\bar{\Omega}\to \mathbb{R}$ is a continuous function. 
Let $\tilde{f}$ be any
bounded uniformly continuous extension of $f$ to $\mathbb{R}^n$ and $\omega$ be the least 
concave majorant of the modulus 
of continuity  of $\tilde{f}$ which is itself a modulus of continuity. 
Let $K=\bar{\Omega}\setminus U$ and define for $M>0$
\begin{equation*} 
	f_K^{M,\infty}(x)=\left\{\begin{array}{ll} 
		 \displaystyle 	f(x),	& \displaystyle  x\in K,\\
		 \displaystyle 	M,	& \displaystyle x\in U,\\
		 \displaystyle 	+\infty,& \displaystyle  x\in \mathbb{R}^n\setminus\bar{\Omega},
			\end{array}\right.
			\quad
	f_K^{-M,-\infty}(x)=\left\{\begin{array}{ll} 
		 \displaystyle 	f(x),	& \displaystyle  x\in K,\\
		 \displaystyle 	-M,	& \displaystyle  x\in U,\\
		 \displaystyle 	-\infty,& \displaystyle  x\in \mathbb{R}^n\setminus\bar{\Omega}.
			\end{array}\right.
\end{equation*}
Then the average approximation in $\bar\Omega$ defined by
\begin{equation}\label{Eq.Inp.AvAprx}
	A^{M;\,\infty}_\lambda(f_K)(x)=
		\frac{1}{2}\left(C^l_\lambda(f_K^{M,+\infty})(x)+ C^u_\lambda(f_K^{-M,-\infty})(x)\right)
\end{equation}
for $x\in \bar\Omega$ satisfies
\begin{equation*} 
	|A^{M,\infty}_\lambda(f_K)(x)-f(x)|\leq \omega(r_c(x)+a/\lambda+\sqrt{b/\lambda})
\end{equation*} 
for all $x\in\bar{\Omega}$, where $r_c(x)$ is the convex density radius of $x\in\bar\Omega$ 
with respect to $K$.
\end{coro}

\begin{nota}
The average approximation defined by \eqref{Eq.Inp.AvAprx} is the same average 
approximation as defined on the bounded domain $\bar\Omega$
\begin{equation*}
	A^{M}_\lambda(f_K;\,\bar{\Omega})(x)=
	\frac{1}{2}\left(C^l_\lambda(f_K^{M};\,\bar{\Omega})(x)+
			C^u_\lambda(f_K^{-M};\,\bar{\Omega})(x)
		\right)
\end{equation*}
for $x\in\bar\Omega$, where $f_K^{M}(x)$ and $f_K^{-M}(x)$ are defined by \eqref{Eq.Def.ExtFnct}, restricted to
$\bar\Omega.$ We can also state the average approximation under the Dirichlet boundary
condition in a similar way. We leave this to interested readers.
\end{nota}

%%%%%%%%%%%%%%%%%%%%%%%%%%%%%%%%%%%%%%%%%%%%%%%%%%%%%%%%%%%%%%%%%%%%%%%%%%%%%%%%%%%%%%%%%%%%%%%%%%%%
%%%%%%%%%%%%%%%%%%%%%%%%%%%%%%%%%%%%%%%%%%%%%%%%%%%%%%%%%%%%%%%%%%%%%%%%%%%%%%%%%%%%%%%%%%%%%%%%%%%%
%%%%%%%%%%%%%%%%%%%%%%%%%%%%%%%%%%%%%%%%%%%%%%%%%%%%%%%%%%%%%%%%%%%%%%%%%%%%%%%%%%%%%%%%%%%%%%%%%%%%

\setcounter{equation}{0}
\section{Prototype Models}\label{Sec.ProtEx}

In this section we present explicitly calculated average approximations
for some particular simple functions of two variables. Recall that such approximations $A_{\lambda}^{\infty}(f_K)$ are obtained by first finding 
lower and upper compensated convex transforms and then taking their arithmetic mean, and  that the 
approximation properties of $A_{\lambda}^{\infty}(f_K)$ hold for $(x,y)\in\co[K]$.
For some examples we  also give  expressions for the constituent lower and upper transforms to help illustrate the construction of the approximations. 
Such examples serve the dual purpose of providing insight into this new
class of approximations based on compensated convexity transforms, and of verifying 
numerical methods for computing such approximations. In fact, in Section \ref{Sec.NumEx} below,
we will see numerical examples that show that, at a sufficient level of magnification, the conditions that occur in practice for the 
approximation of general functions often look essentially like one of these prototypes.

\subsection{Simple prototypes}

%%%%%%%%%%%%%%%%%%%%%%%%%%%%%%%%%%%%%%%%%%%%%%%%%%%%%%%%%%%%%%%%%%%%%%%%%%%%%%%%%%%%%%%%%%%%%
% Example 6.1
%%%%%%%%%%%%%%%%%%%%%%%%%%%%%%%%%%%%%%%%%%%%%%%%%%%%%%%%%%%%%%%%%%%%%%%%%%%%%%%%%%%%%%%%%%%%%

\begin{ex}\label{Ex6.1} 
These two examples give average approximations $A^{\infty}_\lambda(f_K)$ for simple sampled 
functions over non-regular Delaunay cells. In each case, the average approximation is an interpolation of the 
sampled function values.
	
\begin{itemize}
	\item[(i)] Consider the four point set $K=\{(\pm 1,\,0), (0,\,\pm 1)\}$ and define $f_K(1,0)=f_K(0,1)=1$ and
		$f_K(-1,0)=f_K(0,-1)=-1$. The upper and lower compensated convex transforms are 
		then for $\lambda>0$
		\begin{equation*}
			\begin{array}{ll}
			\displaystyle C^l_\lambda(f_K^\infty)(x,y)&=\left\{\begin{array}{ll}
						\displaystyle 2\lambda-1-x+y-\lambda(x^2+y^2),	& \displaystyle \text{if }x\geq -1,\; y\leq 1\text{ and }x\leq y,\\
						\displaystyle 2\lambda-1+x-y-\lambda(x^2+y^2),	& \displaystyle \text{if }y\geq -1,\; x\leq 1 \text{ and } x\geq y,\\
						\displaystyle +\infty,				& \displaystyle \text{if }|x|>1\text{ or }|y|>1;
							\end{array}\right.\\[2ex]
			\displaystyle C^u_\lambda(f_K^{-\infty})(x,y)&=\left\{\begin{array}{ll}
					\displaystyle -2\lambda+1+x+y+\lambda(x^2+y^2),		& \displaystyle \text{if }x\geq -1,\; y\geq -1\text{ and } x+y\leq 0,\\
					\displaystyle -2\lambda+1-x-y+\lambda(x^2+y^2),		& \displaystyle \text{if }x\leq 1,\; y\leq 1\text{ and } x+ y\geq 0,\\
					\displaystyle -\infty,					& \displaystyle \text{if } |x|>1\text{ or } |y|>1\,.
							\end{array}\right.
			\end{array}
		\end{equation*}
		 so that, for $(x,y) \in D:=\co[K]=\{(x,y)\in \mathbb{R}^2:\; |x|\leq 1,\; |y|\leq 1\}$, 
		we have 
		\begin{equation*}
			A^\infty_\lambda(f_K)(x,y)=\left\{\begin{array}{ll}
						\displaystyle y,   & \displaystyle \text{if }x\leq y\text{ and } x+ y\leq 0,\\
						\displaystyle -x,  & \displaystyle \text{if }x\leq y\text{ and } x+y\geq 0,\\
						\displaystyle x	   & \displaystyle \text{if }x\geq y\text{ and } x+y\leq 0,\\
						\displaystyle -y,  & \displaystyle \text{if }x\geq y\text{ and } x+y\geq 0.
					\end{array}\right.		
		\end{equation*}		
		This is the continuous piecewise affine interpolation of $f_K$ inside the square $D$.  The graph of $A^\infty_\lambda(f_K)$ is shown in Figure \ref{Fig.Ex6.1}$(a)$.
	\item[(ii)] Consider the eight point set $K \subset \mathbb{R}^2$ consisting of the eight points on the unit 
		circle  with polar angles $k\pi/4$, $k=0,1,2,\ldots,7$, and define 
		$f_K(\cos(k\pi/4),\sin(k\pi/4))=(-1)^k$. The upper and lower compensated convex 
		transforms are then for $\lambda>0$
		\begin{equation*}
			\begin{array}{ll}
			\displaystyle C^l_\lambda(f_K^\infty)(x,y)&=\left\{\begin{array}{ll}
						\displaystyle \tfrac{\sqrt{2}+1}{\sqrt{2}-1}-\tfrac{2|y|}{\sqrt{2}-1}
							& \displaystyle \text{if }|x|\leq 1,\; |y|\geq 1\;\text{and }
							  \displaystyle |y|+(\sqrt{2}-1)|x|\leq \sqrt{2},\\[1.5ex]
						\displaystyle \tfrac{\sqrt{2}+1}{\sqrt{2}-1}-\tfrac{2|x|}{\sqrt{2}-1}
							& \displaystyle \text{if }|y|\leq -1,\; |x|\geq 1\;\text{and }
							  \displaystyle |x|+(\sqrt{2}-1)|y|\leq \sqrt{2},\\[1.5ex]
						\displaystyle  1 
							& \displaystyle \text{if }|x|\leq 1,\text{ and }|y|\leq 1\\[1.5ex]
						\displaystyle  0		
							& \displaystyle \text{otherwise};
							\end{array}\right.\\[2.5ex]
			\displaystyle C^u_\lambda(f_K^{-\infty})(x,y)&=\left\{\begin{array}{ll}
						\displaystyle \tfrac{\sqrt{2}+1}{\sqrt{2}-1}-\tfrac{\sqrt{2}|x-y|}{\sqrt{2}-1}
							& \displaystyle \text{if }|x+y|\leq \sqrt{2},\; |x-y|\geq \sqrt{2}\;
							  \text{and }\\[1.5ex]
							& \displaystyle \phantom{xxxxx}|x-y|+(\sqrt{2}-1)|x+y|\leq 2,\\[1.5ex]
						\displaystyle \tfrac{\sqrt{2}+1}{\sqrt{2}-1}-\tfrac{\sqrt{2}|x+y|}{\sqrt{2}-1}
							& \displaystyle \text{if }|x+y|\geq \sqrt{2},\; |x+y|\geq \sqrt{2}\;
							  \text{and } \\[1.5ex]
							& \displaystyle  \phantom{xxxxx}|x+y|+(\sqrt{2}-1)|x-y|\leq 2,\\[1.5ex]
						\displaystyle  1 
							& \displaystyle \text{if }|x+y|\leq \sqrt{2}\text{ and }|x-y|\leq \sqrt{2}\\[1.5ex]
						\displaystyle  0		
							& \displaystyle \text{otherwise};
							\end{array}\right.
			\end{array}
		\end{equation*}
		whereas $A^\infty_\lambda(f_K)(x,y)$ is obtained by taking the arithmetic mean of 
		$C^l_\lambda(f_K^\infty)(x,y)$ and $C^u_\lambda(f_K^{-\infty})(x,y)$.
		Figure \ref{Fig.Ex6.1}$(b)$ shows the graph of 
		$A_{\lambda}^{\infty}(f_K)$ in $\co[K]$, which is the inside of the regular octagon with vertices at the eight points of $K$. 
		As in (i), $A_{\lambda}^{\infty}(f_K)$ is a continuous piecewise affine interpolation of $f_K$ in $\co[K]$.
\end{itemize}
\begin{figure}[htbp]
	\centering{$\begin{array}{cc}
		\includegraphics[height=4cm]{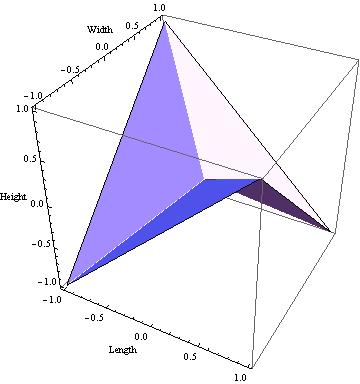}&
		\includegraphics[height=4cm]{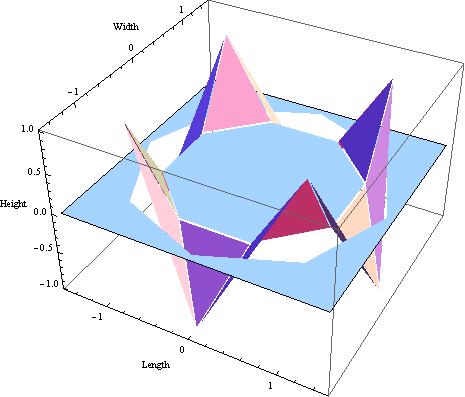}\\
		$(a)$  & $(b)$
	\end{array}$
	}
	\caption{\label{Fig.Ex6.1}  		
	Graphs of the average approximation operators  $A^\infty_\lambda(f_K)$ in Example \ref{Ex6.1}, when $K$ is
 $(a)$   a four point set on the circle of unit radius 
		and  $(b)$  an eight point set on the circle of unit radius. 
		In both (a) and (b), the average approximation operator is an interpolation operator over $\co[K]$.
		}
\end{figure}
\end{ex}

%%%%%%%%%%%%%%%%%%%%%%%%%%%%%%%%%%%%%%%%%%%%%%%%%%%%%%%%%%%%%%%%%%%%%%%%%%%%%%%%%%%%%%%%%%%%%
% Example 6.2
%%%%%%%%%%%%%%%%%%%%%%%%%%%%%%%%%%%%%%%%%%%%%%%%%%%%%%%%%%%%%%%%%%%%%%%%%%%%%%%%%%%%%%%%%%%%%

\begin{figure}[ht]
  \centering{$\begin{array}{cc}
    \includegraphics[height=4cm]{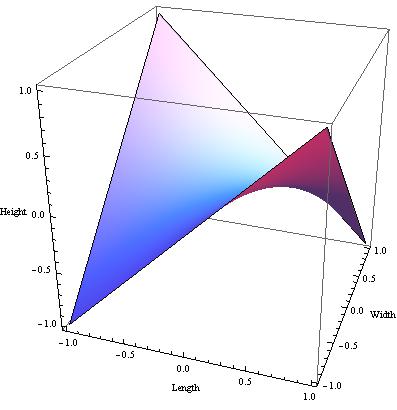}&   
    \includegraphics[height=4cm]{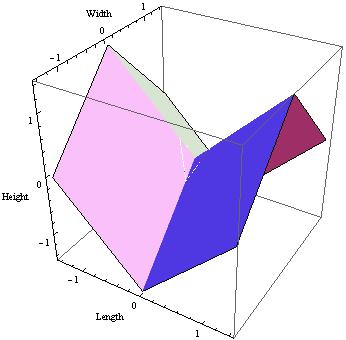}\\
    $(a)$ & $(b)$ 
\end{array}$}
   \caption{\label{Fig.Ex6.2} Graphs of the average approximation operators  $A^\infty_\lambda(f_K)$ in Example 
   \ref{Ex6.2}$(i)$ and $(ii)$, respectively.} \end{figure}

\begin{ex}\label{Ex6.2}
These two examples give  average approximations $A^{\infty}_\lambda(f_K)$ for unbounded sets $K$ with  $\co[K]=\mathbb{R}^2$.
\begin{itemize}
	\item[(i)]  Consider the set $K=\ell_{-}\cup \ell_+$ with $\ell_-=\{(x,x),\,x\in\mathbb{R}\}$,
		$\ell_+=\{(y,-y),\,y\in\mathbb{R}\}$, and define $f_K(x,x)=-x^2$ and $f_K(y,-y)=y^2$. 
		To simplify the calculations,  first consider the scaled and rotated 
		function $g_{\tilde{K}}$ defined on the set $\tilde{K}=\{(x,0), x\in \R\}\cup\{(0,y), y\in \R\}$, with
		$g_{\tilde{K}}(x,0)=-x^2$ and $g_{\tilde{K}}(0,y)=y^2$. 
		Then for $(x,y)\in\mathbb{R}^2$, the lower and upper compensated convex transforms of $g_{\tilde{K}}$ are  
		\begin{equation*}
			  C^l_\lambda(g_{\tilde{K}}^\infty)(x,y)    = y^2+2|x||y|-x^2,\;\;\;\;\;	
			  C^u_\lambda(g_{\tilde{K}}^{-\infty})(x,y) = y^2-2|x||y|-x^2,
		\end{equation*}
	and the average approximation of $g_{\tilde{K}}$ is 
\begin{equation*}
	A^\infty_\lambda (g_{\tilde{K}})(x,y)	=\frac{1}{2}\left(C^l_\lambda(g_{\tilde{K}}^\infty)(x,y)+C^u_\lambda(g_{\tilde{K}}^{-\infty})(x,y)\right)=y^2-x^2\,.
\end{equation*}
	The average approximation $A^\infty_\lambda (f_K)$ of $f_K$ is then obtained from $A^\infty_\lambda (g_{\tilde{K}})$ via a change of 
	variables, and is 
	\begin{equation*}
		A^\infty_\lambda (f_K)(x,y)=\frac{1}{2}\left(A^\infty_\lambda (g_{\tilde{K}})\left(\frac{x+y}{\sqrt{2}},\, \frac{x-y}{\sqrt{2}}\right)\right)=-xy\,.
	\end{equation*}
		Figure \ref{Fig.Ex6.2}$(a)$ shows the graph of $A^\infty_\lambda(f_K)$.
	
	\item[(ii)] Let $K=\{(x,0),\,x\in\mathbb{R}\}\cup \{(0,y),\,y\in\mathbb{R}\}$ and define
			$f_K$ by $f_K(x,0)=|x|$ for $x\in \mathbb{R}$ and $f_K(0,y)=-|y|$ for $y\in \mathbb{R}$.
		For  $(x,y)\in \mathbb{R}^2$, the lower and upper compensated convex transforms of $f_K$ are 
		\begin{equation*}
	\begin{array}{ll}
	\displaystyle  C^l_\lambda(f^\infty_K)(x,y)&=\left\{\begin{array}{ll} 
		\displaystyle	2|x|-\tfrac{1}{4\lambda}-\lambda(x^2+y^2),	& \displaystyle \text{if }|x|+|y|\leq \tfrac{1}{2\lambda},\\[1.5ex]
		\displaystyle  	|x|+2\lambda|x||y|-|y|,				& \displaystyle \text{if }|x|+|y|\geq \tfrac{1}{2\lambda},
						\end{array}
						\right.\\[2ex]
	\displaystyle  C^u_\lambda(f^{-\infty}_K)(x,y)	& = \left\{\begin{array}{ll} 
		\displaystyle	-2|y|+\tfrac{1}{4\lambda}+\lambda(x^2+y^2),	& \displaystyle \text{if }|x|+|y|\leq \frac{1}{2\lambda},\\[1.5ex]
		\displaystyle  	|x|-2\lambda|x||y|-|y|,				& \displaystyle \text{if }|x|+|y|\geq \frac{1}{2\lambda},
						\end{array}
						\right.
	\end{array}
\end{equation*}
	and the average approximation operator is 
\begin{equation*}
	 A^\infty_\lambda(f_K)(x,y)=|x|-|y|\,,
\end{equation*}
which here coincides with the natural interpolation of $f_K$ by the piecewise affine function $f(x,y)=|x|-|y|$, $(x,y)\in\mathbb{R}^2$.
	 The graph of $A^\infty_\lambda(f_K)$ is shown in Figure \ref{Fig.Ex6.2}$(b)$. 
\end{itemize}
\end{ex}

%%%%%%%%%%%%%%%%%%%%%%%%%%%%%%%%%%%%%%%%%%%%%%%%%%%%%%%%%%%%%%%%%%%%%%%%%%%%%%%%%%%%%%%%%%%%%%%%%%%%%%%%%%%%%%%%%

\subsection{Inpainting prototypes}

\noindent Examples \ref{Ex6.3} and \ref{Ex6.4} are prototype models for the inpainting problem. 
Our question is, to what extent our method can preserve singularities on the boundary based on the given boundary values. 
Our calculations show that if the domain is narrow 
and similar singular boundary values appear on both sides of the narrow gap,
the inpainting function $A^\infty_\lambda(f_K)$ can preserve the singular shape across the gap, 
subject to a  $\lambda$-dependent regularisation of the singularity due to the local 
smoothing effect of the compensated convex transforms.

%%%%%%%%%%%%%%%%%%%%%%%%%%%%%%%%%%%%%%%%%%%%%%%%%%%%%%%%%%%%%%%%%%%%%%%%%%%%%%%%%%%%%%%%%%%%%%%%%%%%%%%%%%%%%%%%%

%%%%%%%%%%%%%%%%%%%%%%%%%%%%%%%%%%%%%%%%%%%%%%%%%%%%%%%%%%%%%%%%%%%%%%%%%%%%%%%%%%%%%%%%%%%%%
% Example 6.3
%%%%%%%%%%%%%%%%%%%%%%%%%%%%%%%%%%%%%%%%%%%%%%%%%%%%%%%%%%%%%%%%%%%%%%%%%%%%%%%%%%%%%%%%%%%%%

\begin{ex}\label{Ex6.3}
	
\begin{itemize}	
\item[(i)] For $r>0$, $h>0$, let $K=\{(\pm r, y),\; |y|\leq h\}\subset \mathbb{R}^2$ i.e.  
two parallel line segments a distance $r$ apart (see Figure \ref{Fig.Ex6.3a}(a)), 
and define $f_K(\pm r,y)=1-|y|$. Let $D=\co[K]=\{(x,y)\in \mathbb{R}^2:\; |x|\leq r,\; |y|\leq h\}$. 
Then for $\lambda>1/2h$, 
	\begin{equation*}	
	\begin{array}{ll}
		\displaystyle C^l_\lambda(f^\infty_K)(x,y)  & = \left\{\begin{array}{ll} 
								\displaystyle	1-\tfrac{1}{4\lambda}+\lambda r^2-\lambda x^2-\lambda y^2,	
									& \displaystyle \text{if }|x|\leq r \text{ and } |y|\leq \tfrac{1}{2\lambda},\\[1.5ex]
								\displaystyle	1+\lambda r^2-\lambda x^2-|y|,	
							&  \displaystyle \text{if }  |x|\leq r \text{ and } \tfrac{1}{2\lambda}\leq |y|\leq h,\\[1.5ex]
								\displaystyle	+\infty & \text{otherwise,} 
							\end{array}\right.\\[2ex]
		\displaystyle C^u_\lambda(f^{-\infty}_K)(x,y) & = \left\{\begin{array}{ll} 
									\displaystyle 1 - \lambda r^2 + \lambda x^2 - |y|,
									& \displaystyle \text{if } |x|\leq r \text{ and } |y|\leq h;\\[1.5ex]
								\displaystyle	-\infty & \text{otherwise,} 
							\end{array}\right.
	\end{array}
	\end{equation*}
and for $(x,y)\in D$, the average approximation operator is 
	\begin{equation*}	
		A^\infty_\lambda(f_K)(x,y) =\left\{\begin{array}{ll}
								\displaystyle 1-\tfrac{1}{8\lambda}-\tfrac{\lambda y^2}{2}-\tfrac{|y|}{2}, &  
										\displaystyle \text{if }  |x|\leq r \text{ and } |y|\leq \tfrac{1}{2\lambda},\\[1.5ex]
								\displaystyle 1-|y|,	&  
										\displaystyle \text{if }  |x|\leq r, \text{ and } \tfrac{1}{2\lambda}\leq |y|\leq h\,.
							\end{array}\right.
	\end{equation*}	
The graph of $A^{\infty}_\lambda(f_K)$ is shown in Figure \ref{Fig.Ex6.3a}$(b)$.\\

\noindent  Note that this example shows that if we only sample the two gables $K$ of the roof, 
 the whole roof can be recovered well for any $r>0$ and $h>0$. On the other hand, we will see in the next example that   the situation is more complicated if 
 the other two sides, $(x,\pm h)$ for $|x|\leq r$,  are added to the sample set.
\begin{figure}[htbp]
  \centering{$\begin{array}{cc}
		\includegraphics[width=0.4\textwidth]{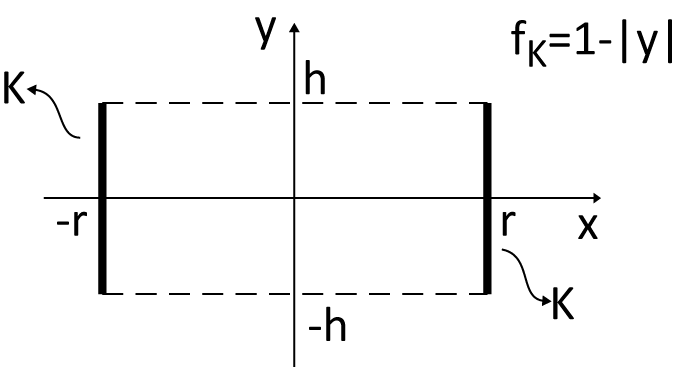}&
		\includegraphics[width=0.35\textwidth]{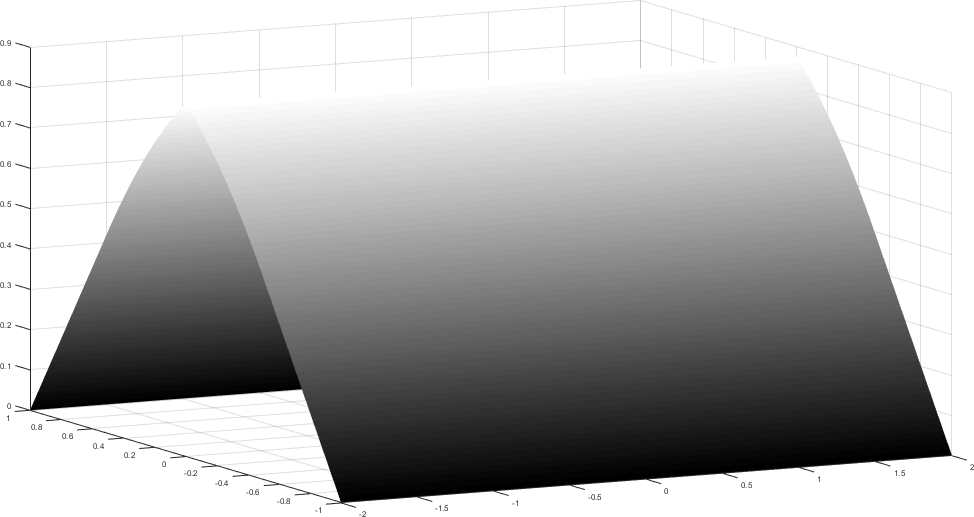}\\
		(a)&(b)
	\end{array}$}
   \caption{\label{Fig.Ex6.3a} Example \ref{Ex6.3}$(i)$.
		$(a)$ The sample set $K$ shown in bold, with the sample function $f_K=1-|y|$.
		$(b)$ Graph of $A^\infty_\lambda(f_K)$ for $\lambda=1$.}
\end{figure}
\item[(ii)] Next let $D=\{(x,y),\; |x|\leq r,\; |y|\leq h\}$ with $h>0$ and $r>0$, take the sample set $K=\partial D=\{(\pm r,y),\; |y|\leq h\}\cup\{(x,\pm h),\; |x|\leq r\}$, and define
	\begin{equation*}
		f_{K}(x,y)=\left\{\begin{array}{ll} 
				\displaystyle h-|y|, & \displaystyle x=\pm r,\; |y|\leq h,\\[1.5ex]
				\displaystyle 	0,   & \displaystyle y=\pm h, \; |x|\leq r\,.
		\end{array}\right.
	\end{equation*}
For large $\lambda$, the shape of $A_\lambda^\infty(f_{K}(x,y))$ in $D$ now depends on whether $h>r$, $h<r$ or $h=r$. 
\begin{itemize}
	\item[(a)] If $h>r$, the two gables of the roof $h-|y|$ at $x = \pm r$ are close to each other and we have a very good approximation of the whole roof $h-|y|$ for $(x,y) \in D$ when $\lambda$ is sufficiently large. For $(x,y) \in D$, the approximation 
	$A_\lambda^\infty(f_{K}(x,y))$ is
\begin{equation*}		
	A^\infty_\lambda(f_K)(x,y)=\left\{\begin{array}{ll}
					\displaystyle h-\tfrac{1}{4\lambda}-\lambda y^2, 
						& \displaystyle \text{if } |y|\leq \tfrac{1}{2\lambda} \text{ and } |x|\leq r,\\[1.5ex]
					\displaystyle h-|y|, & \displaystyle \text{if }  \tfrac{1}{2\lambda}\leq |y|\leq h \text{ and } |x|\leq r\,,
				\end{array}\right.
\end{equation*}
which yields the explicit error estimate
\[
	|A^\infty_\lambda(f_K)(x,y) - f(x,y)| \leq \tfrac{1}{8 \lambda}\,.
\]
In particular, the ridge of the roof is preserved well in this case. 
%%%%%%%%%%%%%%%%%%%%%%%%%%%%%%%%%
	\item[(b)] If $h=r$ and $\lambda>0$ is large,  the roof dips in the middle, while the `ridge' is still preserved. 
	\item[(c)] If $h<r$ and $\lambda>0$ is large, the roof falls inside  $D=\co[K]$ and touches the ground. In this case, the ridge is no longer preserved at all.   	
\end{itemize}	
\noindent In summary, the average approximation can approximate well the non-smooth function given on two sides of $K$ provided the 
two gables are close enough. In this case, we could say that by symmetry we have a behaviour similar to the one seen 
in Example \ref{Ex6.3}$(a)$. As opposite, when the two gables are far apart, i.e. when $h/r<1$, it is somehow the effect of 
$f_K=0$ on the sides $y=\pm h$ to make it feel its presence, by having a zero interpolation in the middle of the domain. We stress again 
that this situation is different from the one seen in Example \ref{Ex6.3}$(a)$ where $f_K$ was sampled only on the sides $x=\pm r$. 
Figure \ref{Fig.Ex6.3b} shows the graphs of $A_{\lambda}^{\infty}$ in each of the three cases, together with the sample set K.
\end{itemize}
\end{ex}

\begin{figure}[htbp]
  \centerline{$\begin{array}{cc}
		   \includegraphics[width=0.45\textwidth]{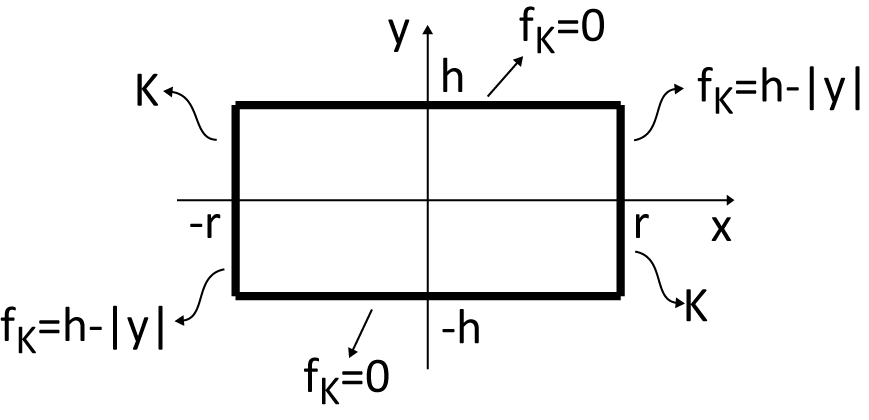}&
		   \includegraphics[width=0.35\textwidth]{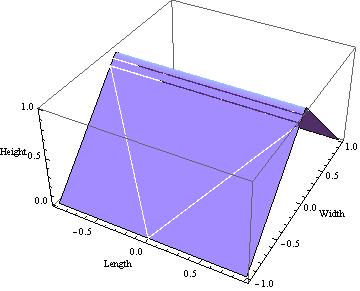}\\
		    (a) & (b) \\
		\includegraphics[width=0.35\textwidth]{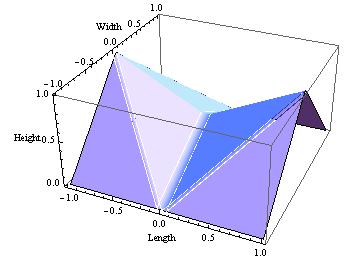}&
		\includegraphics[width=0.35\textwidth]{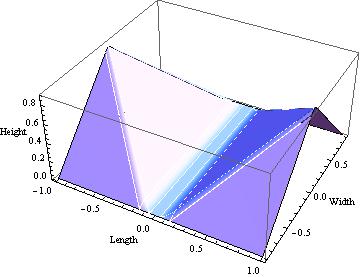}\\
		    (c) & (d)
	\end{array}$
	}   
   \caption{\label{Fig.Ex6.3b} Example \ref{Ex6.3}$(ii)$.
   $(a)$ The sample set $K$ shown in bold, with the sample function $f_K$.
   Average approximation in $D$ for the following parameters: 
   $(b)$ $h=1$, $r=0.9$, $\lambda=10$. 
   $(c)$ $h=1$, $r=1$, $\lambda=10$.
   $(d)$ $h=0.9$, $r=1$, $\lambda=10$.
}
\end{figure}

\bigskip
\noindent A preliminary one-dimensional prototype of the inpainting of a region when the boundary values have discontinuities was 
given in Example \ref{Ex.AprxJump}. We next explore how our inpainting method can preserve jumps in a two-dimensional example. 

%%%%%%%%%%%%%%%%%%%%%%%%%%%%%%%%%%%%%%%%%%%%%%%%%%%%%%%%%%%%%%%%%%%%%%%%%%%%%%%%%%%%%%%%%%%%%
% Example 6.4
%%%%%%%%%%%%%%%%%%%%%%%%%%%%%%%%%%%%%%%%%%%%%%%%%%%%%%%%%%%%%%%%%%%%%%%%%%%%%%%%%%%%%%%%%%%%%

\begin{figure}[htb]
  \centerline{$\begin{array}{ccc}
		\includegraphics[width=0.25\textwidth]{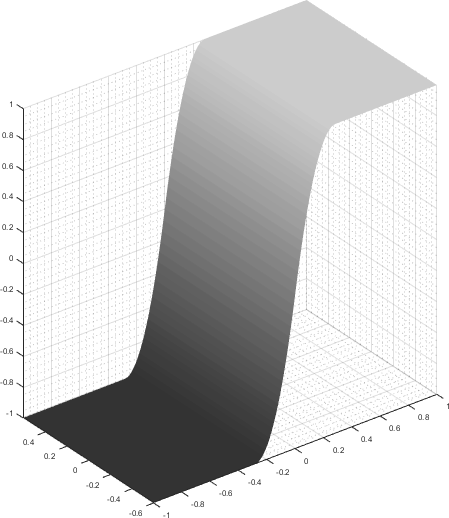} &
		\includegraphics[width=0.30\textwidth]{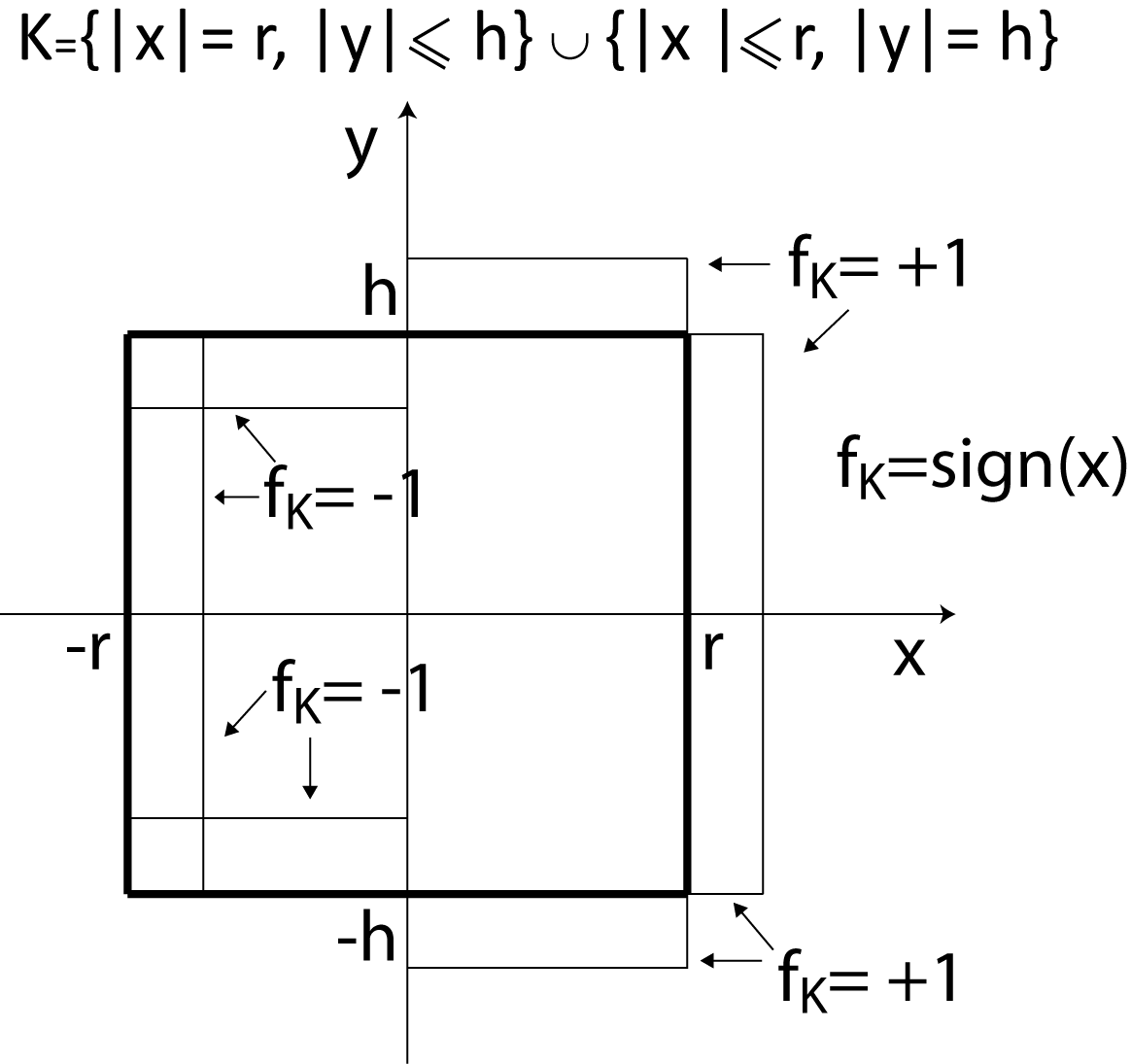}\\
		(a) & (b) 
	\end{array}
	$}
   \caption{\label{Fig.Ex6.4} Example \ref{Ex6.4}.
	Inpainting in the closed set $D=\{(x,y),\; |x|\leq l,\; |y|\leq h\}$ 
	by the boundary value of the sign function on the sample set $K=\partial D$. 
	$(a)$ Graph  of $A^\infty_\lambda(f_{K})$ for $h=0.6$, $r=1$, $\lambda=25$, showing that the jump is preserved across the 
	domain $D$.
	$(b)$ Sample set $K$ shown in bold  with the sampled function $f_K=\sign(x)$. 
}
\end{figure}

\begin{ex}\label{Ex6.4}
 Consider the inpainting of the region
$D=\{(x,y),\; |x|\leq r,\; |y|\leq h\}$, for $r$, $h>0$, in the case of narrow gap, that is, when $h<r$.
The sample set is the boundary of the domain $D$, that is, $K=\partial D$, and the sample function $f_K$ is taken as $f_K(x,y)=\sign(x)$. 
Then for $\lambda>0$ large enough, the average approximation operator is in fact 
given by \eqref{Eq.AprxJump}, that is, for $(x,y)\in D$,
\begin{equation*}
		A^\infty_\lambda(f_{K})(x,y)=\left\{\begin{array}{ll}
	\displaystyle	-1,	& \displaystyle\text{if }x\leq -\sqrt{2/\lambda}\text{ and }|y|\leq h\,,\\[1.5ex]
	\displaystyle	\tfrac{\lambda}{2}(x+\sqrt{2/\lambda})^2-1, & \displaystyle\text{if } -\sqrt{2/\lambda}\leq x\leq 0\text{ and }|y|\leq h\,,\\[1.5ex]
	\displaystyle	1-\tfrac{\lambda}{2}(x-\sqrt{2/\lambda})^2, & \displaystyle\text{if } 0\leq x\leq \sqrt{2/\lambda}\text{ and }|y|\leq h\,,\\[1.5ex]
	\displaystyle	1,   & \displaystyle\text{if }	x\geq \sqrt{2/\lambda}\text{ and }|y|\leq h\,.
				\end{array}\right. 
\end{equation*}
Figure \ref{Fig.Ex6.4}$(a)$ shows the graph of the average approximation $A^\infty_\lambda(f_{K})$ in this case. 
The approximation $A^\infty_\lambda(f_{K})(x,y)$ is different from $\sign(x)$ in the range $[-\sqrt{2/\lambda},\sqrt{2/\lambda}]\times[-h,\,h]$ due to
the smoothing effect of the compensated transform in the neighbourhood of the singularity. 
The width of such a neighbourhood depends on $\sqrt{\lambda}^{-1}$. The full recovery of the sign function in $D$ requires taking the limit $\lim_{\lambda \to \infty} A^{\infty}_{\lambda}(f_K)(x,y)$.\\

\noindent Note that if, on the other hand, $h>r$, the gap is `wide' and the graph of $A^\infty_\lambda(f_{K})$ starts to collapse in the middle of the domain, 
similar to what happens in Example \ref{Ex6.3}(ii)(c). In the collapsed region, the approximation looks like
an affine function connecting the two sides $\{x = \pm r\}$ of $D$ on which $f_K$ is given by the constants $+1$, when $x= +r$,  and $-1$, when $x=-r$. 

%%%%%%%%%%%%%%%%%%%%%%%%%%%%%%%%%%%%%%%%%%%%%%%%%%%%%%%%%%%%%%%%%%
\end{ex}

\subsection{Level-set prototypes}
%
%%%%%%%%%%%%%%%%%%%%%%%%%%%%%%%%%%%%%%%%%%%%%%%%%%%%%%%%%%%%%%%%%%%%%%%%%%%%%%%%%%%%%%%%%%%%%
% Example 6.5
%%%%%%%%%%%%%%%%%%%%%%%%%%%%%%%%%%%%%%%%%%%%%%%%%%%%%%%%%%%%%%%%%%%%%%%%%%%%%%%%%%%%%%%%%%%%%
We next present prototype models for the approximation of functions sampled on contour lines. 

\begin{ex} \label{Ex6.5}
This example examines the behaviour of $A_{\lambda}^{\infty}(f_K)$ when the contour lines of $f$ 
are $(i)$ smooth and $(ii)$ not smooth.
\begin{itemize}
	\item[(i)] For $0<r<R$, let $K=\Gamma_r\cup \Gamma_R$ with $\Gamma_r$ and $\Gamma_R$ circles of radius $r$
	and $R$, respectively, as displayed in Figure \ref{Fig.Ex6.5a}$(a)$, and define the sample function $f_K$  by 
	$f_K(x,y)=0$ for $(x,y)\in \Gamma_r$ and $f_K(x,y)=M>0$ if $(x,y)\in \Gamma_R$. 
		Then for $\lambda>M/(R^2-r^2)$, 
		\begin{equation*}
	\begin{split}
		C^u_\lambda(f^{-\infty}_K)(x,y) & =\left\{\begin{array}{l} 
				\displaystyle	M+\lambda(x^2+y^2-r^2), \quad\displaystyle\text{if } \sqrt{x^2+y^2}\leq r,\\[1.5ex]
				\displaystyle	\lambda(x^2+y^2-R^2)+\tfrac{M+\lambda(R^2-r^2)}{R-r}(R-\sqrt{x^2+y^2}),
					\quad\displaystyle\text{if } r\leq \sqrt{x^2+y^2}\leq R\,,
				\end{array}\right.\\[2ex]
		C^l_\lambda(f^\infty_K)(x,y) & =\left\{\begin{array}{l} 
				\displaystyle	M+\lambda(r^2-x^2-y^2),\quad\displaystyle\text{if }  \sqrt{x^2+y^2}\leq r,\\[1.5ex]
				\displaystyle	\lambda(R^2-x^2-y^2)-\tfrac{\lambda(R^2-r^2)-M}{R-r}(R-\sqrt{x^2+y^2}),\quad
					\displaystyle\text{if } r\leq \sqrt{x^2+y^2}\leq R\,,
				\end{array}\right.
	\end{split}
\end{equation*}
so that for $(x,y)\in D=\co[K]=\{(x,y):x^2+y^2\leq R^2\}$, the average approximation  $A^\infty_\lambda(f_K)$ is 
\begin{equation*}
	A^\infty_\lambda(f_K)(x,y) =\left\{\begin{array}{ll} 
		\displaystyle	M,				& \displaystyle\text{if } \sqrt{x^2+y^2}\leq r,\\[1.5ex]
		\displaystyle	\tfrac{M(R-\sqrt{x^2+y^2)}}{R-r}, & \displaystyle\text{if } r\leq \sqrt{x^2+y^2}\leq R\,.
			\end{array}\right.
\end{equation*}
The graph of $A^\infty_\lambda(f_K)$ is shown in Figure \ref{Fig.Ex6.5a}$(b)$.

\medskip
Note that a common method for the interpolation of function values assigned on contour lines is to solve the Dirichlet problem for the minimal surface equation $\ddiv\frac{Du}{\sqrt{1+|Du|^2}}=0$
	over the annulus domain $r\leq \sqrt{x^2+y^2}\leq R$ with boundary conditions $u(x,y)=0$ if $(x,y)\in \Gamma_r$ 
	and $u(x,y)=M$ if $(x,y)\in \Gamma_R$.
	It is then known that this problem does not have a regular solution \cite{Giu84}. Moreover,  the interpolation 
	obtained by solving the total variation equation $\ddiv\frac{Du}{|Du|}=0$ faces the same type of issue, because to obtain its  
	numerical solution, the denominator $|Du|$ is usually replaced by the term $\sqrt{\epsilon^2+|Du|^2}$, 
	thus obtaining  the scaled minimal surface equation
	$\ddiv\frac{Du}{\sqrt{\epsilon+|Du|^2}}=0$ whose solution, as mentioned above, may not be regular. 
	As a result, these models must be relaxed and one must look for generalised solutions \cite{GMS79}. 
In contrast, the method we propose yields 
	the natural, easy to compute and expected interpolation $A^\infty_\lambda(f_K)$ between the two level lines.
	\begin{figure}[htbp]
		\centerline{$\begin{array}{cc}
					\includegraphics[width=0.3\textwidth]{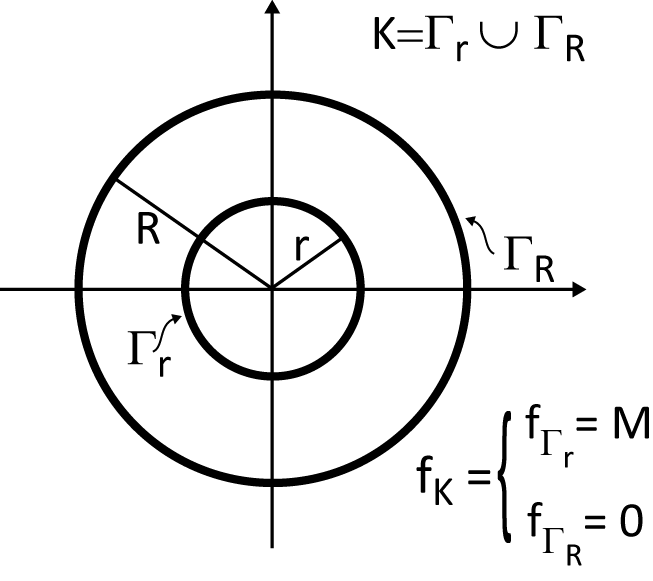}
					&\includegraphics[width=0.25\textwidth]{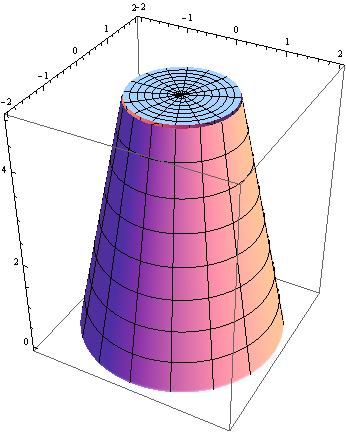}\\
					(a) & (b)
					\end{array}$}
		\caption{\label{Fig.Ex6.5a} Example \ref{Ex6.5}$(i)$.
			$(a)$ Sample set $K$ given by the two circular level lines $\Gamma_r$ and $\Gamma_R$ with 
			$f_K(x,y)=0$ for $(x,y)\in \Gamma_r$ and $f_K(x,y)=M>0$ if $(x,y)\in \Gamma_R$.
			$(b)$ Graph of $A^\infty_\lambda(f_K)$ with $r=1$, $R=2$, $M=5$ and $\lambda=10$.}
	\end{figure}

\item[(ii)] For $a,\,\lambda>0$, consider the sample set $K=K_1\cup K_2$ with $K_1=\{(x,\, y):\; |y|=ax, \; x\geq 0\}$ 
	and $K_2=\{(x,\, y):\; |y|=a(x-\frac{\sqrt{1+a^2}}{a\sqrt{\lambda}}), \; x\geq \frac{\sqrt{1+a^2}}{a\sqrt{\lambda}}\}$, and define the sample function $f_K$ by $f_K(x,y)=1$ for $(x,y)\in K_1$ and $f_K(x,y)=2$ for $(x,y)\in K_2$.
	The set $K$ along with $f_K$ are shown in Figure \ref{Fig.Ex6.5b}$(a)$.
	For $(x,y)\in D=\co[K]=\{(x,y):|y|\leq ax,\; x\geq 0\}$, the  average approximation operator $A^\infty_\lambda(f_K)$ is 	
	\begin{equation*}
		A^\infty_\lambda(f_K)(x,y)=\left\{\begin{array}{l} 
		\displaystyle 1,\quad\displaystyle\text{if } |y|\leq a x \text{ and } 0 \leq x \leq\tfrac{ 1}{a\sqrt{1 + a^2}},\\[1.5ex]
		\displaystyle 1 + \tfrac{\sqrt{1 + a^2}\left(-\tfrac{1}{a\sqrt{1 + a^2}}+ x\right)}{a},\quad\displaystyle\text{if }
			\displaystyle	x\geq \tfrac{1}{a\sqrt{1 + a^2}}\text{ and } 
			\displaystyle	\tfrac{x + a |y|}{\sqrt{1 + a^2}} \leq \tfrac{1}{a\sqrt{\lambda}},\\[1.5ex]
		\displaystyle 2 - \left|\tfrac{1}{\sqrt{\lambda}} + \tfrac{-a x + |y|}{\sqrt{1 + a^2}}\right| \quad
		\displaystyle\text{if } -\tfrac{1}{\sqrt{\lambda}} \leq\tfrac{-a x + |y|}{\sqrt{1 + a^2}} \leq 0\text{ and }
		\displaystyle \tfrac{1}{a\sqrt{\lambda}} \leq \tfrac{x + a|y|}{\sqrt{1 + a^2}},\\[1.5ex]
		\displaystyle 2,\quad \displaystyle\text{if } \tfrac{-a x + |y|}{\sqrt{1 + a^2}} \leq -\tfrac{1}{\sqrt{\lambda}}\text{ and } 
					\displaystyle  x \geq \tfrac{\sqrt{1 + a^2}}{a\sqrt{\lambda}}\,.
			\end{array}\right.
	\end{equation*}
	The graph of $A^\infty_\lambda(f_K)$ is displayed in Figure \ref{Fig.Ex6.5b}$(b)$. Note that the interpolation 
	$A^\infty_\lambda(f_K)$ takes the constant value $1$, which is the value given on the level set $K_1$, inside a triangle next to the corner of $K_1$,
	which is then pieced continuously to $K_2$ by a continuous piecewise affine function.
\end{itemize}
\end{ex}

\begin{figure}[htbp]
	\centerline{$\begin{array}{ccc}
			\includegraphics[width=0.3\textwidth]{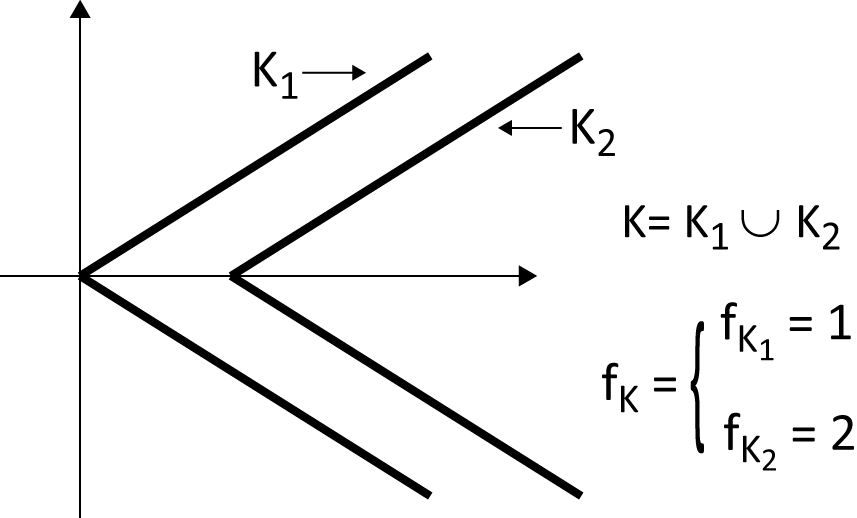}
			&\includegraphics[width=0.35\textwidth]{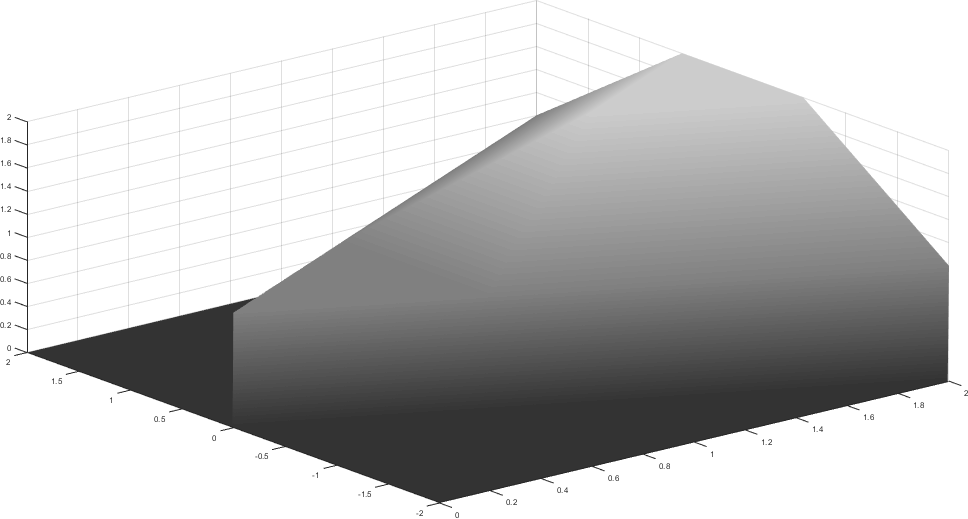}
			&\includegraphics[height=0.25\textwidth,angle=90]{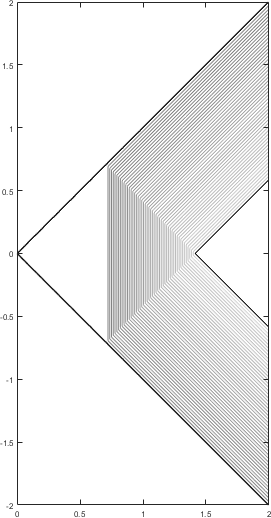}\\
			(a)&(b)&(c)
		\end{array}$}
	\caption{\label{Fig.Ex6.5b} Example \ref{Ex6.5}$(ii)$.
		$(a)$ Sample set $K$ given by two non-smooth level sets $K_1$ and $K_2$ with $a=1$
		and $f_K(x,y)=0$ for $(x,y)\in K_1$ and $f_K(x,y)=2$ if $(x,y)\in K_2$.
		$(b)$ Graph of $A^\infty_\lambda(f_K)$ with $\lambda=1$.
		$(c)$ Isolines of $A^\infty_\lambda(f_K)$.}
\end{figure}

%%%%%%%%%%%%%%%%%%%%%%%%%%%%%%%%%%%%%%%%%%%%%%%%%%%%%%%%%%%%%%%%%%%%%%%%%%%%%%%%%%%%%%%%%%%%%
% Example 6.6
%%%%%%%%%%%%%%%%%%%%%%%%%%%%%%%%%%%%%%%%%%%%%%%%%%%%%%%%%%%%%%%%%%%%%%%%%%%%%%%%%%%%%%%%%%%%%

\noindent We conclude this section with a prototype example of level-set approximation for a function 
with  a jump discontinuity at the point $(0,0)$.

\begin{figure}[ht]
  \centerline{
	$\begin{array}{cc}
		\includegraphics[width=0.3\textwidth]{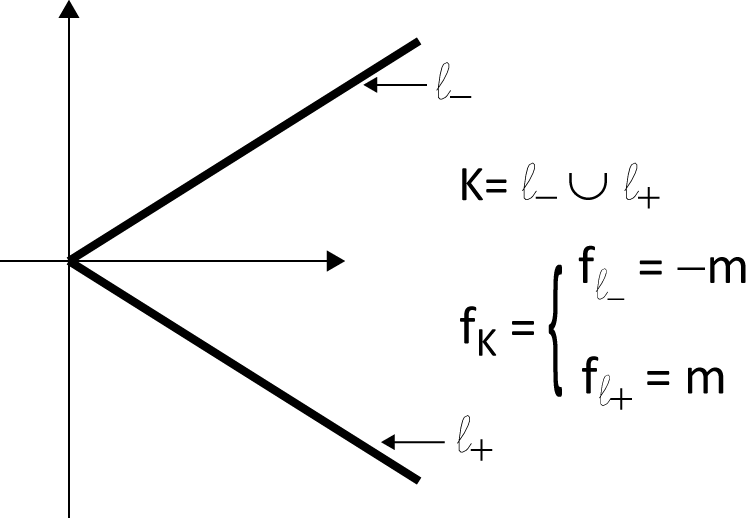}
		&\includegraphics[width=0.3\textwidth]{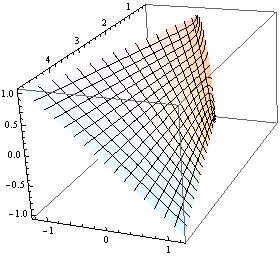}\\
		(a) & (b)
	\end{array}$}
   \caption{\label{Fig.Ex6.6} Example \ref{Ex6.6}.
		$(a)$ Sampled set $K$ with the definition of $f_K$ that presents a discontinuity jump at $(0,0)$.
		$(b)$ Graph of $A_\lambda^\infty(f_K)$ with $\alpha=0.25$, $m=1$ and $\lambda=5$.}
\end{figure}

\begin{ex} \label{Ex6.6}
For $\alpha,\,m>0$, consider the sample set $K$ given by $K=\ell_{+}\cup \ell_{-}$ with $\ell_{+}=\{(x,y),\, y=-\alpha x, \, x>0\}$
and $\ell_-=\{(x,y),\, y=\alpha x, \, x>0\}$, and define $f_K(x,y)=m$ on $\ell_+$ and  
$f_K(x,y)=-m$ on $\ell_+$. 
The set $K$ along with $f_K$ are displayed in Figure \ref{Fig.Ex6.6}$(a)$. 
To describe the average approximation of $f_K$ in  
$\co[K]=\{(x,y),\,|y|\leq \alpha x,\,x>0\}$ which we denote by $S_+$, 
we use a parameterised description of the graph $(x,y, A^{\infty}_\lambda(f_K)(x,y))$ 
in terms of two new parameters. This is to avoid solving
quartic equations when we find the lower and the upper transforms. Let $c_\lambda=2m/\lambda$.
To calculate the lower transform $C^l_\lambda(f_K^\infty)$ in $S_+$
we need to find the common tangent planes for $f^\infty_K(x,y)+\lambda(x^2+y^2)$ of both 
$\ell_+$ and $\ell_-$. 
%%%%%%%%%%%%%%%%%%%%%%%%%%%%%%%%%%%%%%%%%%%%%%%%%%%%%%%%%
We can write the coordinates of the
convex envelope as $(x,y, \co[f^\infty_K(x,y)+\lambda(x^2+y^2)])$ by
\begin{equation*}
	\left(\frac{(1-t_l)\sqrt{s_l^2+c_\lambda}+t_ls_l}{\sqrt{1+\alpha^2}},\,
	\frac{-\alpha(1-t_l)\sqrt{s_l^2+c_\lambda}+\alpha t_ls_l}{\sqrt{1+\alpha^2}},\;
	\lambda s_l^2+2\lambda(1-t_l)c_\lambda -m\right),
\end{equation*}
where $0\leq t_l\leq 1$ and $s_l\geq 0$.
Similarly, the coordinates of $(x,y, \co[\lambda(x^2+y^2)-f^{-\infty}_K(x,y)])$ are
\begin{equation*}
	\left(\frac{(1-t_u)s_ut_u\sqrt{s_u^2+c_\lambda}}{\sqrt{1+\alpha^2}},\,
	\frac{-\alpha(1-t_u)s_u+\alpha t_u\sqrt{s_u^2+c_\lambda}}{\sqrt{1+\alpha^2}},\;
	\lambda s_u^2+2\lambda t_uc_\lambda -m\right),
\end{equation*}
where $0\leq t_u\leq 1$ and $s_u\geq 0$.
However, the $(x,y)$ coordinates in these two cases do not represent the same points. 
Therefore we need to set them equal so that
\begin{equation}\label{Eq.Crd.tltu}
	t_l=\frac{\sqrt{s_u^2+c_\lambda}\left(\sqrt{s_l^2+c_\lambda}-s_u\right)}
	 {\sqrt{s_u^2+c_\lambda}\sqrt{s_l^2+c_\lambda}-s_us_l},\quad
	 t_u=\frac{s_l\left(\sqrt{s_l^2+c_\lambda}-s_u\right)}
	 {\sqrt{s_u^2+c_\lambda}\sqrt{s_l^2+c_\lambda}-s_us_l}\,.
\end{equation}
As $0\leq t_l,\, t_u\leq 1$, we see that $|s_u^2-s_l^2|\leq c_\lambda$. Thus if we let
\begin{equation*}
	x(s_l,s_u)=\frac{(1-t_l)\sqrt{s_l^2+c_\lambda}+t_ls_l}{\sqrt{1+\alpha^2}},\quad
	y(s_l,s_u)=\frac{-\alpha(1-t_l)\sqrt{s_l^2+c_\lambda}+\alpha t_ls_l}{\sqrt{1+\alpha^2}}\,,
\end{equation*}
and
\begin{equation*}
	A_\lambda^\infty(f_K)(s_l,s_u)=\frac{1}{2}\Big(\lambda (s_l^2-s_u^2)+
			2\lambda c_\lambda( (1-t_l-t_u))\Big)\,,
\end{equation*}
the graph of the average approximation of $f_K$ in the sector $S_+$ defined above is
\begin{equation*}
	\Gamma_{S_+,\lambda}=
		\left\{\left(x(s_l,s_u),\, y(s_l,s_u),\, 
			A_\lambda^\infty(f_K)(s_l,s_u) \right),\; 
		s_u\geq 0,\, s_l\geq 0,\; |s_u^2-s_l^2|\leq c_\lambda\right\},
\end{equation*}
where $t_l$ and $t_u$ are given by \eqref{Eq.Crd.tltu}.
\medskip

\noindent Although it is not easy to write the graph in the standard Euclidean system,
 observe that the graph is smooth in the interior region 
$\{(s_l,s_u),\, s_l>0,\, s_u>0, |s^2_l-s^2_u|<c_\lambda\}$. 
By our construction, we also note that the surface $\Gamma_{S_+,\lambda}$ is formed by 
the average of two families of parameterised line segments. Also when $\lambda>0$ is large, 
outside a small sector, say, $S_+^\lambda=\{|y|\leq \alpha, 0<x<2\sqrt{2m/\lambda}\}$, 
our formula is an interpolation in $S_+\setminus S^\lambda_+$.
Figure \ref{Fig.Ex6.6}$(b)$ shows a portion of the graph of $A_\lambda^\infty(f_K)$.
\end{ex}

%%%%%%%%%%%%%%%%%%%%%%%%%%%%%%%%%%%%%%%%%%%%%%%%%%%%%%%%%%%%%%%%%%%%%%%%%%%%%%%%%%%%%%%%%%%%%%%%%%%%
%%%%%%%%%%%%%%%%%%%%%%%%%%%%%%%%%%%%%%%%%%%%%%%%%%%%%%%%%%%%%%%%%%%%%%%%%%%%%%%%%%%%%%%%%%%%%%%%%%%%
%%%%%%%%%%%%%%%%%%%%%%%%%%%%%%%%%%%%%%%%%%%%%%%%%%%%%%%%%%%%%%%%%%%%%%%%%%%%%%%%%%%%%%%%%%%%%%%%%%%%

\setcounter{equation}{0}
\section{Numerical Examples}\label{Sec.NumEx}

\begin{figure}[ht]
  \centerline{\includegraphics[width=0.8\textwidth]{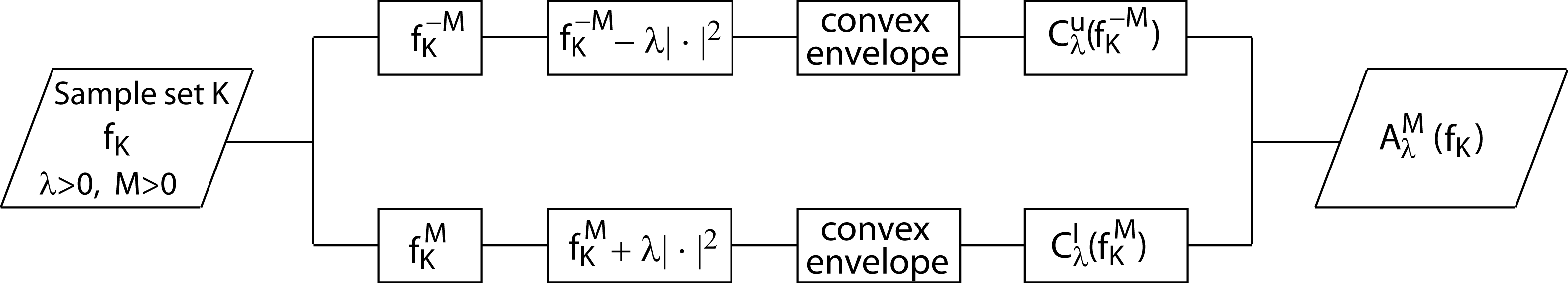}}
   \caption{\label{Fig.Ex7} Flow chart for the numerical evaluation of $A_{\lambda}^M(f_K)$.}
\end{figure}

For more complicated sets $K$ and functions $f_K$, the average approximation operators
$A_{\lambda}^M(f_K)$ and $A_{\lambda}^{\infty}(f_K)$ must be evaluated numerically. Figure \ref{Fig.Ex7}
sketches the steps needed for their implementation. It is noted that the
numerical realization relies mainly on the availability of numerical schemes
for computing the upper and lower transform of a given function, which in turn means the availability 
of schemes to compute the convex envelope of a function. 
Because of the locality property of the compensated convex transforms
(see for instance Theorem 3.10 in \cite{ZOC15a}, where quantitative estimates of the neighbourhood size are also given),
it is possible to develop fast schemes that depend only on the local behaviour of the input function.
This is in sharp contrast to the evaluation of the convex envelope of a function which is a global evaluation.
In the current context, we consider a generalization of 
the scheme introduced in \cite{Obe07} which is briefly summarized in Algorithm \ref{Algo:CnvxEnv} 
and described below. 
Given a uniform grid of points $x_k\in\mathbb{R}^n$, equally spaced with grid size $h$, let us denote by
$S_{x_k}$ the $d-$point stencil of $\mathbb{R}^n$ with center at $x_K$ defined as 
$S_{x_k}=\{x_k+hr, |r|_{\infty}\leq 1, r\in\mathbb{Z}^n\}$ with $|\cdot|_{\infty}$ the $\ell^{\infty}$-norm of 
$r\in\mathbb{Z}^n$ and $d=\#(S)$. At each grid point $x_k$ we compute the 
convex envelope of $f$ at $x_k$ by an iterative scheme where each iteration step $m$ is given by
\[
	(\co f)_m(x_k)=\min\Big\{ 
		f(x_k),\,\sum\lambda_i (\co f)_{m-1}(x_i),\,\,
		\sum \lambda_i=1,\,\lambda_i\geq 0,\,x_i\in S_{x_k}\Big\}
\]
with the minimum taken between $f(x_k)$ and only some convex combinations at the stencil grid points.
For the full algorithmic and implementation details of the
scheme, the convex combinations that one needs to take, and its convergence analysis 
we refer to \cite{ZOC16m}. 

\begin{algorithm}
\begin{algorithmic}[1]
	\STATE{Set $\displaystyle m=1,\,(\co f)_{0}=f,\,\,tol$}
	\STATE{$\displaystyle \epsilon=\|f\|_{L^2}$}
	\WHILE{$\displaystyle \epsilon>tol$}
	\STATE{$\displaystyle \forall x_k, \quad (\co f)_{m}(x_k)=\min\Big\{ 
		f(x_k),\,\sum\lambda_i (\co f)_{m-1}(x_i),\,\,
		\sum \lambda_i=1,\,\lambda_i\geq 0,\,x_i\in S_{x_k}\Big\}$
		}
	\STATE{$\displaystyle \epsilon=\|(\co f)_m-(\co f)_{m-1}\|_{L^2}$}
	\STATE{$m\leftarrow m+1$}
	\ENDWHILE
\end{algorithmic}
\caption{\label{Algo:CnvxEnv} Conceptual implementation of the scheme that computes the convex envelope of $f$.}
\end{algorithm}

In this section, we present some illustrative numerical experiments of the applications described above, namely, for
surface reconstruction from contour lines, point clouds
and image inpainting. For the first two applications,
we discuss examples of approximation of a smooth function, of a continuous but non-differentiable function 
and of a discontinuous function. The quality of the approximation is measured by computing the relative $L^2$-error  
\begin{equation}\label{Eq.RelErr}
	\epsilon=\frac{\|f-A_\lambda^M(f_K)\|_{L^2(\Omega)}}{\|f\|_{L^2(\Omega)}}\,,
\end{equation}
where $f$ is the original function that we want to approximate and $A_\lambda^M(f_K)$ is the average approximation
of the sample $f_K$ of $f$ over $K$.
We mainly postpone a thorough comparison with other state-of-art methods to forthcoming papers, 
just giving  some 
first comparisons with the AMLE method presented in \cite{ACGR02,CMS98}
and applied to surface reconstruction and image inpainting. Image denoising for salt \& pepper noise and  
image inpainting were solved by the TV-model described in 
\cite{CHN05} and in \cite{Get12}, respectively.

We conclude this short introduction by stating that at least for the examples and methods we have considered here, 
we have observed higher accuracy of the $A_{\lambda}^M(f_K)$ interpolant and the 
faster execution time for its numerical evaluation compared to the other methods.

%%%%%%%%%%%%%%%%%%%%%%%%%%%%%%%%%%%%%%%%%%%%%%%%%%%%%%%%%%%%%%%%%%%%%%%%%%%%%%%%%%%%%%%%%%%%%%%%%%%%%%%%%%%%%

\subsection{Surface reconstruction from contour lines}\label{Sec7.1}

%%%%%%%%%%%%%%%%%%%%%%%%%%%%%%%%%%%%%%%%%%%%%%%%%%%%%%%%%%%%%%%%%%%%%%%%%%%%%%%%%%%%%%%%%%%%%%%%%%%%%%%%%%%%%

We describe next some numerical experiments on surface reconstruction from sectional contours. 
This is the problem of reconstructing the graph of a function $f$ by knowing only some
contour lines of $f$, and has applications in medical imaging, computer graphics, reverse engineering 
and terrain modelling, among others. 
The underlying function $f:\R^2\supset\Omega \to \R $ is assumed to have various regularity properties. 
Consider first the reconstruction of an infinitely differentiable function given by the Franke test function \cite{Fra79}, 
and then the reconstruction of functions with less regularity. In addition to the relative 
$L^2$-error $\epsilon$ defined by \eqref{Eq.RelErr},
which gives a measure of how close $A_\lambda^M(f_K)$ is to $f$, we also compute 
\begin{equation}\label{Eq.RelErr.Restr}
	\epsilon_K=\frac{\|f_K-A_\lambda^M(f_K)_K\|_{L^2(K)}}{\|f_K\|_{L^2(K)}}\,,
\end{equation}
where $f_K$ is the sample function and $A_\lambda^M(f_K)_K$ the restriction of 
$A_\lambda^M(f_K)$ to $K$, to assess the quality of $A_\lambda^M(f_K)$ as an interpolant of $f_K$. We will thus verify that in the examples where $f$ is continuous, 
the average approximation $A_\lambda^M(f_K)$ represents an interpolation of $f_K$, consistently  
with the theoretical results established in Section \ref{Sec.LSAprx}. 

%%%%%%%%%%%%%%%%%%%%%%%%%%%%%%%%%%%%%%%%%%%%%%%%%%%%%%

\subsubsection{Franke test function}\label{Sec7.1.1}

%%%%%%%%%%%%%%%%%%%%%%%%%%%%%%%%%%%%%%%%%%%%%%%%%%%%%%

\begin{figure}[htbp]
  \centerline{$\begin{array}{ccc}
		\includegraphics[width=0.50\textwidth]{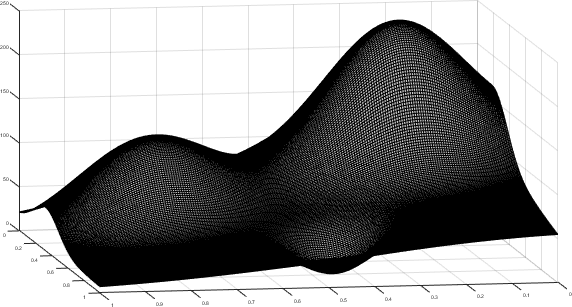}
		&\includegraphics[width=0.25\textwidth]{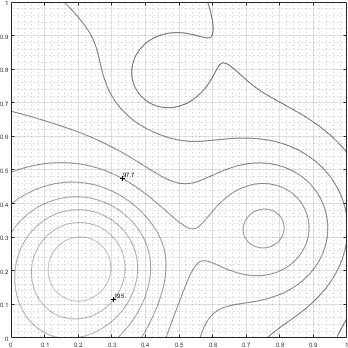}
		&\includegraphics[width=0.25\textwidth]{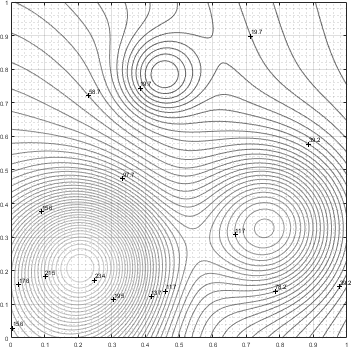}\\
		(a)&(b)&(c)
	     \end{array}$
	     }
   \caption{\label{Fig.NumEx7.1.1a}  Example \ref{Sec7.1.1}.
	$(a)$ Graph of the Franke test function $f$ defined by Equation \eqref{Eq.FrankeFnct}.
	$(b)$ Sample set $K$ of $10$-contour lines of $f$ at equally spaced heights equal to $(\max(f)-\min(f))/10$, 
		defining the sample function $f_K$.
	$(c)$ Sample set $K$ of $50$-contour lines of $f$ at equally spaced heights equal to $(\max(f)-\min(f))/50$, 
		defining the sample function $f_K$.
		}
\end{figure}

\begin{figure}[htbp]
  \centerline{$	\begin{array}{cc}
			\includegraphics[width=0.50\textwidth]{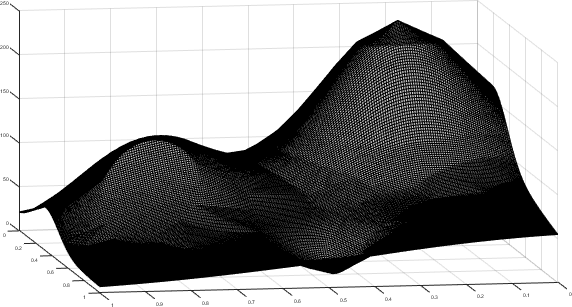}&
			\includegraphics[width=0.25\textwidth]{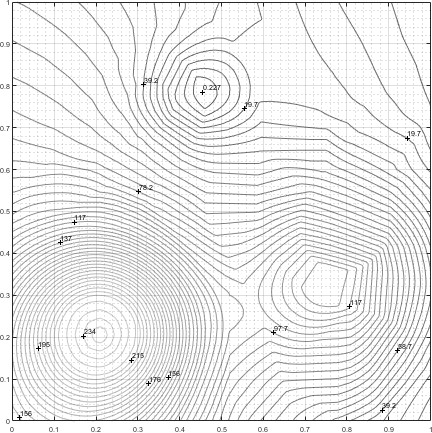}\\
			(a) & (b) \\
			\includegraphics[width=0.50\textwidth]{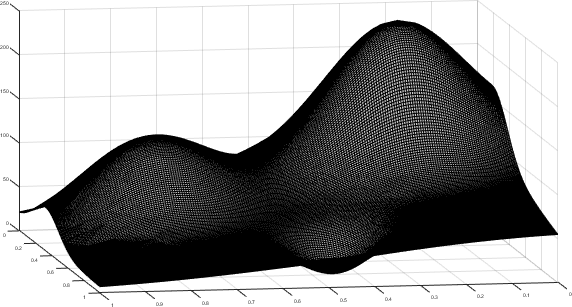}&
			\includegraphics[width=0.25\textwidth]{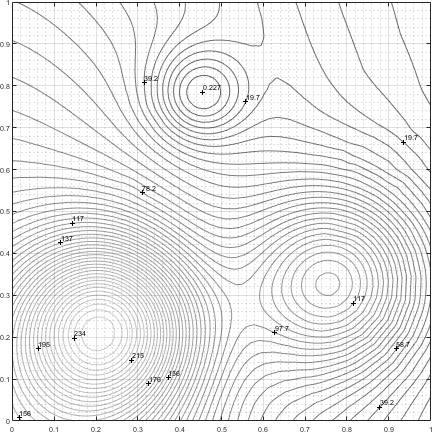}\\
			(c) & (d) 
		\end{array}$
	     }
   \caption{\label{Fig.NumEx7.1.1b} Example \ref{Sec7.1.1}.
	$(a)$ Graph of the interpolation function $A_{\lambda}^M(f_K)$ computed for $\lambda=1\cdot 10^4$, $M=1\cdot 10^5$,
		and corresponding to the set $K$ of $10$-contour lines of $f$ displayed in 
		Figure \ref{Fig.NumEx7.1.1a}$(b)$. 
		Relative $L^2$-Errors: $\epsilon=0.01986$, $\epsilon_K=3.33\cdot 10^{-15}$.  
	$(b)$ Isolines of $A_{\lambda}^M(f_K)$ at equally spaced heights equal to $(\max(f)-\min(f))/50$.
	$(c)$ Graph of the interpolation function $A_{\lambda}^M(f_K)$ computed for $\lambda=1\cdot 10^4$, $M=\cdot 10^5$,
		and corresponding to the set $K$ of $50$-contour lines of $f$ displayed in 
		Figure \ref{Fig.NumEx7.1.1a}$(d)$.
		Relative $L^2$-Errors: $\epsilon=0.0021$, $\epsilon_K=2.62\cdot 10^{-15}$. 
	$(d)$ Isolines of $A_{\lambda}^M(f_K)$ at equally spaced heights equal to $(\max(f)-\min(f))/50$.
	}
\end{figure}

\begin{figure}[htbp]
  \centerline{$\begin{array}{cc}
	\includegraphics[width=0.50\textwidth]{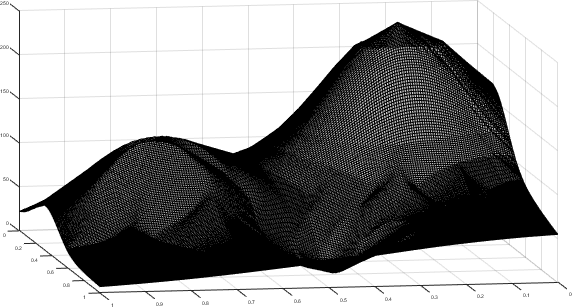}&
	\includegraphics[width=0.25\textwidth]{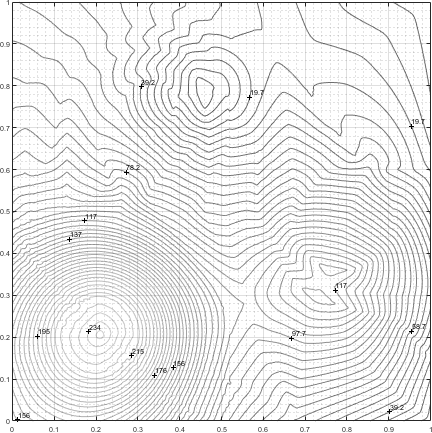}\\
		(a) & (b) \\
	\includegraphics[width=0.50\textwidth]{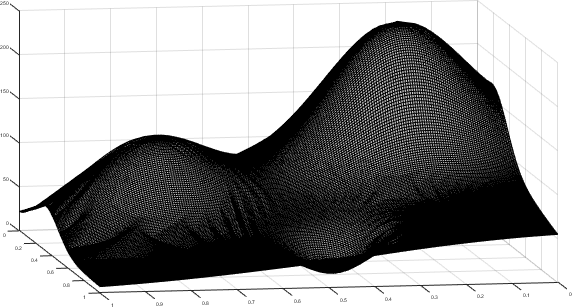}&
	\includegraphics[width=0.25\textwidth]{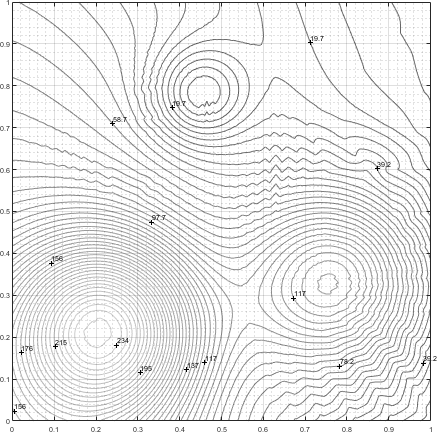}\\
		(c) & (d) 
   \end{array}$
   }
   \caption{\label{Fig.NumEx7.1.1c} Example \ref{Sec7.1.1}.
	$(a)$ Graph of the AMLE interpolation function of $f_K$ with $K$ the set of $10$-contour lines of
		$f$ displayed in Figure \ref{Fig.NumEx7.1.1a}$(b)$.
		Relative $L^2$-Error $\epsilon=0.0338$.  
	$(b)$ Isolines of the AMLE interpolation function of $f_K$ at equally spaced heights equal to $(\max(f)-\min(f))/50$.
	$(c)$ Graph of the AMLE interpolation function of $f_K$ with $K$ the set of $50$-contour lines of
		$f$ displayed in Figure \ref{Fig.NumEx7.1.1a}$(d)$.
		Relative $L^2$-Error $\epsilon=0.0101$.
	$(d)$ Isolines of the AMLE interpolation function of $f_K$ at equally spaced heights equal to $(\max(f)-\min(f))/50$}.
\end{figure}

The Franke function was introduced in \cite{Fra79} as one of the test functions for the evaluation of methods for scattered data 
interpolation \cite{Fra82}. The function consists of two Gaussian peaks and a sharper Gaussian dip superimposed on a surface sloping toward 
the first quadrant \cite{Fra79} and is defined by
\begin{equation}\label{Eq.FrankeFnct}
	\begin{split}
		f(x,y)	& = \frac{3}{4}e^{-\left((9x-2)^2+(9y-2)^2\right)/4}+\frac{3}{4}
			   e^{-\left((9x+1)^2/49+(9y+1)^2/10\right)}+
			   \frac{1}{2}e^{-\left((9x-7)^2/4-(9y-3)^2\right)/4}\\[1.5ex]
			& - \frac{1}{5}e^{-\left((9x-4)^2+(9y-7)^2\right)}\,.
	\end{split}
\end{equation}
Consider $f$ defined in the unit square $\Omega=]0,\,1[^2$. Its graph is displayed in Figure \ref{Fig.NumEx7.1.1a}$(a)$.
Approximations using two different sets of contour lines have been computed by applying the methods 
described in this paper and by the AMLE model introduced in \cite{CMS98} and applied in \cite{ACGR02} to the interpolation 
of digital elevation models.
The two sets of contour lines 
consist of $10$ and $50$ equally spaced level lines, respectively. 
Given the smoothness of $f$, the isolines are also smooth curves. The two sample sets are 
displayed in Figure \ref{Fig.NumEx7.1.1a}$(b)$ and Figure \ref{Fig.NumEx7.1.1a}$(c)$, respectively,
whereas the graph of the corresponding average approximations $A_{\lambda}^M(f_K)$ are shown 
in Figure \ref{Fig.NumEx7.1.1b}$(a)$ and Figure \ref{Fig.NumEx7.1.1b}$(c)$. 
Figure \ref{Fig.NumEx7.1.1b}$(b)$ and Figure \ref{Fig.NumEx7.1.1b}$(d)$ display, on the other hand,
the corresponding contour lines which, compared to the same equally spaced level lines of $f$ displayed
in Figure \ref{Fig.NumEx7.1.1a}$(c)$ show a good quality of the reconstruction given by $A_{\lambda}^M(f_K)$. This is 
also confirmed by the values of the relative $L^2$-error $\epsilon$ equal to $0.01986$
and $0.00218$ for the two sample sets $K$ of contour lines, respectively. 
Note the clear reduction of  error by increasing the density of the data set.
For the two average approximations, the value of $\epsilon_K$ is of the order of $10^{-15}$, confirming that 
the average approximation $A_\lambda^M(f_K)$ interpolates exactly $f_K$. 

Figure \ref{Fig.NumEx7.1.1c} displays the reconstruction obtained by the AMLE method. 
The numerical results were obtained by using the MatLab code described in \cite{PS16}.
In this case, for a number of iterations equal to $10^6$, we found a relative $L^2$-error higher than the one 
generated by $A_{\lambda}^M(f_K)$ with $\epsilon$ equal to $0.0338$ and $0.0101$ for the two sample set $K$ 
of $10$ and $50$ level lines, respectively. Consistently with the findings of \cite{LMS13}, also here we 
find that the AMLE interpolation generates additional kinks which are not present in $f$ and might be the
cause for the reduced quality of the approximation compared to $A_{\lambda}^M(f_K)$. 

%%%%%%%%%%%%%%%%%%%%%%%%%%%%%%%%%%%%%%%%%%%%%%%%%%%%%%

\subsubsection{Continuous piecewise affine function}\label{Sec7.1.2}

%%%%%%%%%%%%%%%%%%%%%%%%%%%%%%%%%%%%%%%%%%%%%%%%%%%%%%

\begin{figure}[htbp]
  \centerline{$\begin{array}{cc}
		\includegraphics[width=0.30\textwidth]{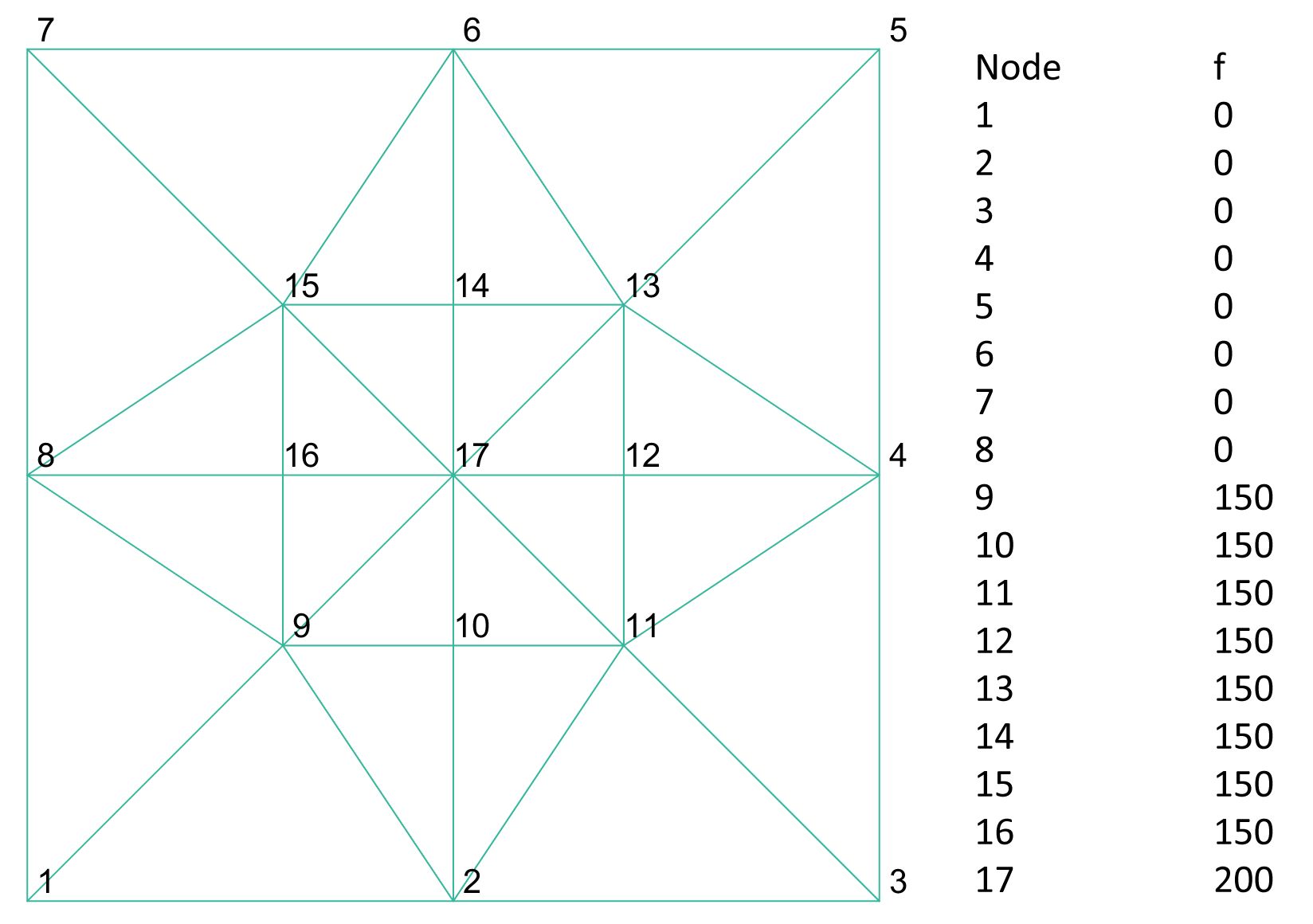}&
		\includegraphics[width=0.50\textwidth]{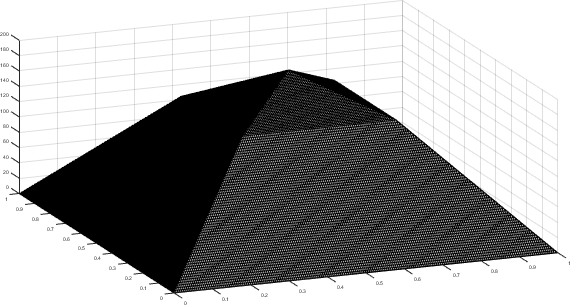}\\
		    (a)&(b)\\
		\includegraphics[width=0.25\textwidth]{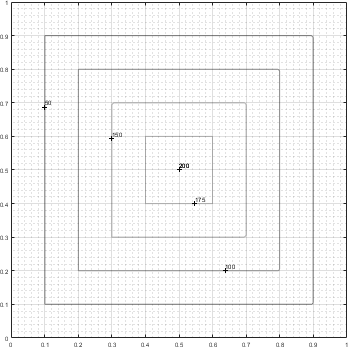}&
		\includegraphics[width=0.25\textwidth]{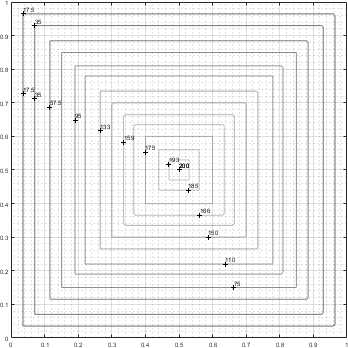}\\
		    (c)&(d)
		\end{array}$}
   \caption{\label{Fig.NumEx7.1.2a}  Example \ref{Sec7.1.2}.
	$(a)$ Triangulation with nodal values used to construct a continuous piecewise affine function. 
	$(b)$ Graph of $f$ associated with the triangulation defined in $(a)$.
	$(c)$ Sample set $K$ of $6$-contour line of $f$, defining the sample function $f_K$. 
	$(d)$ Sample set $K$ of $15$-contour line of $f$, defining the sample function $f_K$.
	}
\end{figure}

\begin{figure}[htbp]
  \centerline{$\begin{array}{cc}
\includegraphics[width=0.50\textwidth]{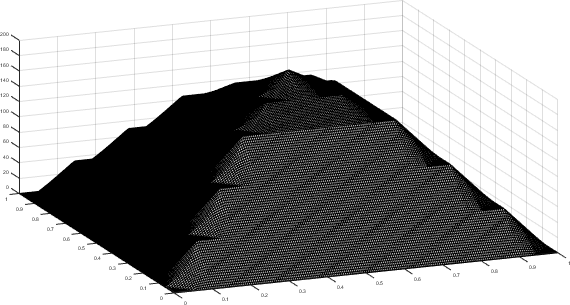}&
\includegraphics[width=0.25\textwidth]{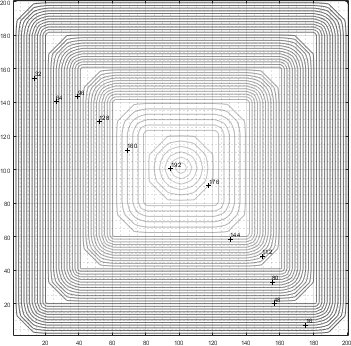}\\
	    (a)&(b)\\
\includegraphics[width=0.50\textwidth]{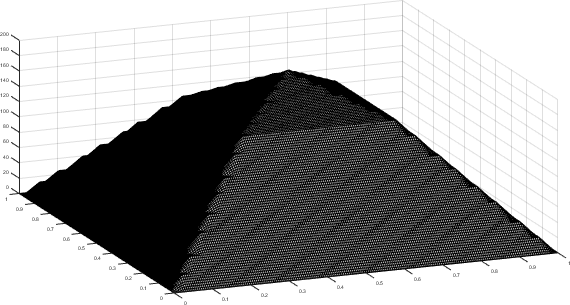}&
\includegraphics[width=0.25\textwidth]{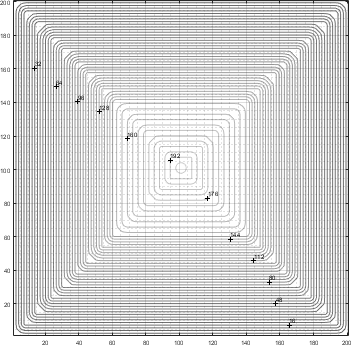}\\
	    (c)&(d)
    \end{array}$}
	\caption{\label{Fig.NumEx7.1.2b}  Example \ref{Sec7.1.2}.
	$(a)$ Graph of the interpolation function $A_{\lambda}^M(f_K)$ with $K$ given in Figure \ref{Fig.NumEx7.1.2a}$(c)$, 
		and $\lambda=1\cdot 10^5$, $M=1\cdot 10^5$, $tol=10^{-9}$.
		Relative $L^2$-Errors: $\epsilon=0.019302$, $\epsilon_K=4.50\cdot 10^{-16}$. 
	$(b)$ Isolines of $A_{\lambda}^M(f_K)$ at equally spaced heights equal to $(\max(f)-\min(f))/50$.		
	$(c)$ Graph of the interpolation function $A_{\lambda}^M(f_K)$  with $K$ given in Figure \ref{Fig.NumEx7.1.2a}$(d)$, 
		and $\lambda=1\cdot 10^5$, $M=1\cdot 10^5$, $tol=10^{-9}$. 
		Relative $L^2$-Errors: $\epsilon=0.004805$, $\epsilon_K=8.68\cdot 10^{-16}$.
	$(d)$ Isolines of $A_{\lambda}^M(f_K)$ at equally spaced heights equal to $(\max(f)-\min(f))/50$.
		}
\end{figure}

\begin{figure}[htbp]
  \centerline{$\begin{array}{cc}
	\includegraphics[width=0.50\textwidth]{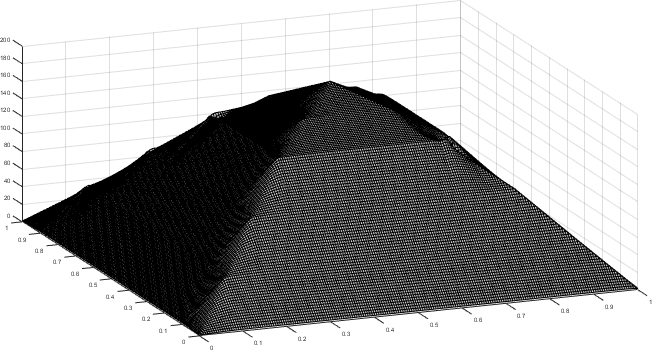}&
	\includegraphics[width=0.25\textwidth]{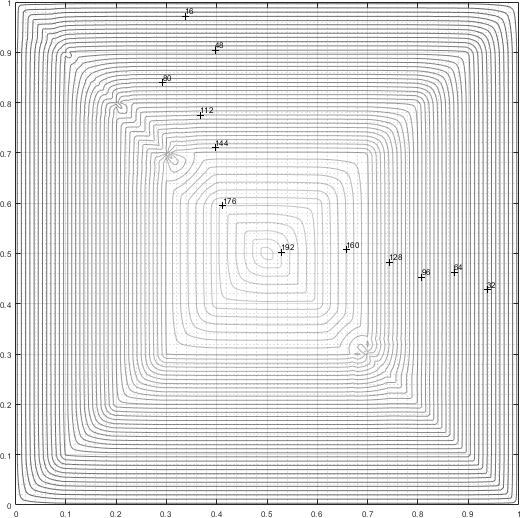}\\
	    (a)&(b)\\
	\includegraphics[width=0.50\textwidth]{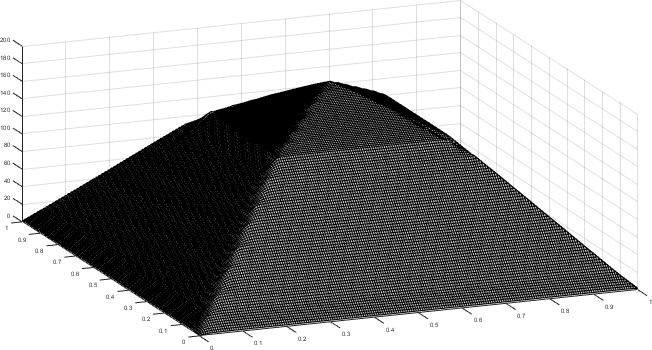}&
	\includegraphics[width=0.25\textwidth]{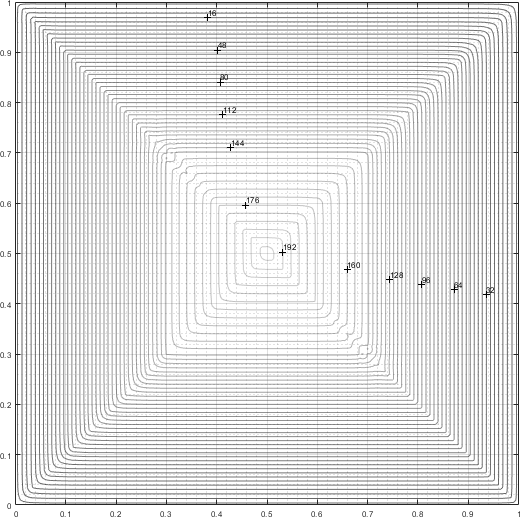}\\
	    (c)&(d)
    \end{array}$}
	\caption{\label{Fig.NumEx7.1.2c} Example \ref{Sec7.1.2}.
	$(a)$ Graph of the AMLE interpolation function of $f_K$ with $K$ the set of $6$-contour lines of
		$f$ displayed in Figure \ref{Fig.NumEx7.1.2a}$(c)$.
		Relative $L^2$-Error $\epsilon=0.01675$. 
	$(b)$ Isolines of the AMLE interpolant at equally spaced heights equal to $(\max(f)-\min(f))/50$.		
	$(c)$ Graph of the AMLE interpolation function of $f_K$ with $K$ the set of $15$-contour lines of
		$f$ displayed in Figure \ref{Fig.NumEx7.1.2a}$(d)$.
		Relative $L^2$-Error $\epsilon= 0.0071297$. 
	$(d)$ Isolines  of the AMLE interpolant at equally spaced heights equal to $(\max(f)-\min(f))/50$.
	}
\end{figure}

We describe now the approximation of the continuous piecewise affine function $f$ associated with the triangulation shown 
in Figure \ref{Fig.NumEx7.1.2a}$(a)$ where also the node values of $f$ are given while Figure \ref{Fig.NumEx7.1.2a}$(b)$
displays the graph of $f$. Two different sample sets of contour lines 
have been considered. One consists of $6$ isolines whereas the other one is formed by $15$ isolines.
The isolines are not equally spaced and are displayed in Figure \ref{Fig.NumEx7.1.2a}$(c)$
and Figure \ref{Fig.NumEx7.1.2a}$(d)$, respectively, whereas 
the graphs of the corresponding average approximations $A_{\lambda}^M(f_K)$ are shown in Figure \ref{Fig.NumEx7.1.2b}$(a)$ 
and Figure \ref{Fig.NumEx7.1.2b}$(c)$ along with the isolines corresponding to $50$ equally spaced isolevels.
In this example the isolines are not smooth curves so that locally, around their singularities, 
for the interpretation of the results, it can be useful to recall and compare with the behaviour of the average approximation
described in the Prototype Example \ref{Ex6.5}$(ii)$ in Section \ref{Sec.ProtEx}. 
The average approximation displays a step which reduces by increasing the 
number of isolines. Note that these steps are also visible in the Matlab display of the graph of the function $f$,
thus they are  errors of the interpolation scheme that is used. 
We find that for the reconstruction of the function sampled on the $6$-contour line set, 
the relative $L^2$-error $\epsilon$ is equal to $0.019302$. This value reduces to $0.004805$ for the reconstruction of the 
function sampled on the $15$-contour line set $K$. For both these two examples,  
it is confirmed that the average approximation $A_\lambda^M(f_K)$ interpolates $f_K$ given that the computed 
value of $\epsilon_K$ is of the order of $10^{-16}$.  

The AMLE method appears yielding 
slightly better results for the reconstruction from the 
sample set $K$ of $6$ contour lines. In this case, we find a relative $L^2$-error $\epsilon$ equal to $0.01675$, slightly 
lower than the one produced by $A_{\lambda}^M(f_K)$. Figure \ref{Fig.NumEx7.1.2c}$(a)$ displays the graph of the AMLE 
interpolant which does not contain steps along the edges of the pyramid, whereas Figure \ref{Fig.NumEx7.1.2c}$(b)$ shows its 
isolines for $50$ levels of equally spaced heights. For the AMLE interpolant of the 
sample set $K$ of $15$ contour lines, whose graph is displayed in Figure \ref{Fig.NumEx7.1.2c}$(c)$
and the isolines in Figure \ref{Fig.NumEx7.1.2c}$(d)$, the relative $L^2$-error $\epsilon$ is equal to $0.00713$,
which is slightly higher than the one produced by $A_{\lambda}^M(f_K)$ for the same sample set $K$. Note also here
the appearance of additional kinks in the graph of the AMLE interpolant which might reduce the global 
quality of the AMLE approximation compared to $A_{\lambda}^M(f_K)$.

%%%%%%%%%%%%%%%%%%%%%%%%%%%%%%%%%%%%%%%%%%%%%%%%%%%%%%

\subsubsection{Discontinuous piecewise affine function}\label{Sec7.1.3}

%%%%%%%%%%%%%%%%%%%%%%%%%%%%%%%%%%%%%%%%%%%%%%%%%%%%%%

\begin{figure}[htbp]
  \centerline{$\begin{array}{cc}
		\includegraphics[width=0.35\textwidth]{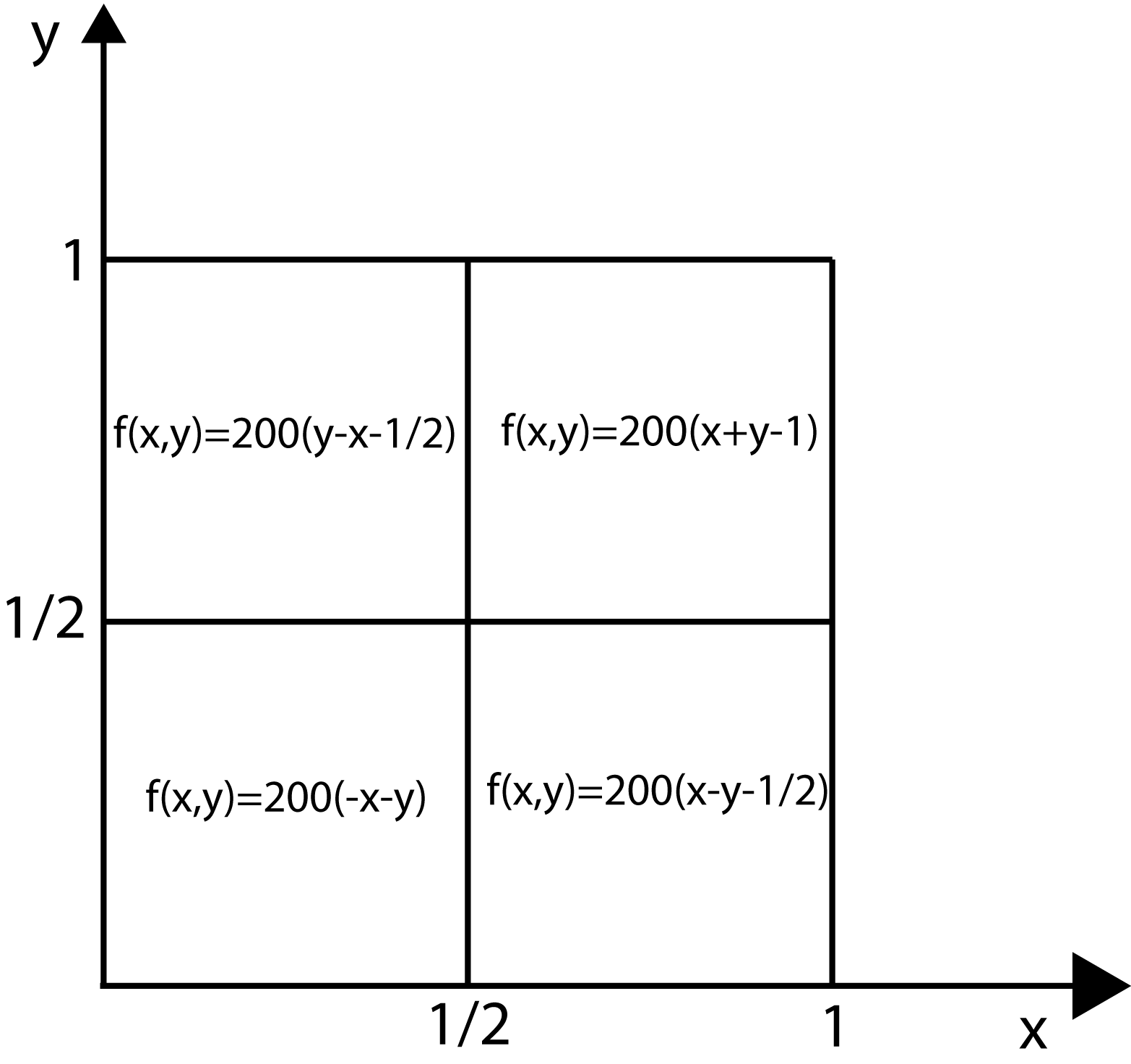}
		&\includegraphics[width=0.50\textwidth]{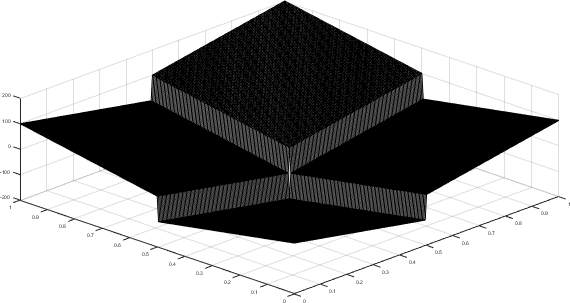}\\
		(a)&(b)\\
		\includegraphics[width=0.25\textwidth]{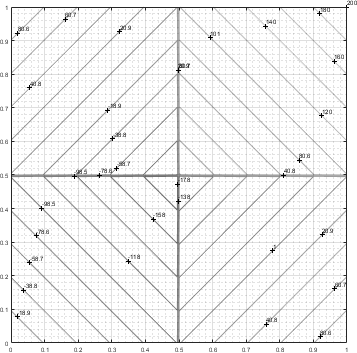}&
		\includegraphics[width=0.25\textwidth]{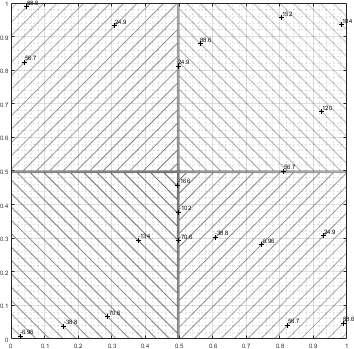}\\
		(c)&(d)
    	      \end{array}$}
	\caption{\label{Fig.NumEx7.1.3a}  Example \ref{Sec7.1.3}.
	$(a)$ Equations of each affine part of $f$. 
	$(b)$ Graph of $f$.
	$(c)$ Sample set $K$ of $20$-contour line of $f$ at equally spaced heights equal to $(\max(f)-\min(f))/20$, 
		defining the sample function $f_K$. 	
	$(d)$ Sample set $K$ of $50$-contour lines of $f$at equally spaced heights equal to $(\max(f)-\min(f))/50$, 
		defining the sample function $f_K$;
	}
\end{figure}

\begin{figure}[htbp]
	\centerline{$\begin{array}{cc}
	\includegraphics[width=0.50\textwidth]{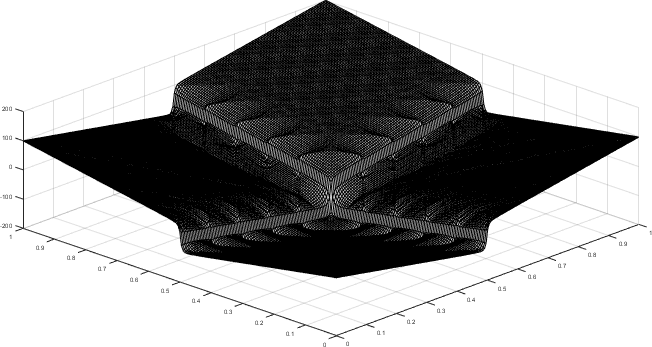}&
	\includegraphics[width=0.25\textwidth]{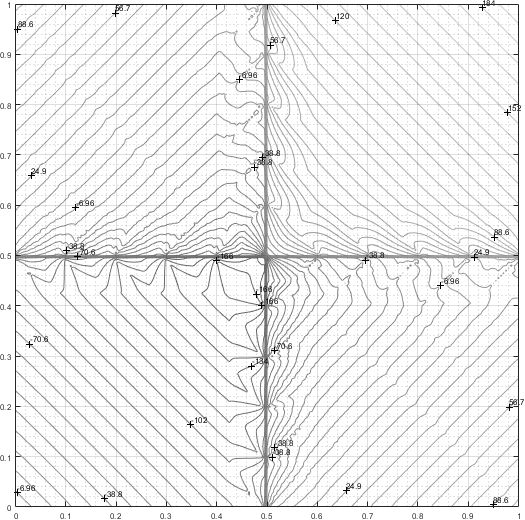}\\
	    (a)&(b)\\
	\includegraphics[width=0.50\textwidth]{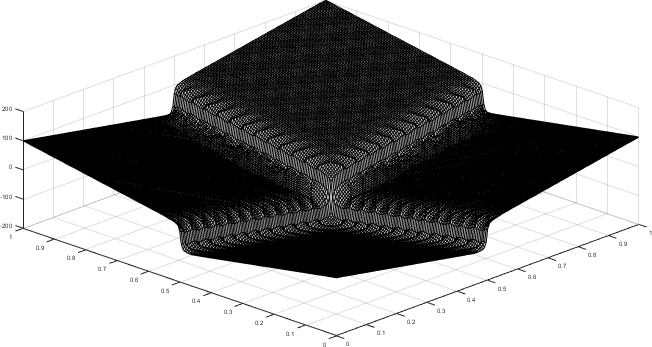}&
	\includegraphics[width=0.25\textwidth]{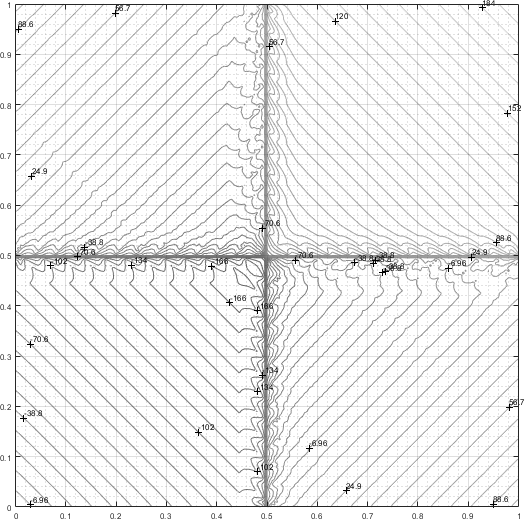}\\
	    (c)&(d)
		\end{array}
		$}
	\caption{\label{Fig.NumEx7.1.3b}  Example \ref{Sec7.1.3}.
	$(a)$ Graph of the AMLE interpolation function of $f_K$ with $K$ the set of $20$-contour lines of
		$f$ displayed in Figure \ref{Fig.NumEx7.1.3a}$(c)$.
		Relative $L^2$-Error $\epsilon=0.1071 $.
	$(b)$ Isolines of the AMLE interpolant at equally spaced heights equal to $(\max(f)-\min(f))/50$.		
	$(c)$ Graph of the AMLE interpolation function of $f_K$ with $K$ the set of $50$-contour lines of
		$f$ displayed in Figure \ref{Fig.NumEx7.1.3a}$(d)$.
		Relative $L^2$-Error $\epsilon=0.06738$.
	$(d)$ Isolines of the AMLE interpolant  at equally spaced heights equal to $(\max(f)-\min(f))/50$.		
	}
\end{figure}

The approximation of discontinuous functions has not been covered by the theoretical developments of
Section \ref{Sec.LSAprx}, where we assumed $f$ to be continuous. Now we present
a test case where we examine how our average approximation performs 
numerically and verify that also in this case $A_{\lambda}^M(f_K)$ represents a continuous 
interpolation of $f_K$.
We consider the following discontinuous piecewise affine function
\begin{equation*}
	f:(x,y)\in]0,\,1[^2\to 200, \;\;\; f(x,y) = \left\{\begin{array}{lll}
		\displaystyle x+y-1,& \displaystyle\text{if }1/2\leq x\leq 1, &\displaystyle 1/2 \leq y\leq 1 \\[1.5ex]
		\displaystyle x-y-1/2& \displaystyle\text{if }1/2\leq x\leq 1, &\displaystyle 0 \leq y< 1/2 \\[1.5ex]
		\displaystyle -x+y-1/2& \displaystyle\text{if }0\leq x< 1/2, &\displaystyle 1/2 \leq y\leq 1 \\[1.5ex]
		\displaystyle -x-y& \displaystyle\text{if }0\leq x<1/2, &\displaystyle 0 \leq y< 1/2 \\[1.5ex]
	\end{array}\right.
\end{equation*}
whose graph is displayed in Figure \ref{Fig.NumEx7.1.3a}$(b)$ while Figure \ref{Fig.NumEx7.1.3a}$(a)$
shows the equation of $f$ in each of its affine parts.

We compare the reconstruction of $f$ for two sample sets $K$, one  formed by $20$ equally spaced isolines
and the other by $100$  equally spaced isolines. Such sets are displayed in Figure \ref{Fig.NumEx7.1.3a}$(c)$
and Figure  \ref{Fig.NumEx7.1.3a}$(d)$, respectively. Notably, for both sample sets $K$,
$A_{\lambda}^M(f_K)$ coincides exactly with the original function $f$. We find, indeed, for 
both sample sets $K$, $\epsilon$ and $\epsilon_K$ the order $10^{-15}$ by taking 
$\lambda=10^7$, $M=10^6$. This occurs because of an exact sampling of 
the discontinuity jump, thus we are able to reproduce exactly the affine parts of $f$, consistently with the 
theoretical findings of Section \ref{Sec.LSAprx}. Furthermore, given the high value of $\lambda$ and recalling the 
behaviour of the jump in the Prototype Example \ref{Ex.AprxJump}, we are able to describe the sharp discontinuity.

For the case where we do not have an exact sampling of the discontinuity jump,  
we refer to Example \ref{Sec7.2.3} concerning the surface
reconstruction from point clouds with sampling points not necessarily on the discontinuity.

A different behaviour is displayed by the AMLE interpolation. Consistently with the observations in 
\cite{LMS13}, the level lines of the AMLE interpolant are smooth \cite{Sav05}, thus discontinuities  cannot be 
recovered. A better visual appreciation of this fact is obtained  by looking at the graphs of the AMLE interpolant shown in  
Figure \ref{Fig.NumEx7.1.3b}$(a)$ and Figure \ref{Fig.NumEx7.1.3b}$(c)$ for the two sample sets $K$,
and at their isolines displayed in Figure \ref{Fig.NumEx7.1.3b}$(b)$ and Figure \ref{Fig.NumEx7.1.3b}$(d)$, respectively.
The isolines at the two sides of the jump should `end' in the discontinuity but they 
are somehow enforced to join each other by the continuous isolines of the AMLE interpolant. In this case
we find values of the relative $L^2$-error $\epsilon$, with $\epsilon=0.1071$
and $\epsilon=0.06738$ for the two sample sets, respectively. 

Table \ref{Tab:CntLin} summarizes the relative $L^2$-errors of $A_{\lambda}^M(f_K)$
and the AMLE interpolant for the examples considered in this section.

%%%%%%%%%%%%%%%%%%%%%%%%%%%%%%%%%%%%%%%%%%%%%%%%%%%%%%%%%%%%%%%%%%%%%%%%%%%%%%%%%%%%%%%%%%%%%%%%%%%%%%%%%%%%%%%

\begin{table}[tbhp]
\centerline{
\begin{tabular}{cc|c|c|}
\cline{3-4}
			    &     & \multicolumn{2}{c|}{ $\epsilon$} \\ \hline
\multicolumn{1}{ |c| }{$f$} & $K$ & $A_{\lambda}^M(f_K)$ &  AMLE      \\ \hline
\multicolumn{1}{ |c }{\multirow{2}{*}{F}} & \multicolumn{1}{ |c| }{$10$ level lines}    & $0.0199$		& $0.0338$	 \\ \cline{2-4}
\multicolumn{1}{ |c }{} & \multicolumn{1}{ |c| }{$50$ level lines}			& $0.0021$		& $0.0101$	 \\ \hline
\multicolumn{1}{ |c }{\multirow{2}{*}{CPA}} & \multicolumn{1}{ |c| }{$6$ level lines}	& $0.0193$		& $0.0167$	 \\ \cline{2-4}
\multicolumn{1}{ |c }{} & \multicolumn{1}{ |c| }{$15$ level lines}			&  $0.0048$ 		& $0.0071$	 \\ \hline
\multicolumn{1}{ |c }{\multirow{2}{*}{DPA}} & \multicolumn{1}{ |c| }{$20$ level lines}	& $8.7\cdot 10^{-15}$	& $0.1071$	 \\ \cline{2-4}
\multicolumn{1}{ |c }{} & \multicolumn{1}{ |c| }{$100$ level lines}			& $1.5\cdot 10^{-16}$	& $0.0674$	 \\ \hline
\end {tabular}
}
\caption{\label{Tab:CntLin} Summary of the accuracy of the compensated convexity based interpolant $A_{\lambda}^M(f_K)$ 
and of the AMLE interpolant for the Examples considered in Section \ref{Sec7.1}.
Legenda: K Sample set. $\epsilon$ Relative $L^2$-error. $\epsilon_K$ Relative $L^2$-error on the sample set $K$. 
	F Franke test function (Example~\ref{Sec7.1.1}). 
	CPA Continuous piecewise affine function (Example~\ref{Sec7.1.2}). 
	 DPA Discontinuous piecewise affine function (Example~\ref{Sec7.1.3}).
	 }
\end{table}

%%%%%%%%%%%%%%%%%%%%%%%%%%%%%%%%%%%%%%%%%%%%%%%%%%%%%%%%%%%%%%%%%%%%%%%%%%%%%%%%%%%%%%%%%%%%%%%%%%%%%%%%%%%%%%%

\subsection{Scattered data approximation}\label{Sec7.2}

%%%%%%%%%%%%%%%%%%%%%%%%%%%%%%%%%%%%%%%%%%%%%%%%%%%%%%%%%%%%%%%%%%%%%%%%%%%%%%%%%%%%%%%%%%%%%%%%%%%%%%%%%%%%%%%
\begin{figure}[htbp]
  \centerline{$\begin{array}{cc}
	\includegraphics[width=0.25\textwidth]{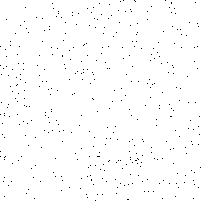}&
	\includegraphics[width=0.25\textwidth]{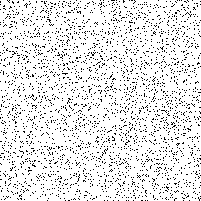}\\
	(a) & (b)	
    \end{array}$}
   \caption{\label{Fig.NumEx7.2} 
	Set $K$ of sample points of a grid of $201\times 201$ points in $]0,\,1[^2$ for two levels of 
	sampling density:
	$(a)$ Coarse sampling with $400$ grid points out of $40401$. 
	$(b)$ Dense sampling with $4061$ grid points out of $40401$.
	}
\end{figure}

We turn now to some numerical experiments on scattered data approximation. In particular, in the terminology of 
\cite{LF99}, we consider the problem of function reconstruction from point clouds, where 
the sample points that form the set $K$ do not meet any particular condition as to spacing or density.
As in the previous section, the set of test problems consists of three test functions with different regularity:
an infinitely differentiable function given by the Franke test function, a continuous piecewise affine function and 
a discontinuous piecewise affine function. The three test functions are all to be approximated 
in $\Omega=]0,\,1[^2$. In the numerical implementation of the method, the domain $\Omega$ is discretized with 
a grid of $201\times201$ points and the two sample sets $K$ are obtained by sampling the grid points using a 
random number generator with different levels of density. The two sample sets $K$, corresponding to a coarse and 
a dense sampling, are displayed in Figure \ref{Fig.NumEx7.2}$(a)$ and Figure \ref{Fig.NumEx7.2}$(b)$, respectively.
The reason for taking such a regular discretization of $\Omega$ is because the numerical 
scheme we use to compute the convex envelope 
(see Algorithm \ref{Algo:CnvxEnv}), is particularly suitable for applications to image processing
where such discrete geometry is related to the image resolution.

For the measure of the global quality of the approximation $A_{\lambda}^M(f_K)$ we compute the  relative $L^2$-error
$\epsilon$ defined by Eq. \eqref{Eq.RelErr} whereas we will use the relative $L^2$-error $\epsilon_K$ 
defined  by

\begin{equation}\label{Eq.RelErrResSc}
	\epsilon_K=\displaystyle \frac{\sqrt{\displaystyle \sum_{k\in K} |f(x_k)-A_{\lambda}^M(f_K)(x_k)|^2}}{\sqrt{\displaystyle\sum_{k\in K}|f(x_k)|^2}}
\end{equation}

\noindent to assess the quality of $A_{\lambda}^M(f_K)$ as interpolant of $f_K$. In this case too, we will find that 
the average approximation $A_{\lambda}^M(f_K)$ is an interpolation of $f_K$, consistently with the 
theoretical findings of Section \ref{Sec.AprxSctD}. We then conclude this section by giving an example of
digital elevation model reconstruction starting from real data, and another of
salt \& pepper noise removal as an application of scattered data approximation to image processing.

%%%%%%%%%%%%%%%%%%%%%%%%%%%%%%%%%%%%%%%%%%%%%%%%%%%%%%

\subsubsection{Franke test function}\label{Sec7.2.1}

%%%%%%%%%%%%%%%%%%%%%%%%%%%%%%%%%%%%%%%%%%%%%%%%%%%%%%
\begin{figure}[htbp]
  \centerline{$\begin{array}{cc}
	\includegraphics[width=0.55\textwidth]{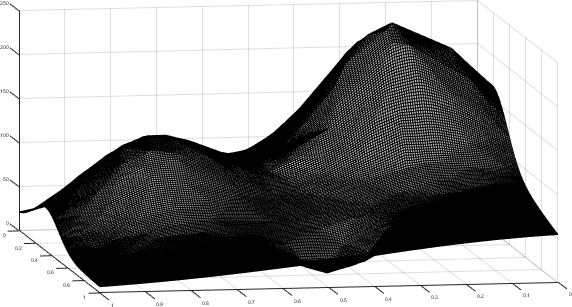}&
	\includegraphics[width=0.25\textwidth]{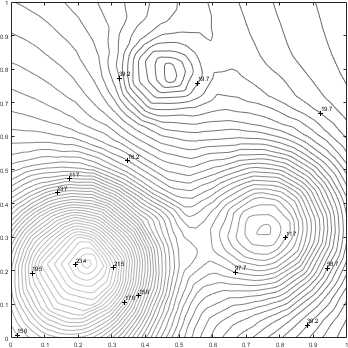}\\
	(a) & (b)	\\
	\includegraphics[width=0.55\textwidth]{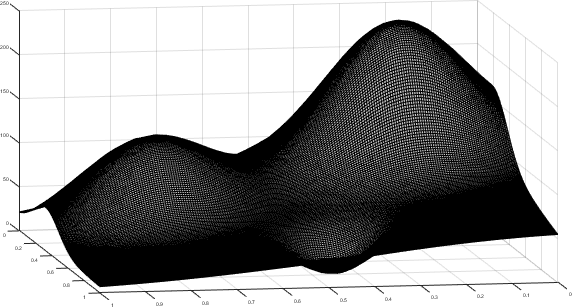}&
	\includegraphics[width=0.25\textwidth]{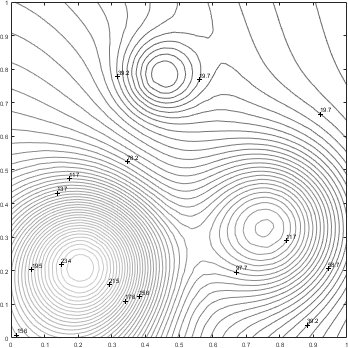}\\
	(c) & (d)	
    \end{array}$
    }
   \caption{\label{Fig.NumEx7.2.1a} Example \ref{Sec7.2.1}. 
	 $(a)$ Graph of $A_{\lambda}^M(f_K)$ for $\lambda=1\cdot 10^4$, $M=1\cdot 10^5$ and the set $K$ of Figure \ref{Fig.NumEx7.2}$(a)$. 
	 Relative $L^2$-Errors: $\epsilon=0.020252$, $\epsilon_K= 5.31\cdot 10^{-15}$.
	 $(b)$ Isolines of  $A_{\lambda}^M(f_K)$ at equally spaced heights equal to $(\max(f)-\min(f))/50$.
	 $(c)$ Graph of $A_{\lambda}^M(f_K)$ for $\lambda=5\cdot 10^3$, $M=1\cdot 10^5$  and the set $K$ of Figure \ref{Fig.NumEx7.2}$(b)$ 
	 Relative $L^2$-Errors: $\epsilon=0.0015548$, $\epsilon_K= 4.13\cdot 10^{-15}$.
	 $(d)$ Isolines of  $A_{\lambda}^M(f_K)$ at equally spaced heights equal to $(\max(f)-\min(f))/50$.
	}
\end{figure}

\begin{figure}[htbp]
  \centerline{$\begin{array}{cc}
	\includegraphics[width=0.50\textwidth]{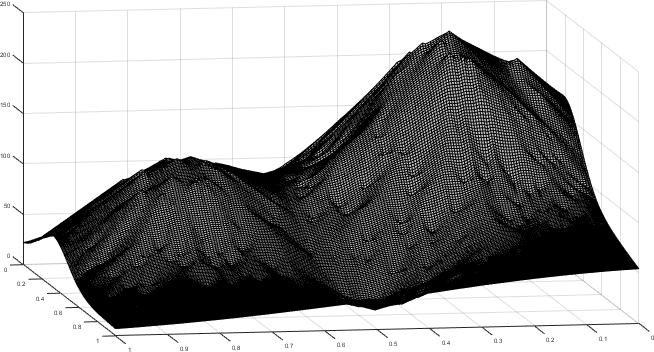}&
	\includegraphics[width=0.25\textwidth]{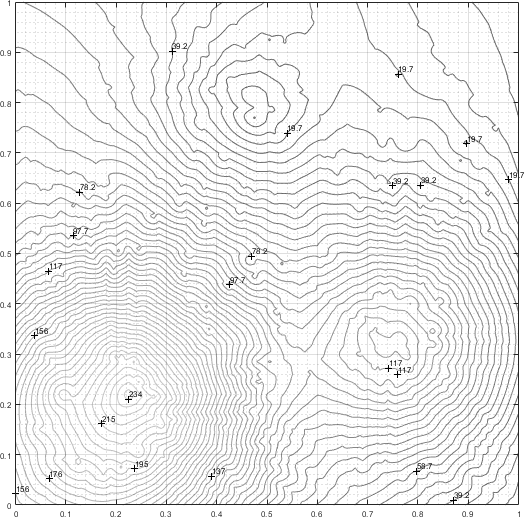}\\
	(a) & (b)	\\
	\includegraphics[width=0.50\textwidth]{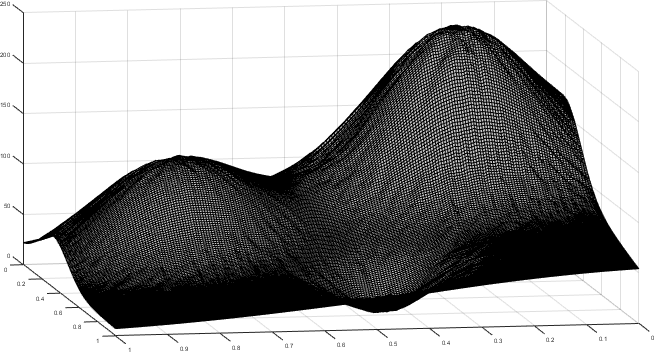}&
	\includegraphics[width=0.25\textwidth]{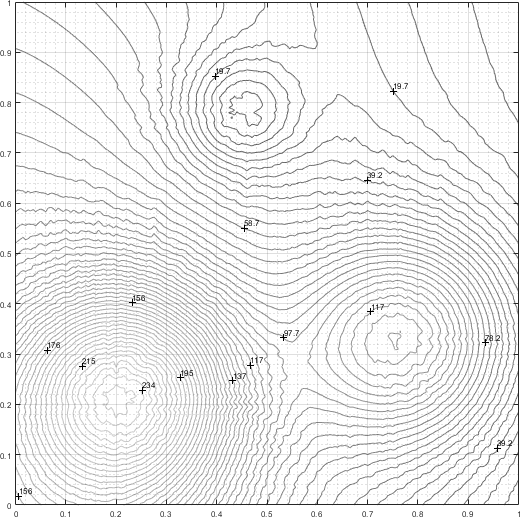}\\
	(c) & (d)	
    \end{array}$
    }
   \caption{\label{Fig.NumEx7.2.1b} Example \ref{Sec7.2.1}.
	 $(a)$ Graph of the AMLE interpolation function of $f_K$ with $K$ the set of scattered points 
		displayed in Figure \ref{Fig.NumEx7.2}$(a)$. 
		Relative $L^2$-Error: $\epsilon=0.05764$.
	 $(b)$ Isolines of the AMLE interpolant at equally spaced heights equal to $(\max(f)-\min(f))/50$.
	 $(c)$  Graph of the AMLE interpolation function of $f_K$ with $K$ the set of scattered points 
		displayed in Figure \ref{Fig.NumEx7.2}$(b)$. 
		 Relative $L^2$-Error: $\epsilon=0.010902$.
	 $(d)$ Isolines of the AMLE interpolant at equally spaced heights equal to $(\max(f)-\min(f))/50$.
	}
\end{figure}

In this example, the Franke test function $f$ defined by Eq. \eqref{Eq.FrankeFnct} is sampled over the two sets 
$K$ of scattered points displayed in Figure \ref{Fig.NumEx7.2}$(a)$ and Figure \ref{Fig.NumEx7.2}$(b)$, respectively.
For the resulting sample functions $f_K$ we compute the corresponding average approximations $A_{\lambda}^M(f_K)$
whose graphs are displayed in Figure \ref{Fig.NumEx7.2.1a}, along with the respective isolines. Specifically,
the comparison of the isolines of $A_{\lambda}^M(f_K)$ displayed in Figure \ref{Fig.NumEx7.2.1a}$(b)$ and 
in Figure \ref{Fig.NumEx7.2.1a}$(d)$ for the coarse and dense sample sets $K$, respectively,
with the isolines of the Franke function $f$ displayed in Figure \ref{Fig.NumEx7.1.1a}$(d)$, allows a visual 
appreciation of the quality of the reconstruction. This is also confirmed by the computed values of the relative $L^2$-error 
$\epsilon$. For the coarse sample set we get 
$\epsilon=0.0206$ whereas,  for the denser sample set, $\epsilon=0.00157$. Finally, also in this case, we verify that 
$A_{\lambda}^M(f_K)$ is an interpolant of $f_K$ given that for both approximations the relative $L^2$-error $\epsilon_K$
defined by Eq. \eqref{Eq.RelErrResSc} is of the order of $10^{-15}$.

The AMLE method as introduced in \cite{CMS98} can be applied also in this case for the interpolation 
of isolated points. In fact, this is one of its particular feature out of the pde based interpolators. The graphs
of the AMLE interpolants for the two sample sets are displayed in Figure \ref{Fig.NumEx7.2.1b},
which contains also the plot of the corresponding isolines for $50$ level lines of equally spaced heights. 
The plot of these isolines, once compared with the same isolines of $f$ displayed in 
Figure~\ref{Fig.NumEx7.1.1a}$(c)$, allows a visual 
assessment of the quality of the reconstruction. As in the Example \ref{Sec7.1.1} concerning with the 
reconstruction from contour lines, we note also here the introduction of artificial artefacts in the form of krinks 
in the graph of the interpolant, which, in contrast, are not present in the graph of $A_{\lambda}^M(f_K)$. 
For the coarse and dense sampling set we find that the relative $L^2$-error of the AMLE interpolant 
amounts to $\epsilon=0.05764$ and $\epsilon=0.010902$, respectively, which are slightly higher than the values 
produced by $A_{\lambda}^M(f_K)$.

%%%%%%%%%%%%%%%%%%%%%%%%%%%%%%%%%%%%%%%%%%%%%%%%%%%%%%

\subsubsection{Continuous piecewise affine function}\label{Sec7.2.2}

%%%%%%%%%%%%%%%%%%%%%%%%%%%%%%%%%%%%%%%%%%%%%%%%%%%%%%
\begin{figure}[htbp]
  \centerline{$\begin{array}{cc}
	\includegraphics[width=0.5\textwidth]{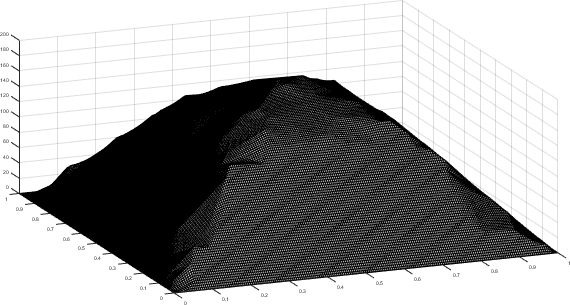}
	&\includegraphics[width=0.25\textwidth]{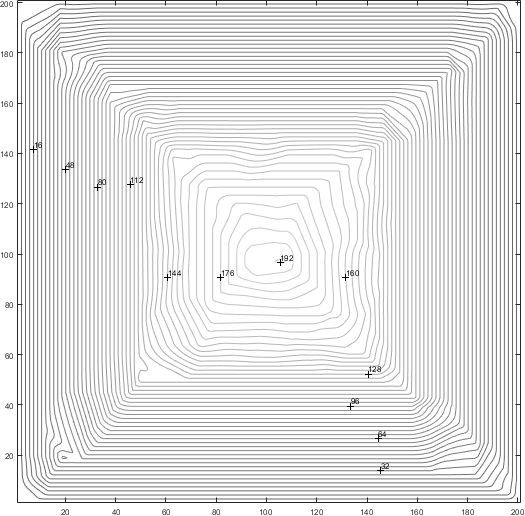}\\
	(a) & (b)	\\
	\includegraphics[width=0.5\textwidth]{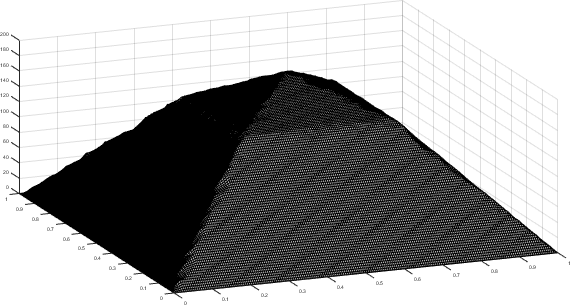}
	&\includegraphics[width=0.25\textwidth]{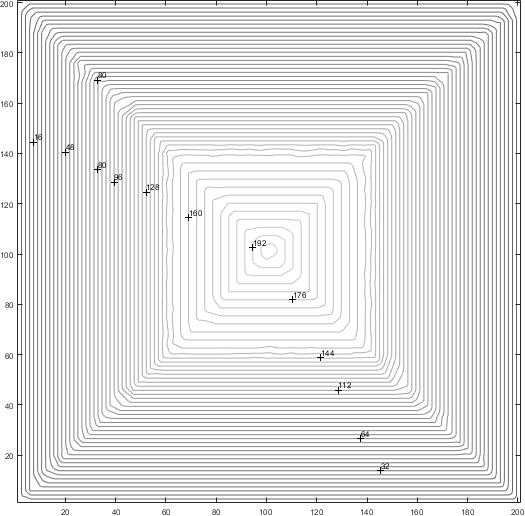}\\
	(c) & (d)	
    \end{array}$
}
   \caption{\label{Fig.NumEx7.2.2a} Example \ref{Sec7.2.2}.
	 $(a)$ Graph of $A_{\lambda}^M(f_K)$ for $\lambda=5\cdot 10^4$, $M=1\cdot 10^5$ and the set $K$ of Figure \ref{Fig.NumEx7.2}$(a)$. 
	 Relative $L^2$-Errors: $\epsilon=0.021574$, $\epsilon_K= 4.4626\cdot 10^{-16}$.
	 $(b)$ Isolines of  $A_{\lambda}^M(f_K)$ at equally spaced heights equal to $(\max(f)-\min(f))/50$.
	 $(c)$ Graph of $A_{\lambda}^M(f_K)$ for $\lambda=5\cdot 10^4$, $M=1\cdot 10^5$  and the set $K$ of Figure \ref{Fig.NumEx7.2}$(b)$. 
	 Relative $L^2$-Errors: $\epsilon=0.003914$, $\epsilon_K=6.2983\cdot 10^{-16}$.
	 $(d)$ Isolines of $A_{\lambda}^M(f_K)$ at equally spaced heights equal to $(\max(f)-\min(f))/50$.
	}
\end{figure}

\begin{figure}[htbp]
  \centerline{$\begin{array}{cc}
	\includegraphics[width=0.5\textwidth]{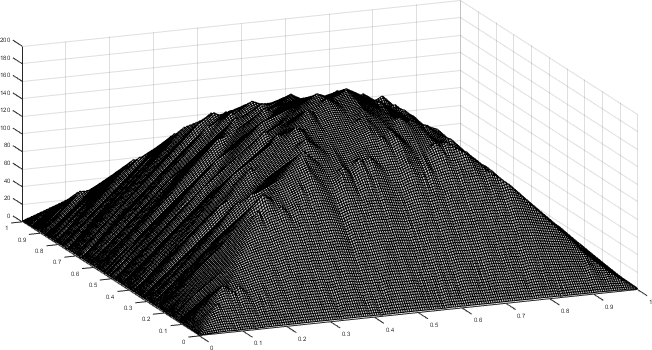}
	&\includegraphics[width=0.25\textwidth]{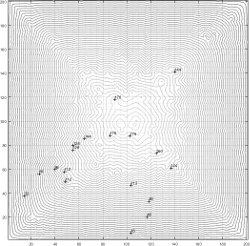}\\
	(a) & (b)	\\
	\includegraphics[width=0.5\textwidth]{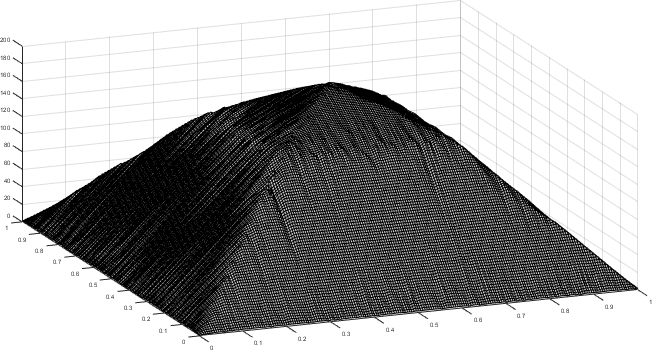}
	&\includegraphics[width=0.25\textwidth]{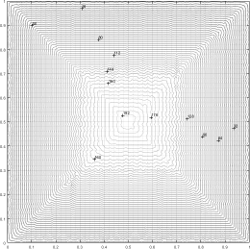}\\
	(c) & (d)	
    \end{array}$
}
   \caption{\label{Fig.NumEx7.2.2b} Example \ref{Sec7.2.2}.
	 $(a)$ Graph of the AMLE interpolation function of $f_K$ with $K$ the set of scattered points 
		displayed in Figure \ref{Fig.NumEx7.2}$(a)$.  
		Relative $L^2$-Error: $\epsilon=0.053594$.
	 $(b)$ Isolines of the AMLE interpolant at equally spaced heights equal to $(\max(f)-\min(f))/50$.
	 $(c)$ Graph of the AMLE interpolation function of $f_K$ with $K$ the set of scattered points 
		displayed in Figure \ref{Fig.NumEx7.2}$(b)$.  
		Relative $L^2$-Error: $\epsilon=0.012515$.
	 $(d)$ Isolines of the AMLE interpolant at equally spaced heights equal to $(\max(f)-\min(f))/50$.
	}
\end{figure}

The continuous piecewise affine function $f$ introduced in Section \ref{Sec7.1.2} is evaluated here over the two sample 
sets $K$ of Figure \ref{Fig.NumEx7.2}$(a)$ and Figure \ref{Fig.NumEx7.2}$(b)$, defining two test cases of sample 
function $f_K$. The graph of the corresponding average approximation $A_{\lambda}^M(f_K)$ is displayed in   
Figure \ref{Fig.NumEx7.2.2a} along with the respectives isolines whereas Figure \ref{Fig.NumEx7.2.2b}
shows those of the AMLE interpolating along with its isolines of equally spaced heights. The drawing of the isolines allows a visual 
assessment of the quality of the reconstruction if these are compared to the isolines of the original function 
$f$ displayed in Figure \ref{Fig.NumEx7.1.2b}$(c)$. 
A first observation about the graphs of $A_{\lambda}^M(f_K)$ is the nearly 
absence of the steps along the edges of the pyramid due to the constraint enforced by the fixed contour lines, 
on the contrary  the graphs of the AMLE  interpolant
present, even for this example, artefacts in the form of artificial krinks and valleys.
The relative $L^2$-error $\epsilon$ produced by  $A_{\lambda}^M(f_K)$
is equal to $0.0215$
for the coarse sample set and  to $0.00390$ for the denser sample set, whereas it is 
 $\epsilon=0.053594$ and  $\epsilon=0.012515$ for the AMLE interpolant of the coarse and dense sample set, respectively.
Compared with the reconstruction of $f$ from contour lines, where the sample points can be considered to be somehow organized, we 
observe that both the reconstructed function $A_{\lambda}^M(f_K)$ 
and the AMLE interpolant appear to be less regular, which reflects the fact that 
the sample points are scattered over $\Omega$ without any requirement of spacing or density. This effect clearly 
reduces by increasing the density of the sample points, though for the AMLE interpolant 
we note that the relative $L^2$-errors for the 
two cases of sampling density remains of the same order of magnitude. For  this example too, we finally verify that 
$A_{\lambda}^M(f_K)$ is an interpolation of $f_K$ given that the relative $L^2$-error $\epsilon_K$ is of the order
$10^{-16}$ for both the two test cases.

%%%%%%%%%%%%%%%%%%%%%%%%%%%%%%%%%%%%%%%%%%%%%%%%%%%%%%

\subsubsection{Discontinuous piecewise affine function}\label{Sec7.2.3}

%%%%%%%%%%%%%%%%%%%%%%%%%%%%%%%%%%%%%%%%%%%%%%%%%%%%%%

\begin{figure}[htbp]
  \centerline{$\begin{array}{cc}
	\includegraphics[width=0.55\textwidth]{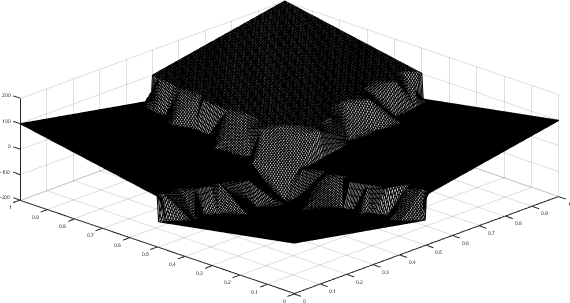}
	&\includegraphics[width=0.25\textwidth]{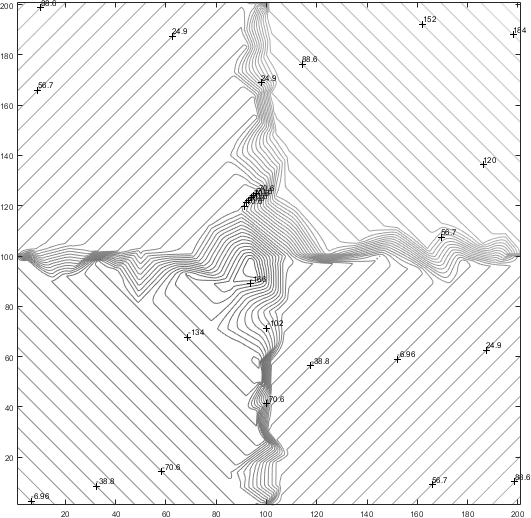}\\
	(a) & (b)	\\
	\includegraphics[width=0.55\textwidth]{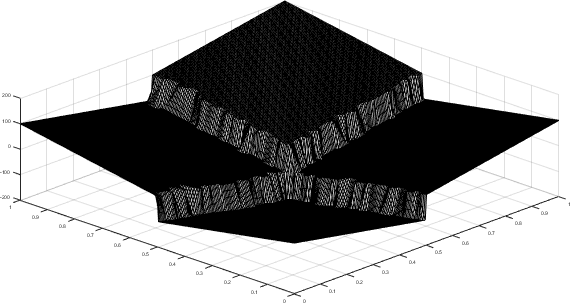}
	&\includegraphics[width=0.25\textwidth]{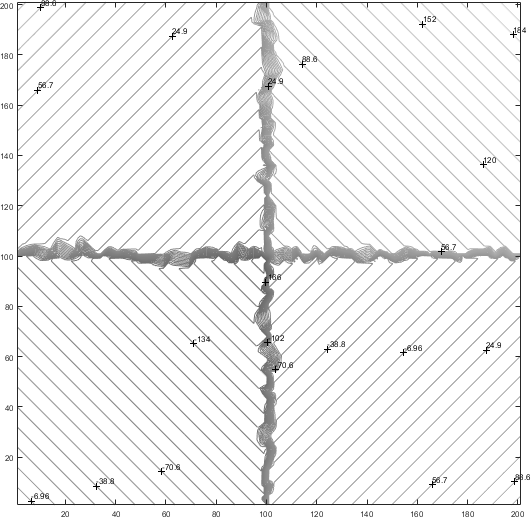}\\
	(c) & (d)	
    \end{array}$
}
   \caption{\label{Fig.NumEx7.2.3a} Example \ref{Sec7.2.3}.
	 $(a)$ Graph of $A_{\lambda}^M(f_K)$ for $\lambda=1\cdot 10^7$, $M=1\cdot 10^5$ and the set $K$ of Figure \ref{Fig.NumEx7.2}$(a)$. 
	 Relative $L^2$-Errors: $\epsilon=0.16729$, $\epsilon_K=1.2849\cdot 10^{-16}$.
	 $(b)$ Isolines of  $A_{\lambda}^M(f_K)$ at equally spaced heights equal to $(\max(f)-\min(f))/50$.
	 $(c)$ Graph of $A_{\lambda}^M(f_K)$ for $\lambda=1\cdot 10^7$, $M=1\cdot 10^5$  and the set $K$ of Figure \ref{Fig.NumEx7.2}$(b)$. 
	 Relative $L^2$-Errors: $\epsilon=0.088589$, $\epsilon_K=1.459\cdot 10^{-16}$.
	 $(d)$ Isolines of  $A_{\lambda}^M(f_K)$ at equally spaced heights equal to $(\max(f)-\min(f))/50$.
	}
\end{figure}

\begin{figure}[htbp]
  \centerline{$\begin{array}{cc}
	\includegraphics[width=0.5\textwidth]{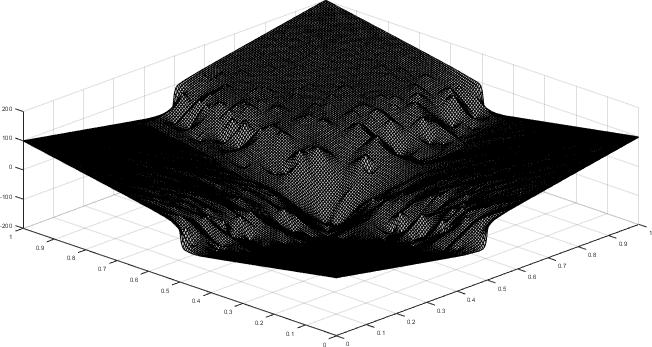}
	&\includegraphics[width=0.25\textwidth]{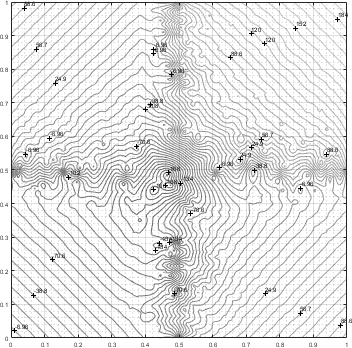}\\
	(a) & (b)	\\
	\includegraphics[width=0.5\textwidth]{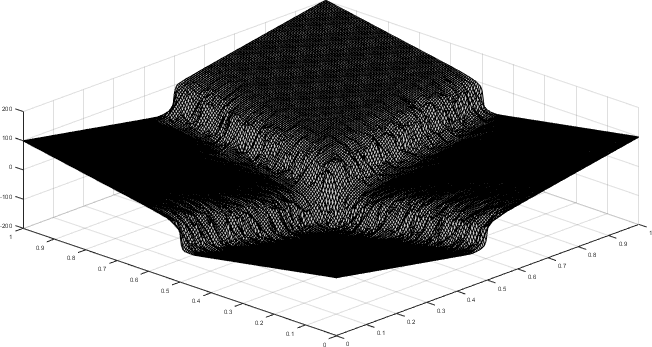}
	&\includegraphics[width=0.25\textwidth]{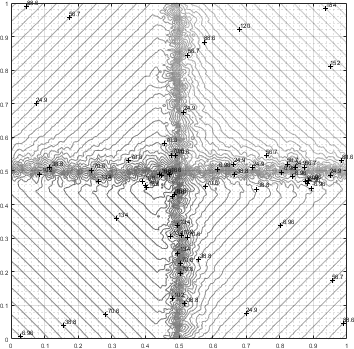}\\
	(c) & (d)	
    \end{array}$
}
   \caption{\label{Fig.NumEx7.2.3b} Example \ref{Sec7.2.3}.
	 $(a)$ Graph of the AMLE interpolation function of $f_K$ with $K$ the set of scattered points 
		displayed in Figure \ref{Fig.NumEx7.2}$(a)$.   
	 Relative $L^2$-Error: $\epsilon=0.22577$.
	 $(b)$ Isolines of the AMLE interpolant at equally spaced heights equal to $(\max(f)-\min(f))/50$.
	 $(c)$ Graph of the AMLE interpolation function of $f_K$ with $K$ the set of scattered points 
		displayed in Figure \ref{Fig.NumEx7.2}$(b)$.   
	 Relative $L^2$-Error: $\epsilon=0.13897$.
	 $(d)$ Isolines of the AMLE interpolant at equally spaced heights equal to $(\max(f)-\min(f))/50$.
	}
\end{figure}

The discontinuous piecewise affine function $f$ introduced in Section \ref{Sec7.1.3} is evaluated here over the two sample 
sets $K$ displayed in Figure \ref{Fig.NumEx7.2}$(a)$ and Figure \ref{Fig.NumEx7.2}$(b)$, to form 
two sample functions $f_K$ corresponding to a coarse and a dense sample set, respectively. The graph of 
$A_{\lambda}^M(f_K)$ is displayed in Figure \ref{Fig.NumEx7.2.3a} for the two cases, along with their isolines,
whereas Figure \ref{Fig.NumEx7.2.3b} shows the graph of the AMLE interpolants
along with their isolines with equally spaced heights. Also here, 
it is useful to compare such isolines with those of the original function $f$ displayed in 
Figure \ref{Fig.NumEx7.1.3b}$(d)$ for a visual assessment of the quality of the reconstructions. 
Unlike the reconstruction of $f$ from contour lines, where we had the exact sampling of the discontinuity which 
was coincident with the grid lines,  here we note an irregular
behaviour for $A_{\lambda}^M(f_K)$
around the discontinuities of $f$. Such irregular behaviour  reduces by increasing the sampling density, especially 
if such density increase occurs in the neighborhood of the singularities. 
On the other hand, the AMLE interpolant displays around the singularities a behaviour similar 
to the one obtained from the 
contour lines, with the difference that now the transition from one affine part of $f$ to the other appears 
to be smoother. As for the accuracy of the reconstructions, for $A_{\lambda}^M(f_K)$ 
we find that $\epsilon=0.173$
for the coarse sample set and $\epsilon=0.0901$ for the denser sample set, whereas 
the relative $L^2$-error
$\epsilon_K$ on both sample sets $K$ is  of the order of $10^{-16}$, 
confirming that again,  $A_{\lambda}^M(f_K)$ is an interpolant of $f_K$. For the AMLE interpolant, 
even in this case, we find higher values for the relative $L^2$-error, with $\epsilon=0.22577$ and $\epsilon=0.13897$ 
for the coarser and denser sample set, respectively. We note also the introduction of artificial artefacts in 
the graph of the AMLE interpolant.

The relative $L^2$-errors obtained for scattered data approximation using $A_{\lambda}^M$
and AMLE interpolation are summarized in Table \ref{Tab:ScaDat} for the examples considered in this section.

%%%%%%%%%%%%%%%%%%%%%%%%%%%%%%%%%%%%%%%%%%%%%%%%%%%%%%

\begin{table}[tbhp]
\centerline{
\begin{tabular}{cc|c|c|c|c|}
\cline{3-4}
			    &     & \multicolumn{2}{c|}{ $\epsilon$} \\ \hline
\multicolumn{1}{ |c| }{$f$} & $K$ & $A_{\lambda}^M(f_K)$ &  AMLE     \\ \hline
\multicolumn{1}{ |c }{\multirow{2}{*}{F}} & \multicolumn{1}{ |c| }{coarse}   & $0.0203$		& $0.0576$	\\ \cline{2-4}
\multicolumn{1}{ |c }{} & \multicolumn{1}{ |c| }{dense}			     & $0.0016$		& $0.0109$	\\ \hline
\multicolumn{1}{ |c }{\multirow{2}{*}{CPA}} & \multicolumn{1}{ |c| }{coarse} & $0.0216$		& $0.0536$	\\ \cline{2-4}
\multicolumn{1}{ |c }{} & \multicolumn{1}{ |c| }{dense}			     & $0.0039$		& $0.0125$	\\ \hline
\multicolumn{1}{ |c }{\multirow{2}{*}{DPA}} & \multicolumn{1}{ |c| }{coarse} & $0.1673$		& $0.2258$	\\ \cline{2-4}
\multicolumn{1}{ |c }{} & \multicolumn{1}{ |c| }{dense}			     & $0.0876$		& $0.1390$	\\ \hline
\end {tabular}
}
\caption{\label{Tab:ScaDat} Accuracy of the interpolation for the examples considered in Section \ref{Sec7.2}.
Legenda: K Sample set. $\epsilon$ Relative $L^2$-error. $\epsilon_K$ Relative $L^2$-error on the sample set $K$. 
	F Franke test function (Example~\ref{Sec7.1.1}). 
	CPA Continuous piecewise affine function (Example~\ref{Sec7.1.2}). 
	 DPA Discontinuous piecewise affine function (Example~\ref{Sec7.1.3}).}
\end{table}

%%%%%%%%%%%%%%%%%%%%%%%%%%%%%%%%%%%%%%%%%%%%%%%%%%%%%%

\subsubsection{DEM Reconstruction}\label{Sec7.2.4}

%%%%%%%%%%%%%%%%%%%%%%%%%%%%%%%%%%%%%%%%%%%%%%%%%%%%%%
\begin{figure}[htbp]
  \centerline{$\begin{array}{ccc}
	\includegraphics[width=0.50\textwidth]{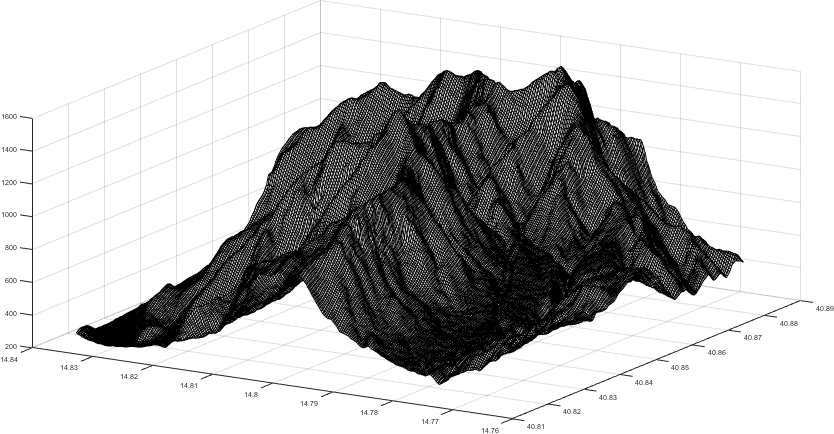}
	&\includegraphics[width=0.225\textwidth]{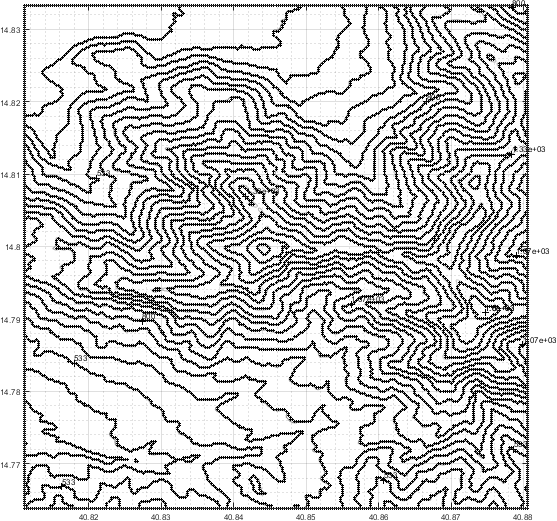}
	&\includegraphics[width=0.225\textwidth]{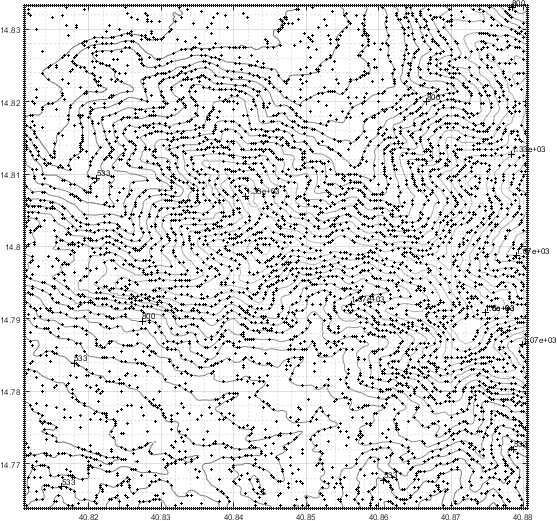}\\
	(a) & (b) & (c)
    \end{array}$
}
   \caption{\label{Fig.NumEx7.2.4a} Example \ref{Sec7.2.4}. Reconstruction of real-world digital elevation maps.
	 $(a)$ Ground truth model from USGS-STRM1 data relative to the area with geographical coordinates;
	 $[\mathrm{N}\,40^{\circ}48'50'',\,\mathrm{N}\,40^{\circ}52'50'']\times[\mathrm{E}\,14^{\circ}45'50'',\,\mathrm{E}\,14^{\circ}50'00'']$. 
	 $(b)$ Sample set $K_1$ formed by only level lines at regular height interval of $66\,\mathrm{m}$. 
		The set $K_1$ contains $19\%$ of the ground truth points.
	 $(c)$ Sample set $K_2$ formed by taking randomly $30\%$ of the
		points belonging to the level lines of the set $K_1$ and scattered points 
		corresponding to $5\%$ density. The sample set $K_2$ contains $9\%$ of the 
		ground truth points.
	}
\end{figure}

\begin{figure}[htbp]
  \centerline{$\begin{array}{cc}
	\includegraphics[width=0.50\textwidth]{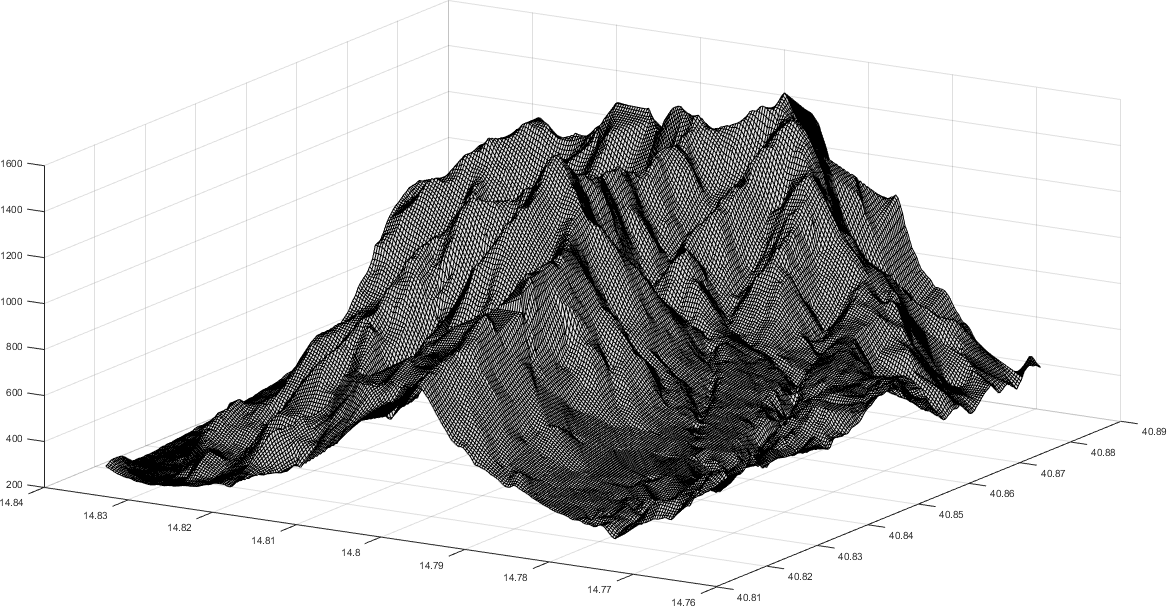}
	&\includegraphics[width=0.50\textwidth]{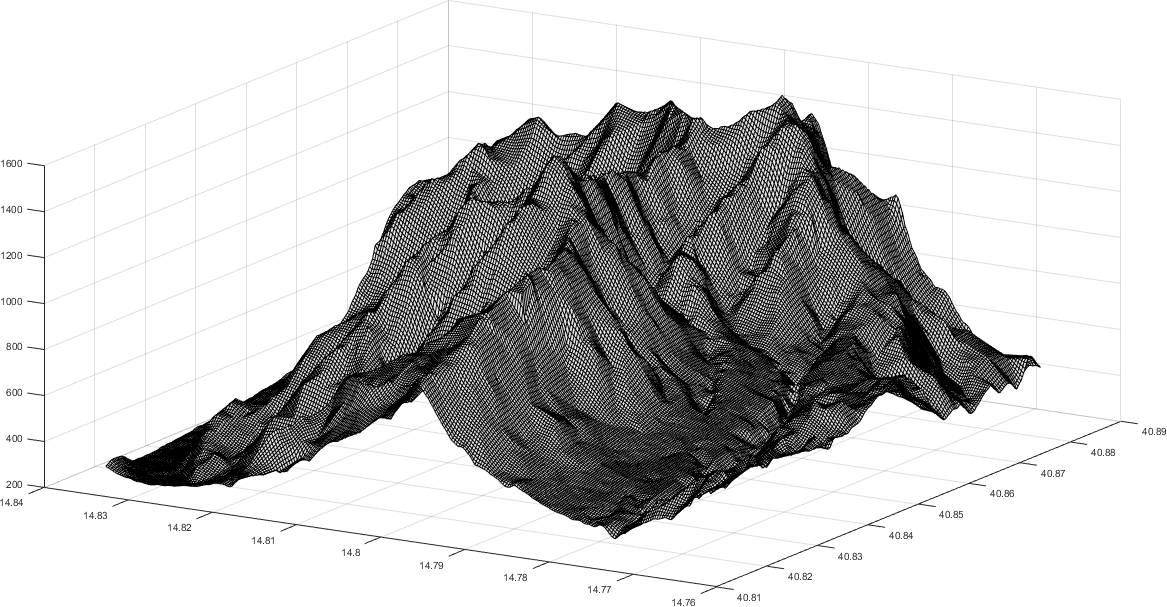}\\
	(a) & (b)\\
	\includegraphics[width=0.225\textwidth]{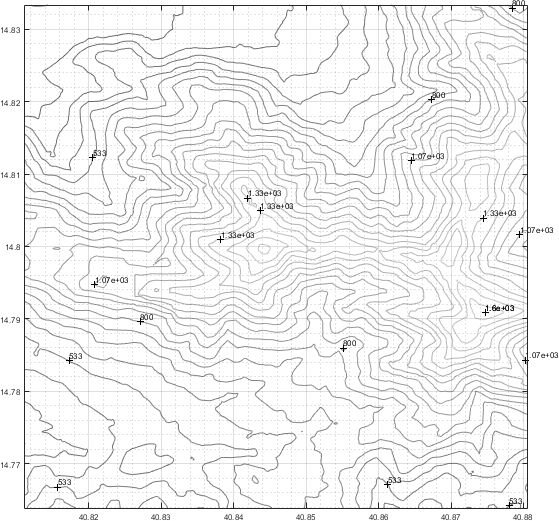}
	&\includegraphics[width=0.225\textwidth]{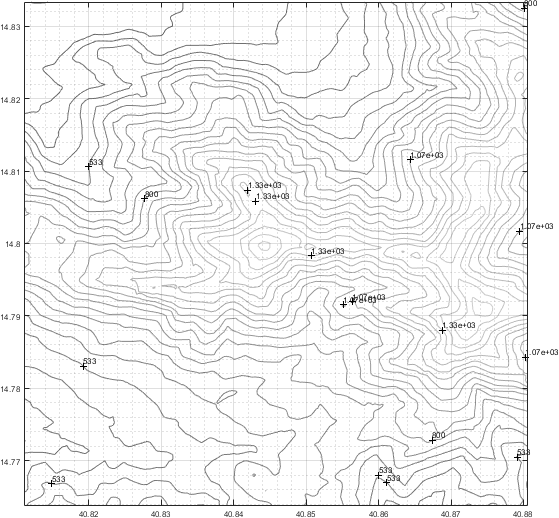}\\
	(c) & (d)
    \end{array}$
}
   \caption{\label{Fig.NumEx7.2.4b} Example \ref{Sec7.2.4}.
   Reconstruction of real-world digital elevation maps.
	 $(a)$ Graph of $A_{\lambda}^M(f_K)$ for sample set $K_1$. 
		Parameters: $\lambda=1\cdot 10^3$, $M=1\cdot 10^6$.   
		Relative $L^2$-Errors: $\epsilon=0.01560$, $\epsilon_K=0$.
	 $(b)$ Graph of $A_{\lambda}^M(f_K)$ for sample set $K_2$. 
		Parameters: $\lambda=1\cdot 10^3$, $M=1\cdot 10^6$.   
		Relative $L^2$-Errors: $\epsilon=0.0117$, $\epsilon_K=0$.
	 $(c)$ Isolines of $A_{\lambda}^M(f_K)$ from sample set $K_1$
		at regular heights of $66\,\mathrm{m}$.
	 $(d)$ Isolines of $A_{\lambda}^M(f_K)$ from sample set $K_2$
		 at regular heights of $66\,\mathrm{m}$.
	}
\end{figure}

\begin{figure}[htbp]
  \centerline{$\begin{array}{cc}
	\includegraphics[width=0.50\textwidth]{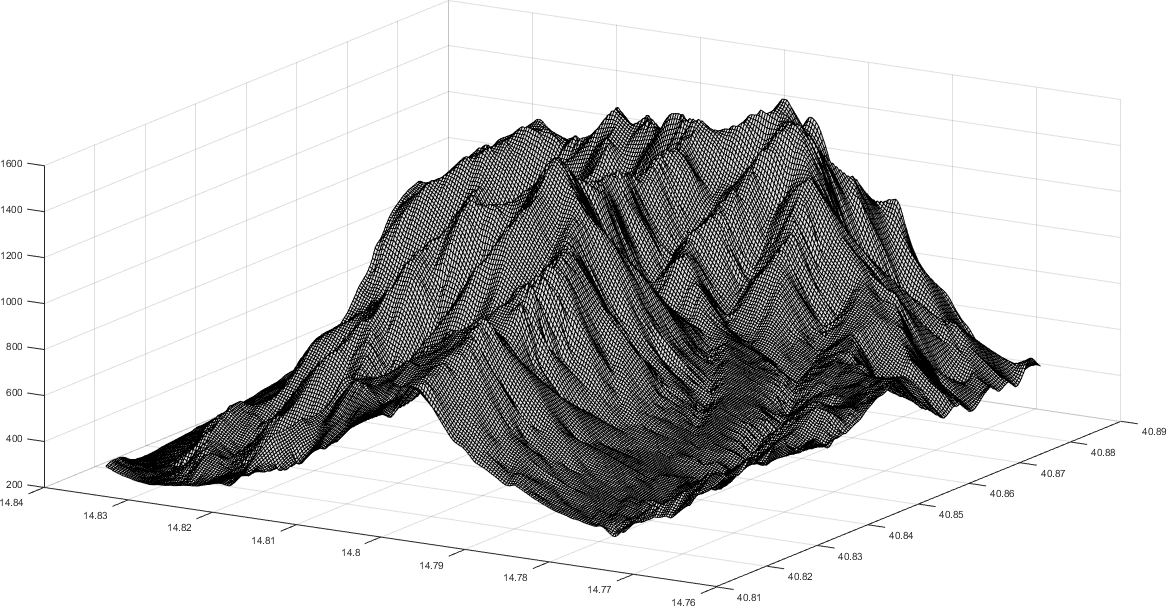}
	&\includegraphics[width=0.50\textwidth]{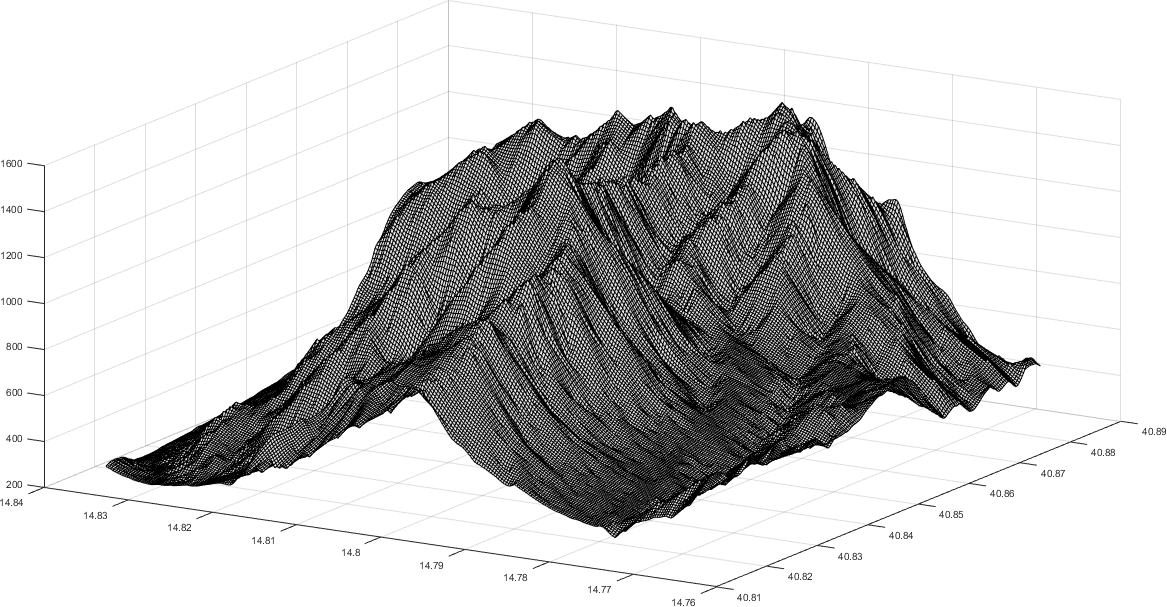}\\
	(a) & (b)\\
	\includegraphics[width=0.225\textwidth]{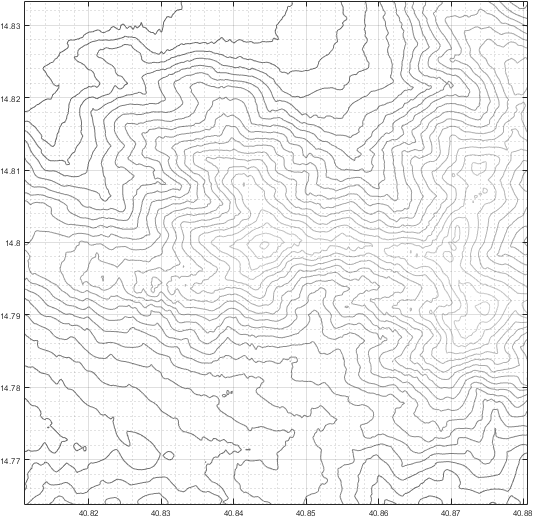}
	&\includegraphics[width=0.225\textwidth]{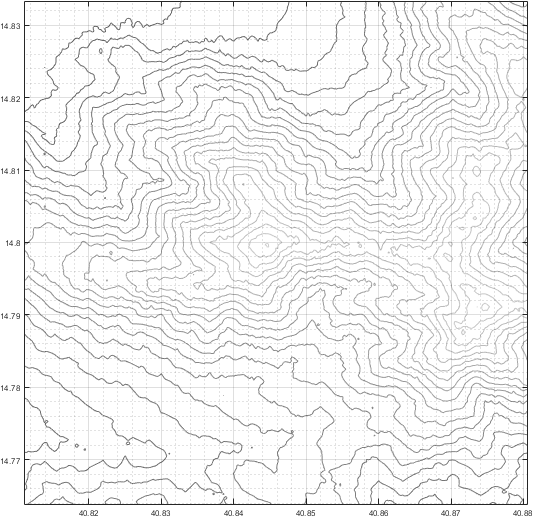}\\
	(c) & (d)
    \end{array}$
}
   \caption{\label{Fig.NumEx7.2.4c} Example \ref{Sec7.2.4}.
   Reconstruction of real-world digital elevation maps.
	 $(a)$ Graph of the AMLE Interpolant from set $K_1$.   
		Relative $L^2$-Error: $\epsilon=0.0214$.
	 $(b)$ Graph of the AMLE Interpolant from set $K_2$.   
		Relative $L^2$-Error: $\epsilon=0.0226$.
	 $(c)$ Isolines of the AMLE Interpolant from sample set $K_1$
		at regular heights of $66\,\mathrm{m}$.
	 $(d)$ Isolines of the AMLE Interpolant from sample set $K_2$
		 at regular heights of $66\,\mathrm{m}$.
	}
\end{figure}

%%%%%%%%%%%%%%%%%%%%%%%%%%%%%%%%%%%%%%%%%%%%%%%%%%%%%%

\begin{table}[tbhp]
\centerline{
\begin{tabular}{c|c|c|}
\cline{2-3}
    & \multicolumn{2}{c|}{ $\epsilon$} \\ \hline
\multicolumn{1}{ |c| }{Sample set}	& $A_{\lambda}^M(f_K)$  &  AMLE     \\ \hline
\multicolumn{1}{ |c| }{$K_1$}		& $0.0156$ & $0.02137$	\\ \hline
\multicolumn{1}{ |c| }{$K_2$}		& $0.0117$ & $0.02261$	\\ \hline
\end {tabular}
}
\caption{\label{Tab:RealDat} Relative $L^2$-error for the DEM Reconstruction from the two sample sets
using the $A_{\lambda}^M(f_K)$ and the AMLE interpolant.}
\end{table}

We consider here the problem of producing a Digital Elevation Map from a sample of the 
the NASA SRTM global digital 
elevation model of Earth land. The data provided by the 
National Elevation Dataset \cite{GEMHC09} contain geographical coordinates (latitude, longitude and elevation)
of points sampled at one arc-second intervals in latitude and longitude. For our experiments,
we choose the region defined by the coordinates 
$[\mathrm{N}\,40^{\circ}48'50'',\,\mathrm{N}\,40^{\circ}52'50'']\times[\mathrm{E}\,14^{\circ}45'50'',\,\mathrm{E}\,14^{\circ}50'00'']$
extracted from the SRTM1 cell $N40E014.hgt$ \cite{SRTM1}. Such region consists of an area 
with extension $7.413\,\mathrm{Km}\times 5.844\,\mathrm{km}$ and height varying between $266\,\mathrm{m}$ and $1600\,\mathrm{m}$, 
with variegated topography features. In the digitization by the US Geological Survey, 
each pixel represents a $30\,\mathrm{m} \times 30\,\mathrm{m}$ patch.
Figure \ref{Fig.NumEx7.2.4a}$(a)$ displays the elevation model from the SRTM1 data
which we refer in the following to as the ground truth model. 
We will take a sample $f_K$ of such data, make the reconstruction using the $A_{\lambda}^M(f_K)$
and the AMLE interpolant, and compare them with the ground truth model.
In the numerical experiments, we consider two sample data, characterized by different 
data density and typo of information. The first, which we refer to as sample set $K_1$, 
consists only of level lines at regular height interval of $66\,\mathrm{m}$ and contains 
the $19\%$ of the ground truth real digital data.
The second sample set, denoted by $K_2$, has been formed by taking randomly the $30\%$ of the
points belonging to the level lines of the set $K_1$ and scattered points corresponding to
$5\%$ density so that the sample set $K_2$ amounts to about $9\%$ of the ground truth points.
The two sample sets $K_1$ and $K_2$ are shown in Figure \ref{Fig.NumEx7.2.4a}$(b)$ and 
Figure \ref{Fig.NumEx7.2.4a}$(c)$, respectively. The graph of the $A_{\lambda}^M(f_K)$ 
interpolant and of the $AMLE$ interpolant for the two sample sets along with the respective 
isolines at equally spaced heighs equal to $66\,\mathrm{m}$, are displayed in 
Figure \ref{Fig.NumEx7.2.4b} and Figure \ref{Fig.NumEx7.2.4c}, respectively, 
whereas Table \ref{Tab:RealDat} contains the values of the 
relative $L^2$-error between such interpolants and the ground truth model.
Though both reconstructions are comparable visually to the ground truth model, 
a closer inspection of the pictures show that the reconstruction from the synthetic data, 
the AMLE interpolant does not reconstruct correctly the mountains peaks, which appear to be smoothed,
and introduce artificial ridges along the slopes of the mountains. In contrast, the 
$A_{\lambda}^M(f_K)$ interpolant appears to better for capturing features of the ground truth model.
Finally, we also note that though the sample set $K_1$ contains a number 
of ground truth points higher than the sample set $K_2$, the reconstruction from $K_2$ appears 
to be better than the one obtained from $K_1$. This behaviour was found for both interpolations,
though it is more notable in the case of the $A_{\lambda}^M(f_K)$ interpolant. By taking scattered data,
we are able to get a better characterization of irregular surfaces, compared
to the one obtained from a structured representation such as provided by the level lines.

%%%%%%%%%%%%%%%%%%%%%%%%%%%%%%%%%%%%%%%%%%%%%%%%%%%%%%

\subsubsection{Salt \& Pepper Noise Removal}\label{Sec7.2.5}

%%%%%%%%%%%%%%%%%%%%%%%%%%%%%%%%%%%%%%%%%%%%%%%%%%%%%%

\begin{figure}[htbp]
	\centerline{$\begin{array}{cc}
		\includegraphics[height=5cm]{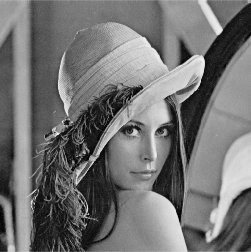}&
		\includegraphics[height=5cm]{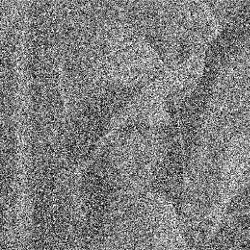}\\
		(a)&(b)\\
		\includegraphics[height=5cm]{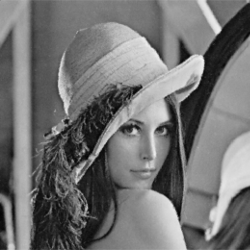}&
		\includegraphics[height=5cm]{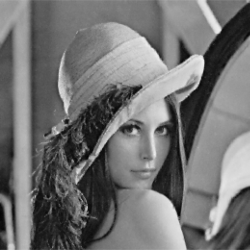}\\
		 (c)&(d)
	\end{array}$}
	\caption{\label{Fig.NumEx7.2.5} Example \ref{Sec7.2.5}.
		$(a)$ Original image with size $512\times 512$;
		$(b)$ Original image covered by a salt \& pepper noise density of $70\%$. $\mathrm{PSNR}=6.99\,\mathrm{dB}$; 
		$(c)$ Restored image $A_{\lambda}^M(f_K)$ with $K$ the set of the pixels not corrupted by the salt 
			\& pepper noise when the corrupted image is enlarged symmetrically by two pixels on each side,
			$\lambda=15$ and $M=1E13$. 
			$\mathrm{PSNR}=31.865\,\mathrm{dB}$. If the boundary pixels were noise-free,
			the corresponding restored image would have $\mathrm{PSNR}=31.910\,\mathrm{dB}$. 
		$(d)$ Restored image by the two-stage TV-based method described in \cite{CCM07,CHN05}
		with $K$  the set of the pixels not corrupted by the salt 
			\& pepper noise when the corrupted image is enlarged symmetrically by two pixels on each side. 
			$\mathrm{PSNR}=31.134\,\mathrm{dB}$.
			If the boundary pixels were noise-free, 
			the corresponding restored image would have $\mathrm{PSNR}=31.175\,\mathrm{dB}$.
	}
\end{figure}

As an application of  scattered data approximation to image processing, 
we consider here the restoration of 
an image corrupted by salt \& pepper noise. This is an impulse type noise that is caused, for instance, by 
malfunctioning pixels in camera sensors or faulty memory locations
in hardware, so that information is lost at the faulty pixels and the 
corrupted pixels are set alternatively
to the minumum or to the maximum value of the range of the image values. 
When the noise density  is low, about less than $40\%$, the median filter is quite effective for restoring the 
image. However, this filter loses  its denoising power for higher noise density given that 
details and features of the original image are smeared out. In those cases, other techniques must be 
applied; one possibility is the two-stage TV-based method proposed in \cite{CHN05}.
In the following numerical experiments,  we consider the image displayed in Figure \ref{Fig.NumEx7.2.5}$(a)$ 
with size $512\times 512$ pixels, damaged by $70\%$ salt \& pepper noise. The resulting corrupted image is displayed in 
Figure \ref{Fig.NumEx7.2.5}$(b)$ where only  $78643$ pixels out of the total $262144$ pixels 
carry true information.
The true image values represent our sample function $f_K$ whereas the set of
the true pixels forms our sample set $K$. 
To assess the restoration performance we use 
the peak signal-to-noise ratio ($\mathrm{PSNR}$) which is expressed in the units of $\mathrm{dB}$ and, for an $8-$bit image,
is defined by
\begin{equation}
	\mathrm{PSNR}=10\log_{10}\displaystyle \frac{255^2}{\frac{1}{mn}\sum_{i,j}|f_{i,j}-r_{i,j}|^2}
\end{equation}
where $f_{i,j}$ and $r_{i,j}$ denote the pixels values of the original and restored image, respectively,
and $m,\,n$ denote the size of the image $f$. 
In our numerical experiments, we have considered the following cases. 
The first one assumes the set $K$ to be given by the noise-free interior pixels of the corrupted image together with the boundary pixels of the original image.
In the second case,  $K$ is just the set of the noise-free pixels of the corrupted image, without any special 
consideration on the image boundary pixels.  
In analysing this second case, to reduce the boundary effects produced by the application of Algorithm \ref{Algo:CnvxEnv}, 
we have applied our method to an enlarged image and then 
restricted the resulting restored image to the original domain. The enlarged image has been obtained by padding a 
fixed number of pixels before the first image element and after the last image element along each dimension, making 
mirror reflections with respect to the boundary. The values used for padding are all from the corrupted image.
In our examples, we have considered two versions of enlarged images, 
obtained by padding the corrupted image with 2 pixels and 10 pixels, respectively. 
Table \ref{Table:PSNR} compares the values of the $\mathrm{PSNR}$ of the restored images by our method and the 
TV-based method applied to the corrupted image with noise-free boundary and to the two versions of the enlarged images 
with the boundary values of the enlarged images given by the padded noisy image data. 
We observe that there are no important variations in the denoising result between the different methods of treating the image boundary.
This is also reflected by the close value of the $\mathrm{PSNR}$ of the resulting restored images.
For $70\%$ salt \& pepper noise, Figure \ref{Fig.NumEx7.2.5}$(c)$ displays the restored image $A_{\lambda}^M(f_K)$ with $K$ equal to the true 
set that has been enlarged by two pixels, whereas Figure \ref{Fig.NumEx7.2.5}$(d)$ shows the restored image 
by the TV-based method \cite{CCM07,CHN05} using the same set $K$. 
Although the visual quality of the images restored from $70\%$ noise corruption is comparable between our method 
and the TV-based method, the $\mathrm{PSNR}$ using our method is higher than that for the TV-based method in all 
of the experiments reported in Table \ref{Table:PSNR}.
An additional advantage of our method is its speed. 
Our method does not require initialisation which is in contrast with the two-stage TV-based method, for which the initialisation, for instance, 
is given by the restored image using an adaptive median filter.

\begin{figure}[htbp]
	\centerline{$\begin{array}{cc}
		\includegraphics[height=5cm]{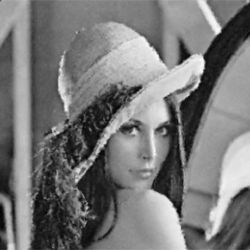}&
		\includegraphics[height=5cm]{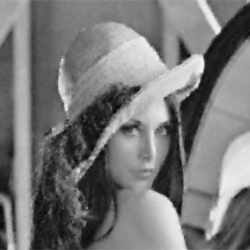}\\
		(a)&(b)\\
		\includegraphics[height=5cm]{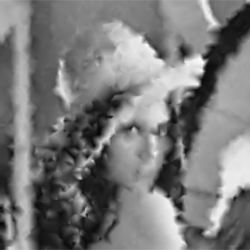}&
		\includegraphics[height=5cm]{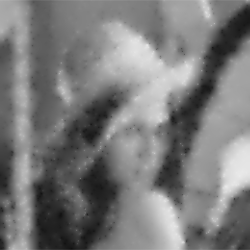}\\
		 (c)&(d)
	\end{array}$}
	\caption{\label{Fig.NumEx7.2.5a} Example \ref{Sec7.2.5}.
		  Restoration of $90\%$ corrupted image ($\mathrm{PSNR}=5.901\,\mathrm{dB}$) by:	
			$(a)$ Restored image $A_{\lambda}^M(f_K)$,  with $K$ the set of the pixels not corrupted by the salt 
				\& pepper noise when the corrupted image is enlarged symmetrically by two 
				pixels on each side, $\lambda=15$ and $M=1E13$. 
				$\mathrm{PSNR}=27.506\,\mathrm{dB}$.  
			$(b)$ Restored Image by the two-stage TV-based method described in \cite{CCM07,CHN05}
				with the same set $K$ as in $(a)$. 
				$\mathrm{PSNR}=26.564\,\mathrm{dB}$.
		  Restoration of $99\%$ corrupted image ($\mathrm{PSNR}=5.492 \,\mathrm{dB}$) by:	
			$(c)$ Restored image $A_{\lambda}^M(f_K)$,  with $K$ the set of the pixels not corrupted by the salt 
				\& pepper noise when the corrupted image is enlarged symmetrically by ten pixels 
				on each side, $\lambda=15$ and $M=1E13$. 
				$\mathrm{PSNR}=21.972\,\mathrm{dB}$.  
			$(d)$ Restored Image by the two-stage TV-based method described in \cite{CCM07,CHN05}
				with the same set $K$ as in $(c)$. 
				$\mathrm{PSNR}=20.492\,\mathrm{dB}$.
	}
\end{figure}

Finally, to demonstrate the performance of our method in some extreme cases of very sparse data, we consider 
cases of noise density equal to $90\%$ and $99\%$.
Figure \ref{Fig.NumEx7.2.5a} displays the restored image by the compensated convexity based method and by the TV-based method for cases where 
$K$ are padded by two pixels and ten pixels for $90\%$ and $99\%$ noise level, respectively. As far as the visual quality of the restored 
images is concerned, and to the extent that such judgement can make sense given the high level of noise density,
the inspection of Figure \ref{Fig.NumEx7.2.5a} seems to indicate that 
$A_{\lambda}^M(f_K)$ gives a better approximation of details than the TV-based restored image. This is also reflected by the values
of the $\mathrm{PSNR}$ index in Table \ref{Table:PSNR}.

%%%%%%%%%%%%%%%%%%%%%%%%%%%%%%%%%%%%%%%%%%%%%%%%%%%%%%%%%%%%%%%%%%%%%%%%%%%%%%%%%%%%%%%%%%%%%%%%%%%%%%%%%%%%%%%

\begin{table}[tbhp]
\centerline{
\begin{tabular}{c|c|c|c|c|c|c|}
	\cline{2-7}			& \multicolumn{6}{c|}{$\mathrm{PSNR}$} \\ \cline{2-7}
	\cline{2-7}			& \multicolumn{2}{c|}{$K$ with noise-free boundary}	& \multicolumn{2}{c|}{$K$ padded  by two pixels}		& \multicolumn{2}{c|}{$K$ padded  by ten pixels}  \\ \hline
\multicolumn{1}{|c|}{Noise Density}	& $A_{\lambda}^M(f_K)$	& TV				& $A_{\lambda}^M(f_K)$	& TV				& $A_{\lambda}^M(f_K)$	& TV			  \\ \hline
\multicolumn{1}{|c|}{$70\%$ ($6.990\,\mathrm{dB}$)}	& $31.910\,\mathrm{dB}$	& $31.175\,\mathrm{dB}$	& $31.865\,\mathrm{dB}$	& $31.134\,\mathrm{dB}$	& $31.869\,\mathrm{dB}$	& $31.136\,\mathrm{dB}$	\\ \hline
\multicolumn{1}{|c|}{$90\%$ ($5.901\,\mathrm{dB}$)}	& $27.574\,\mathrm{dB}$	& $26.625\,\mathrm{dB}$	& $27.506\,\mathrm{dB}$	& $26.564\,\mathrm{dB}$	& $27.513\,\mathrm{dB}$	& $26.566\,\mathrm{dB}$	\\ \hline
\multicolumn{1}{|c|}{$99\%$ ($5.492\,\mathrm{dB}$)}	& $22.076\,\mathrm{dB}$	& $20.595\,\mathrm{dB}$	& $21.761\,\mathrm{dB}$	& $20.469\,\mathrm{dB}$	& $21.972\,\mathrm{dB}$	& $20.492\,\mathrm{dB}$	\\ \hline
\end {tabular}
}
\caption{\label{Table:PSNR} 
	Comparison of $\mathrm{PSNR}$ of the restored images by the compensated convexity based method ($A_{\lambda}^M(f_K)$) 
	and by the two-stage TV-based method (TV), for different sets $K$.
}
\end{table}

%%%%%%%%%%%%%%%%%%%%%%%%%%%%%%%%%%%%%%%%%%%%%%%%%%%%%%%%%%%%%%%%%%%%%%%%%%%%%%%%%%%%%%%%%%%%%%%%%%%%%%%%%%%%%%%

%%%%%%%%%%%%%%%%%%%%%%%%%%%%%%%%%%%%%%%%%%%%%%%%%%%%%%%%%%%%%%%%%%%%%%%%%%%%%%%%%%%%%%%%%%%%%%%%%%%%%%%%%%%%%%%

\subsection{Image inpainting}\label{Sec7.3}

%%%%%%%%%%%%%%%%%%%%%%%%%%%%%%%%%%%%%%%%%%%%%%%%%%%%%%%%%%%%%%%%%%%%%%%%%%%%%%%%%%%%%%%%%%%%%%%%%%%%%%%%%%%%%%%

\begin{figure}[htbp]
  \centerline{$\begin{array}{cc}
	\includegraphics[height=5cm]{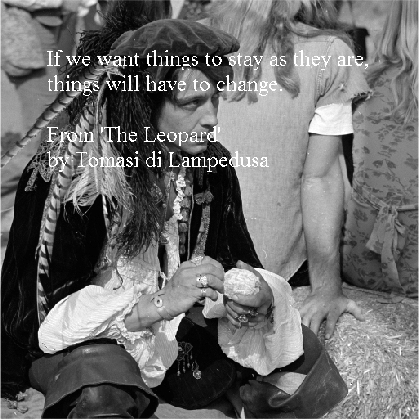}&
	\includegraphics[height=5cm]{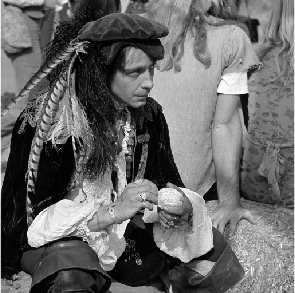}\\
	(a)&(b)\\
	\includegraphics[height=5cm]{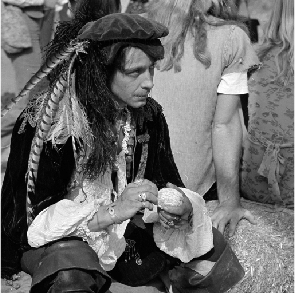}&
	\includegraphics[height=5cm]{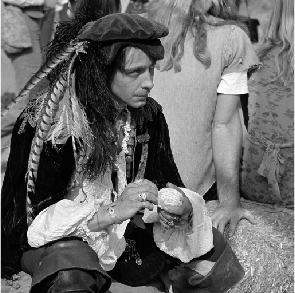}\\
	(c)&(d)
    \end{array}$
}
   \caption{\label{Fig.NumEx7.3a} Example \ref{Sec7.3}.
	Inpainting of the text overprinted on an image: 
	$(a)$ Original image with overprinted text.
	$(b)$ Restored image $A_{\lambda}^M(f_K)$ with $K$ the set to be inpainted, $\lambda=250$ and $M=1\cdot 10^4$. 
		 Computed value for $\mathrm{PSNR}=42.2066 \,\mathrm{dB}$; Relative $L^2$-error $\epsilon=0.016139$.
	$(c)$ Restored image by the AMLE method described in \cite{Sch15,PS16}. 
		 Computed value for $\mathrm{PSNR}=39.4405 \,\mathrm{dB}$.
		 Relative $L^2$-error $\epsilon=0.022192$.
	$(d)$ Restored image by the Split Bregman inpainting method described in \cite{Get12}. 
		 Computed value for $\mathrm{PSNR}=41.0498\,\mathrm{dB}$.
		 Relative $L^2$-error $\epsilon=0.018438$.
	}
\end{figure}

\begin{figure}[htbp]
  \centerline{$\begin{array}{ccc}
		\includegraphics[height=2.5cm]{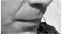}&
		\includegraphics[height=2.5cm]{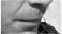}\\
		(a)&(b)\\
		\includegraphics[height=2.5cm]{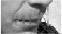}&
		\includegraphics[height=2.5cm]{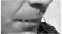}&\\
		(c)&(d)
    \end{array}$
}
   \caption{\label{Fig.NumEx7.3b} Example \ref{Sec7.3}.
		Comparison of a detail of the original image with the 
		corresponding detail of the restored images according to the 
		compensated convexity method and the TV-based method:  
		$(a)$ Lips detail of the original image without overprinted text.
		$(b)$ Lips detail of the restored image $A_{\lambda}^M(f_K)$.
		$(c)$ Lips detail of the AMLE-based restored image.
		$(d)$ Lips detail of the TV-based restored image.		
	}
\end{figure}

As an example of image inpainting, we consider the problem of removing text overprinted on the image displayed 
in Figure \ref{Fig.NumEx7.3a}$(a)$. If we denote by $P$ the set of pixels containing the overprinted text, and 
by $\Omega$ the domain of the whole image, then $K=\Omega\setminus P$ is the set of the true pixels and 
the inpainting problem is in fact the problem of reconstructing the image over $P$ from knowing $f_K$, if we denote by 
$f$ the original image values. To assess the performance of our reconstruction compared to state-of-art inpainting methods, 
we compare our method 
with the total variation based image inpainting method solved by the split Bregman method described in \cite{Get12}
and with the AMLE inpainting reported in \cite{Sch15}.
The restored image $A_{\lambda}^M(f_K)$ obtained by our compensated convexity method is displayed in 
Figure \ref{Fig.NumEx7.3a}$(b)$, the restored image by the AMLE method is shown 
in Figure \ref{Fig.NumEx7.3a}$(d)$ whereas \ref{Fig.NumEx7.3a}$(c)$ presents 
the restored image by the the split Bregman inpainting method. 
All the restored images look visually quite good. However, if we use 
 the $\mathrm{PSNR}$ as a measure of the quality of the restoration, we find that $A_{\lambda}^M(f_K)$
has a value of $\mathrm{PSNR}$ equal to $ 42.2066\,\mathrm{dB}$,  the split Bregman inpainting
restored image gives a value for $\mathrm{PSNR}=41.0498\,\mathrm{dB}$,
whereas the AMLE restored image has $\mathrm{PSNR}$ equal to $ 39.4405\,\mathrm{dB}$.

Finally, to assess how well $A_{\lambda}^M(f_K)$ is able to preserve image details 
and not to introduce unintended effects such as image blurring and staircase effects, Figure \ref{Fig.NumEx7.3b}
displays details of the original image and of the restored images by the three methods. 
Once again, the good performance of $A_{\lambda}^M(f_K)$ can be appreciated visually.

%%%%%%%%%%%%%%%%%%%%%%%%%%%%%%%%%%%%%%%%%%%%%%%%%%%%%%%%%%%%%%%%%%%%%%%%%%%%%%%%%%%%%%%%%%%%%%%%%%%%
%%%%%%%%%%%%%%%%%%%%%%%%%%%%%%%%%%%%%%%%%%%%%%%%%%%%%%%%%%%%%%%%%%%%%%%%%%%%%%%%%%%%%%%%%%%%%%%%%%%%
%%%%%%%%%%%%%%%%%%%%%%%%%%%%%%%%%%%%%%%%%%%%%%%%%%%%%%%%%%%%%%%%%%%%%%%%%%%%%%%%%%%%%%%%%%%%%%%%%%%%

\setcounter{equation}{0}
\section{Proofs of the Main Results}\label{Sec.Proofs}

%%%%%%%%%%%%%%%%%%%%%%%%%%%%%%%%%%%%%%%%%%%%%%%%%%%%%%%%%%%%%%%%%%%%%%%
%%%%%%%%%%%%%%%%%%%%%%%%%%%%%%%%%%%%%%%%%%%%%%%%%%%%%%%%%%%%%%%%%%%%%%%%%%
\begin{proof}(\cref{Prp.Rest}) 
We write $(x,y)\in \mathbb{R}^{n+m}$
with $x\in \mathbb{R}^n$ and $y\in \mathbb{R}^m$. We only prove the result for the
upper transform as the proof of the lower transform is similar. By the definition of the 
upper transform, we have
\[
	\co[\lambda|\cdot|^2-f](x)=\lambda|x|^2-C^u_\lambda(f(x)),\quad x\in \mathbb{R}^n\,.
\]
We show that $\co[\lambda|\cdot|^2-f](x)$ is also the convex envelope of the function 
$\lambda(|x|^2+|y|^2)-g^{-M}(x,y)$ restricted to $z=0$. By definition,
\[
	\lambda|x|^2-C^u_\lambda(f(x))=\co[\lambda|\cdot|^2-f](x)\leq \lambda|x|^2-f(x)\leq \lambda(|x|^2+|y|^2)-g^{-M}(x,y)
\]
as $f(x)\geq g^{-M}(x,y)$ for all $x\in\mathbb{R}^n$ and $y\in \mathbb{R}^m$. 
Thus for $y=0$,
\[
	\co[\lambda|\cdot|^2-f](x)\leq \co[\lambda(|x|^2+|y|^2)-g^{-M}(x,y)]|_{y=0}\,.
\]
On the other hand,
\[
	\co[\lambda(|x|^2+|y|^2)-g^{-M}(x,y)]|_{y=0}\leq
	\lambda|x|^2-g^{-M}(x,0)= \lambda|x|^2-f(x)\,.
\]
Since the restriction of a convex function to a linear subspace remains convex, we also see that
\[
	\co[\lambda(|x|^2+|y|^2)-g^{-M}(x,y)]|_{y=0}\leq \co[\lambda|\cdot|^2-f](x)\,.
\]
Thus
\[
	\co[\lambda(|x|^2+|y|^2)-g^{-M}(x,z)]|_{y=0}= \co[\lambda|\cdot|^2-f](x)\,,
\]
hence the conclusion follows. 
\end{proof}

%%%%%%%%%%%%%%%%%%%%%%%%%%%%%%%%%%%%%%%%%%%%%%%%%%%%%%%%%%%%%%%%%%%%%%%%%%%%%%%%%%%%%%%

\begin{proof}(\cref{Thm.LSaprx}) 
Note first that it follows from the fact that $a_0 < a_1 < \cdots < a_m$, $m \in \mathbb{N}$, that  $V_{a_i} \subset V_{a_j}$ 
for all $0\leq i<j\leq m$.
Also, by the translation invariant property of compensated convex
transforms,  we may assume without loss of generality that $x_0=0$, so that
$$C^l_\lambda(f^M_K)(0)=\co[f^M_K+\lambda|\cdot|^2](0), \;\;\;\; 
C^u_\lambda(f^{-M}_{K})(0)=\co[\lambda|\cdot|^2-f^{-M}_{K}](0).$$

%%%%%%%%%%%%%%%%%%%%%%%%%%%%%%%%%%%%%%%%%%%%%%%%%%%%%%%%%%%%%%%%%%%%%%%%%%
%%%%%%%%%%%%%%%%%%%%%%%%%%%%%%%%%%%%%%%%%%%%%%%%%%%%%%%%%%%%%%%%%%%%%%%%%%

\noindent \textit{(i):}   Suppose that  $x_0=0 \in \Gamma_{a_k}$ and consider the
constant function $\ell(x)=a_k$. Clearly $a_k=f^M_{K}(0)+\lambda|0|^2$. Next we
 show that $a_k\leq f^M_{K}(x)+\lambda|x|^2$ for $x\in \Gamma_{a_j}$ for $j\neq k$.
Thus we need to prove that $a_k\leq a_j+\lambda|x|^2$.  Since $0\in \Gamma_{a_k}$ and
$x\in\Gamma_{a_j}$, we have $|x|^2\geq\delta_0^2$. Under our assumption on $\lambda$,
we see that $a_k\leq a_j+\lambda|x|^2$ holds. Since $a_k<M$, we have $a_k\leq f^M_{K}(x)+\lambda|x|^2$ 
for all $x\in \mathbb{R}^n$, hence
$C^l_\lambda(f^M_K)(0)=a_k$. Similarly we can show that $C^u_\lambda(f^{-M}_K)(0)=a_k$, so that 
$ A^M_\lambda(f_{K})(0)=a_{k}$.\\

%%%%%%%%%%%%%%%%%%%%%%%%%%%%%%%%%%%%%%%%%%%%%%%%%%%%%%%%%%%%%%%%%%%%%%%%%%
%%%%%%%%%%%%%%%%%%%%%%%%%%%%%%%%%%%%%%%%%%%%%%%%%%%%%%%%%%%%%%%%%%%%%%%%%%

\noindent \textit{(ii):} Since (i) clearly ensures that \eqref{ineq} holds whenever $f(x_0) = a_i$ 
for some $0\leq i\leq m$, it remains to consider $x_0=0$ such that
$a_i < f(x_0) < a_{i+1}$ for some $0\leq i\leq m-1.$ Now  define
\begin{equation}\label{Eq.Prf.01} 
	f^M_{K^-_i}(x)=\left\{\begin{array}{ll}
			\displaystyle f^M_{K}(x),	& \displaystyle x\notin \Gamma_{a_{i+1}}\,,\\[1.5ex]
			\displaystyle  a_i,		& \displaystyle x\in \Gamma_{a_{i+1}}\,;
		\end{array}\right.\qquad
	f^M_{K^+_i}(x)=\left\{\begin{array}{ll}
			\displaystyle f^M_{K}(x), & \displaystyle x\notin \Gamma_{a_{i}}\,,\\[1.5ex]
			\displaystyle a_{i+1},	  & \displaystyle x\in \Gamma_{a_{i}}\,.
		\end{array}\right.
\end{equation}
Clearly $f^M_{K^-_i}(x)\leq f^M_{K}(x)\leq f^M_{K^+_i}(x)$ and $f^{-M}_{K^-_i}(x)\leq f^{-M}_{K}(x)\leq f^{-M}_{K^+_i}(x)$
for $x\in\mathbb{R}^n$,  so that
\begin{equation}\label{Eq.Prf.02} 
	C^l_\lambda(f^M_{K^-_i})(x)\leq C^l_\lambda(f^M_{K})(x)\leq C^l_\lambda(f^M_{K^+_i})(x), \quad
	C^u_\lambda(f^{-M}_{K^-_i})(x)\leq C^u_\lambda(f^{-M}_{K})(x)\leq C^u_\lambda(f^{-M}_{K^+_i})(x),
	\quad x\in \mathbb{R}^n
\end{equation}
and hence by definition,
\begin{equation}\label{Eq.Prf.03} 
	A^M_\lambda(f_{K^-_i})(x)\leq A^M_\lambda(f_{K})(x)\leq A^M_\lambda(f_{K^+_i})(x),\quad x\in 
	\mathbb{R}^n\,.
\end{equation}
Next we will prove that
\begin{equation}\label{Eq.Prf.04}  
	A^M_\lambda(f_{K^-_i})(0)=a_i,\quad
	A^M_\lambda(f_{K^+_i})(0)=a_{i+1}\,.
\end{equation}

\smallskip

\noindent We first show that $\co[f^M_{K^-_i}+\lambda|\cdot|^2](0)\geq a_i$.
Clearly $a_i\leq a_i+\lambda|x|^2=f^M_{K^-_i}(x)+\lambda|x|^2$ for $x\in \Gamma_{a_i}\cup \Gamma_{a_{i+1}}$. 
For $x\in \Gamma_{a_j}$ with $j\neq i,\, i+1$,  $a_i\leq a_j+ \lambda|x|^2$   if $a_i-a_j\leq \lambda|x|^2.$ 
This inequality holds if
$a_m-a_0\leq \lambda\delta_0^2$, that is, for $\lambda\geq  (a_m-a_0)/\delta_0^2$ which is  
what we have assumed. The inequality $|x|\geq \delta_0$ for
$x\in \Gamma_{a_j}$ can be proved by applying the intermediate value theorem to $f$. If $j<i$, 
as $f(0)>a_i$ and $f(x)=a_j<a_i$, we have, by the intermediate value
theorem, that there is some $\xi\in (0,\,1)$ such that $f(\xi x)=a_i$, that is, $\xi x\in  \Gamma_{a_i}$. 
Thus $|x|>(1-\xi)|x|=|x-\xi x|\geq \delta_0$ as $x\in\Gamma_{a_j}$ and
$\xi x\in  \Gamma_{a_i}$. If $j>i+1$, we have $f(0)<a_{i+1}$ and $f(x)=a_j>a_{i+1}.$
Again we can use the same method to show that $|x|\geq \delta_0$.

\medskip
\noindent By definition of the convex envelope, we see that there is an affine function $\ell$ such that 
$\ell(x)\leq f^M_{K^-_i}(x)+\lambda|x|^2$ for $x\in\mathbb{R}^n$ and
$\ell(0)=\co[f^M_{K^-_i}+\lambda|\cdot|^2](0)$. From the proof above, we see that $\ell(0)\geq a_i$.
Furthermore, if we let
$K_l=\{x\in \mathbb{R}^n,\; \ell(x)=f^M_{K^-_i}(x)+\lambda|x|^2\}$, then $0\in \co[K_l]$ and  
$\ell(x)=\co[f^M_{K^-_i}+\lambda|\cdot|^2](x)$ for $x\in \co[K_l]$.

\medskip
\noindent By \cite[Proposition 3.3]{ZCO16}, we see that $K_l\subset K$. 
Now we show that $K_l\subset \Gamma_{a_i}\cup \Gamma_{a_{i+1}}$. 
If this is not the case, then $K_l\cap \Gamma_{a_k}\neq \varnothing$ for some $k \not\in \{i, i+1\}$. We consider two different cases:
\textit{(a)}: $k<i$ and \textit{(b)}: $k>i+1$. For the case \textit{(a)}, we see that there is some 
$x^\ast\in K_l\cap \Gamma_{a_k}$. Thus $\ell(x^\ast)=a_k+\lambda|x^\ast|^2$. 
As $f(0)>a_i$ and $f(x^\ast)=a_k<a_i$,  similar to the proof above, by
the intermediate value theorem,  we have that there is some $\xi\in (0,\,1)$ such that $f(\xi x^\ast)=a_i$.
Therefore, $\xi x^\ast\in \Gamma_{a_i}$ so that $\ell(\xi x^\ast)\leq f^M_{K^-_i}(\xi x^\ast)$.
This implies
\begin{equation}\label{Eq.Prf.05} 
	(1-\xi)\ell(0)+\xi \ell(x^\ast)\leq a_i+\lambda|\xi x^\ast|^2\,.
\end{equation}
As $\ell(0)\geq a_i$ and $\ell(x^\ast)=a_k+\lambda|x^\ast|^2$ so that \eqref{Eq.Prf.05} implies that
\begin{equation}\label{Eq.Prf.06}
	( 1-\xi)a_i+\xi\left(a_k+\lambda|x^\ast|^2\right)\leq a_i+\lambda|\xi x^\ast|^2
\end{equation}
that is 
\begin{equation}\label{Eq.Prf.07}
	\xi(1-\xi)\lambda|x^\ast|^2\leq \xi(a_i-a_k)\,. 
\end{equation}	
Thus we have found that for  $0<\xi<1$ 
\begin{equation}\label{Eq.Prf.08}
	\lambda(1-\xi)|x^\ast|^2\leq (a_i-a_k)\,. 
\end{equation}
Since $\lambda(1-\xi)|x^\ast|^2\geq \lambda(1-\xi)^2|x^\ast|^2\geq \lambda\delta_0^2$ 
and $a_i-a_k\leq a_m-a_0$, we have
$\lambda\delta_0^2\leq a_m-a_0$, which contradicts our assumption on $\lambda$.

\medskip
\noindent If the case \textit{(b)} occurs, we have $f(0)<a_{i+1}$ and $f(x^\ast)=a_k>a_{i+1}$. Again by the 
intermediate value theorem, there is some $\xi\in (0,\,1)$ such that $f(\xi x^\ast)=a_{i+1}$.
However note that here  the value of $f^M_{K^-_i}$ on $\Gamma_{a_{i+1}}$ is $a_{i}$.
Therefore a similar argument to that for case \textit{(a)} will lead to a contradiction. Thus in both cases
we have proved that $K_l\subset \Gamma_{a_i}\cup \Gamma_{a_{i+1}}$.

\medskip
\noindent Now we consider $C^u_\lambda(f^{-M}_{K^-_i})(0)=\co[\lambda|\cdot|^2-f^{-M}_{K^-_i}](0)$. 
Let $\hat{\ell}$ be the affine function such that $\hat{\ell}(x)\leq \lambda|x|^2-f^{-M}_{K^-_i}(x)$, 
$\hat{\ell}(0)=\co[\lambda|\cdot|^2-f^{-M}_{K^-_i}](0)$ and
let $K_u=\{x\in K,\; \hat{\ell}(x)= \lambda|x|^2-f^{-M}_{K^-_i}(x)\}$. Again we have $\hat{\ell}(0)\geq -a_i$
and we can also show that $K_u\subset  \Gamma_{a_i}\cup \Gamma_{a_{i+1}}$.
By the definition of the convex envelope, we have
\begin{equation}\label{Eq.Prf.09} 
	\begin{array}{l}
		\displaystyle \co[f^M_{K^-_i}+\lambda|\cdot|^2](0)\\
		\displaystyle \phantom{xxx} = \inf\left\{\sum^{n+1}_{k=1}
	\displaystyle \lambda_k\Big(f^M_{K^-_i}(x_k)+\lambda|x_k|^2\Big),\; x_k\in \mathbb{R}^n,\; 
	\displaystyle \lambda_k\geq 0,\;\sum^{n+1}_{k=1}\lambda_k=1,\;
	\displaystyle \sum^{n+1}_{k=1}\lambda_kx_k=0\right\}\\
			\phantom{xxx}\displaystyle = \inf\left\{\sum^{n+1}_{k=1}
	\displaystyle \lambda_k\Big(f^M_{K^-_i}(x_k)+\lambda|x_k|^2\Big),\; x_k\in K_l,\; 
	\displaystyle \lambda_k\geq 0,\;\sum^{n+1}_{k=1}\lambda_k=1,\;
	\displaystyle \sum^{n+1}_{k=1}\lambda_kx_k=0\right\}\\
		\displaystyle \phantom{xxx} = \inf\left\{\sum^{n+1}_{k=1}
	\displaystyle \lambda_k\Big(f^M_{K^-_i}(x_k)+\lambda|x_k|^2],\; x_k\in K_l\cup K_u,\; 
	\displaystyle \lambda_k\geq 0,\;\sum^{n+1}_{k=1}\lambda_k=1,\;
	\displaystyle \sum^{n+1}_{k=1}\lambda_kx_k=0\right\}\\
		\displaystyle \phantom{xxx} = a_i+
	\displaystyle \inf\left\{\sum^{n+1}_{k=1}\lambda_k\lambda|x_k|^2,\; x_k\in K_l\cup K_u,\; 
	\displaystyle \lambda_k\geq 0,\;\sum^{n+1}_{k=1}\lambda_k=1,\;
	\displaystyle \sum^{n+1}_{k=1}\lambda_kx_k=0\right\}\\
		\displaystyle \phantom{xxx} =: a_i+C_0\,.
\end{array}
\end{equation}
Similarly, we have $\co[\lambda|\cdot|^2-f^{-M}_{K^-_i}](0)=-a_i+C_0$, and  hence
\begin{equation}\label{Eq.Prf.10}
	A^M_\lambda(f_{K^-_i})(0)=\frac{1}{2}\Big(\co[f^M_{K^-_i}+\lambda|\cdot|^2](0)-
						  \co[\lambda|\cdot|^2-f^{-M}_{K^-_i}](0)
					\Big) = a_i\,.
\end{equation}
By using the same argument as above, we can also show that
$A^M_\lambda(f_{K^+_i})(0)=a_{i+1}$ and this proves \eqref{Eq.Prf.04}. \\

%%%%%%%%%%%%%%%%%%%%%%%%%%%%%%%%%%%%%%%%%%%%%%%%%%%%%%%%%%

%%%%%%%%%%%%%%%%%%%%%%%%%%%%%%%%%%%%%%%%%%%%%%%%%%%%%%%%%%

\noindent \textit{(iii):}  Suppose $f(0)<a_0$. If we let $\ell$ be the affine function such that 
$\ell(x)\leq f^M_K(x)+\lambda|x|^2$, $\ell(0)=\co[f^M_K+\lambda|\cdot|^2](0)$
and let $K_l=\{x\in \co[K],\; \ell(x)=f^M_K(x)+\lambda|x|^2\}$, then in this special
case we only need to show that $K_l\subset \Gamma_{a_0}$. As $a_0<a_1<\ldots<a_m$, we only need to
rule out one possibility that $K_l\cap \Gamma_i\neq\emptyset$ for any $0<i\leq m$. By following 
the arguments of the proof of $(ii)(b)$, we can show that $K_l\subset \Gamma_0$.
Similarly we can also show that $K_u\subset \Gamma_0$, where $K_u = \{ x \in \co[K], \hat{\ell}(x) = \lambda|x|^2 - f^{-M}_K(x)\}$ 
for the affine function $\hat{\ell}$ such that $\hat{\ell}(x) \leq \lambda|x|^2 - f^{-M}_K(x)$ and 
$\hat{\ell}(0) = \co[\lambda| \cdot|^2 - f^{-M}_K](0)$. 
The proof is then similar to that of part $(ii)$. 
Note that here we do not have to introduce functions $f^M_{K^+_0}$ and $f^M_{K^-_0}$ as in $(ii)$ given that the
condition we have is $f(0)<a_0$ while in $(ii)$ we had  $a_i<f(0)<a_{i+1}$. 
\end{proof}

%%%%%%%%%%%%%%%%%%%%%%%%%%%%%%%%%%%%%%%%%%%%%%%%%%%%%%%%%%%%%%%%%%%%%%%%%%%%%%%%%%%%%%%%%%%%%%
%%%%%%%%%%%%%%%%%%%%%%%%%%%%%%%%%%%%%%%%%%%%%%%%%%%%%%%%%%%%%%%%%%%%%%%%%%%%%%%%%%%%%%%%%%%%%%

\begin{proof}(\cref{Prop.AprxLS})
\textit{(i):} Without loss of generality, we may assume $x_0=0\in \Omega_i$. 
Now note  that Corollary \ref{Cor.AprxCnt}, applied with $f, r$ and $R$ given  by $\tilde{f}$, $R$ and $R+1$ respectively, gives that
\begin{equation}
\label{cor2.7-output}
| A^M_{\lambda}(\tilde{f}_{K_{R+1}})(0) - \tilde{f}(0) | \leq \tilde{\omega} \left( r_c(0) + \frac{\tilde{a}}{\lambda} + \sqrt{\frac{2 \tilde{b}}{\lambda}}\right).
\end{equation}
Then since $0 \in \Omega_i \subset V_{a_m}$, it follows that $\tilde{f}(0) = f(0)$, and also that $r_c(0) \leq d_i(0)$, 
by \eqref{Eq.Ineq01}. To prove \eqref{Eq.Prop.AprxLS.1}, it thus remains to show that 
$A^M_{\lambda}(\tilde{f}_{K_{R+1}})(0) = A^M_{\lambda}(f_K)(0)$. 
To see this, note first that by arguments similar to those in the proof of \cite[Theorem 3.7]{ZCO16}, we have that
\begin{equation}\label{Eq.Prf.AprxLS.01}
	C^l_\lambda(\tilde{f}^M_{K_{R+1}})(0)=\sum^{n^\ast}_{k=1}\lambda_k(\tilde{f}_{K_R}(x_k)+\lambda|x_k|^2)
\end{equation}
for some $2\leq n^\ast\leq n+1$, $\lambda_k>0$, $x_k\in K_{R+1}$, $k=1,2,\dots,n^\ast$, 
with $\Sigma_{k=1}^{n^{\ast}} \lambda_k = 1$ and $\Sigma_{k=1}^{n^{\ast}} \lambda_k x_k = 0$. 
Now if $x_k \in K$ for each $1 \leq k \leq n^{\ast}$, then $\tilde{f}_{K_{R+1}}(x_k) = f_K(x_k)$, and hence
$C^l_{\lambda}(\tilde{f}^M_{K_{R+1}})(0) = C^l_{\lambda}({f}^M_{K})(0)$. So suppose, for contradiction, 
that $x_{k_0} \in K_{R+1} \setminus K = B^c(0; R+1)$. Then there exists an affine function $\ell$ such that 
\[
	\ell(y) \leq \tilde{f}^M_{K_{R+1}}(y) + \lambda |y|^2 \;\; \mbox{for all}
	\;\; y \in \mathbb{R}^n, \;\;\; \ell(x_k) = \tilde{f}^M_{K_{R+1}}(x_k) + \lambda |x_k|^2, \;\; 1 \leq k \leq n^{\ast},
\]
so that 
\[
	\ell(x_{k_0}) = \tilde{f}^M_{K_{R+1}}(x_{k_0})  + \lambda | x_{k_0} |^2 = a_m+1 + \lambda | x_{k_0} |^2 .
\]
Since $\tilde{f}^M_{K_{R+1}}(y) = a_m+1$ for all $y \in B^c(0, R+1)$, $\ell$ must be the unique tangent 
plane to the function  $y \to a_m+1 + \lambda |y|^2$ at $y= x_{k_0}$, namely
\[
	\ell(y) = a_m+1 + \lambda | x_{k_0}|^2 + 2 \lambda x_{k_0} \cdot (y-x_{k_0}), \;\;\; y \in \mathbb{R}^n.
\]
Now it follows from the fact that  this plane does not touch the graph of $y \to a_m+1 + \lambda |y|^2$ 
at any other point that $x_k \not\in B^c(0, R+1)$ for $1 \leq k \leq n^{\ast}$, $k \neq k_0$, and hence, since
$ n^{\ast} \geq 2$, there must exist $x_{\hat{k}}$, $\hat{k} \neq k_0$, with $x_{\hat{k}} \in \Gamma_{a_j}$ 
for some $1 \leq j \leq m$ and $\ell(x_{\hat{k}}) = \tilde{f}^M_{K_{R+1}}(x_{\hat{k}}) + \lambda|x_{\hat{k}}|^2 = a_j + \lambda|x_{\hat{k}}|^2$. 
But then
\[
	a_m+1 + \lambda | x_{k_0}|^2 + 2 \lambda x_{k_0} \cdot x_{\hat{k}}- 2 \lambda | x_{k_0}|^2 = a_j + \lambda |x_{\hat{k}}|^2,
\]
and hence, since $x_{k_0} \in B^c(0; R+1)$ and $x_{\hat{k}} \in B(0, R)$, 
\[
	a_m - a_j + 1 = \lambda (|x_{\hat{k}}|^2 - 2 x_{k_0} \cdot x_{\hat{k}} + | x_{k_0}|^2) = \lambda | x_{\hat{k}} - x_{k_0}|^2 > \lambda,
\]
which contradicts the assumption on $\lambda$. 
Likewise, $C^u_{\lambda}(\tilde{f}^{-M}_{K_{R+1}})(0) = C^u_{\lambda}(\tilde{f}^{-M}_{K})(0)$, 
and hence $A^M_{\lambda}({\tilde{f}}_{K_{R+1}})(0) = A^M_{\lambda}(f_K)(0)$, as required.

\medskip

\noindent \textit{(ii):} The proof of the Lipschitz case follows similar arguments. 
\end{proof}

%%%%%%%%%%%%%%%%%%%%%%%%%%%%%%%%%%%%%%%%%%%%%%%%%%%%%%%%%%%%%%%%%%%%%%%%%%%%%%%%%%%%%%%%%%%%%%
%%%%%%%%%%%%%%%%%%%%%%%%%%%%%%%%%%%%%%%%%%%%%%%%%%%%%%%%%%%%%%%%%%%%%%%%%%%%%%%%%%%%%%%%%%%%%%

\begin{proof}(\cref{Thm.ULAprx})
Similar to the proof of Theorem \ref{Thm.LSaprx}\textit{(i)},
we fix $x_{j_0}\in K$ and let $f_\lambda(x)=\lambda|x-x_{j_0}|^2-f^{-M}_K(x)$.
Define $\ell(x) = -f(x_{j_0})$ for all $x \in \mathbb{R}^n$. Then $\ell$ is a constant function, so is  affine. 
Clearly $\ell(x_{j_0})=f_\lambda(x_{j_0})$. We need to prove that 
\begin{equation}\label{Eq.Ineq}
	\ell(x)\leq f_\lambda(x) 
\end{equation}
for all $x\in \mathbb{R}^n$ so that $\co[f_\lambda](x_{j_0})=\ell(x_{j_0})=-f(x_{j_0})$, 
hence $C^u_\lambda (f^{-M}_K)(x_{j_0})=f(x_{j_0})$.
Inequality \eqref{Eq.Ineq} is equivalent to
\begin{equation*}
	-f(x_{j_0})\leq \lambda|x-x_{j_0}|^2-f^{-M}_K(x),\quad x\in \mathbb{R}^n\,.
\end{equation*}

\noindent If $x\in \mathbb{R}^n\setminus K$, $f_K^{_M}(x)=-M$. Since $-f(x_{j_0}))<M< \lambda|x-x_{j_0}|^2+M$, we clearly have
$\ell(x)\leq f_\lambda(x)$ for all $x\in \mathbb{R}^n\setminus K$.
If $x_j\in K$ and $x_j\neq x_{j_0}$, we need to prove that
\begin{equation*}
	-f( x_{j_0})\leq  \lambda|x_j-x_{j_0}|^2-f(x_j),\quad \text{or equivalently, }\quad 
	f(x_j)-f( x_{j_0})\leq  \lambda|x_j-x_{j_0}|^2\,.
\end{equation*}
Since $\alpha=\min\{|x_i-x_j|,\; x_i,\,x_j\in K,\; x_i\neq x_j\}$, then if $\lambda>L/\alpha$, we have
\begin{equation*}
	f(x_j)-f( x_{j_0})\leq L|x_j-x_{j_0}|\leq \lambda \alpha|x_j-x_{j_0}|\leq \lambda|x_j-x_{j_0}|^2\,,
\end{equation*}
which completes the proof. 
\end{proof}

%%%%%%%%%%%%%%%%%%%%%%%%%%%%%%%%%%%%%%%%%%%%%%%%%%%%%%%%%%%%%%%%%%%%%%%%%%%%%%%%%%%%%%%%%%%%%%
%%%%%%%%%%%%%%%%%%%%%%%%%%%%%%%%%%%%%%%%%%%%%%%%%%%%%%%%%%%%%%%%%%%%%%%%%%%%%%%%%%%%%%%%%%%%%%

\begin{proof}(\cref{Lem.BndAf})
We may write $\ell_s(x)=a\cdot x+b$ with $a\in \mathbb{R}^n$ and
$b\in\mathbb{R}$. We see that $D\ell_s(x)=a$ and we need to give an estimate of $|a|$. 
Since we have $\ell_s(x_i)=f_S(x_i)$ and $|\ell_s(x_i)-\ell_s(x_1)|=|f_S(x_i)-f_S(x_1)|\leq L|x_i-x_1|$, 
we see that $|a\cdot(x_i-x_1)|\leq  L|x_i-x_0|$
for $i=1,2,\dots,k$. As $\dim(\co[S])=n,$ there are at least $n$-vectors, say $\{x_2-x_1,\,\ldots, x_{n+1}-x_1\}$,
which are linearly independent and hence form a basis of $\mathbb{R}^n$. If we let $\{e_1,\,\dots,e_n\}$ be any
orthonormal basis of $\mathbb{R}^n$, there is an $n\times n$ invertible matrix $A=(a_{ij})_{i,j = 1}^n$ such that
$e_i=\sum^{n}_{j=1}a_{ij}(x_{j+1}-x_1)$. Hence
\begin{equation*}
	|a\cdot e_i|\leq \sum^n_{j=1}| a_{ij}||a\cdot (x_{j+1}-x_1)|\leq L\left(\sum^n_{j=1}|a_{ij}|^2\right)^{1/2}
	\left(\sum^n_{j=1}|x_{j+1}-x_1|^2\right)^{1/2}\,.
\end{equation*}
Therefore, the Euclidean norm of $a$ satisfies $|a|\leq L |A|(\sum^n_{j=1}|x_i-x_0|^2)^{1/2}$, 
where $|A|$ denotes the Frobenius norm of the matrix $A$,
and can then take $C_s=|A|(\sum^n_{j=1}|x_i-x_0|^2)^{1/2}$, which completes the proof.
\end{proof}

%%%%%%%%%%%%%%%%%%%%%%%%%%%%%%%%%%%%%%%%%%%%%%%%%%%%%%%%%%%%%%%%%%%%%%%%%%%%%%%%%%%%%%%%%%%%%%
%%%%%%%%%%%%%%%%%%%%%%%%%%%%%%%%%%%%%%%%%%%%%%%%%%%%%%%%%%%%%%%%%%%%%%%%%%%%%%%%%%%%%%%%%%%%%%

\begin{proof}(\cref{Thm.AprxRC})
We prove the result for the upper transform. The proof of the lower
transform follows similar arguments.

\noindent Let us consider the affine function $\lambda r_s^2-\ell_s(x)$. For $x\in S$, clearly
\begin{equation}\label{Eq.Prf.Aprx01}
	\lambda r_s^2-\ell_s(x)=\lambda r_s^2-f_K(x)=\lambda|x-x_s|^2-f^{-M}_K(x)\,.
\end{equation}
If we can show that $\lambda r_s^2-\ell_s(x)<\lambda|x-x_s|^2-f^{-M}_K(x)$ for
$x\in\mathbb{R}^n\setminus S$, then one obtains 
\begin{equation}\label{Eq.Prf.Aprx02}
	\co[\lambda|(\cdot)-x_s|^2-f^{-M}_K](x)=\lambda r_s^2-\ell_s(x)
\end{equation}
for $x\in \co[S]$ and the proof for the upper transform then follows.

\medskip
\noindent We consider two different cases: $(i)$ $x\in K\setminus S$ and $(ii)$ $x\in\mathbb{R}^n\setminus K$.

\medskip
\noindent For the case $(i)$, let $x\in K\setminus S$. We need then to prove that
\begin{equation}\label{Eq.Prf.Aprx03}
	\lambda r_s^2-\ell_s(x)< \lambda|x-x_s|^2-f_K(x)\,,
\end{equation}
or, equivalently, that
\begin{equation}\label{Eq.Prf.Aprx04}
	\lambda r_s^2-\ell_s(x)+f_K(x)< \lambda|x-x_s|^2\,.
\end{equation}
We have the following estimates for the left hand side of \eqref{Eq.Prf.Aprx04}.
\begin{equation}\label{Eq.Prf.Aprx05}
	\begin{split}
		\lambda r_s^2-\ell_s(x)+f_K(x)	& \leq \lambda r_s^2+|\ell_s(x)-\ell_s(x_s)|+|\ell_s(x_s)|+A_0 \\[1.5ex]
						& \leq \lambda r_s^2+C_sL|x-x_s|+C_sLr_s+2A_0\,.
	\end{split}
\end{equation}
We have used the fact that for any $x^\ast\in S$,
\begin{equation}\label{Eq.Prf.Aprx06}
	|\ell_s(x_s)|\leq |\ell_s(x_s)-\ell_s(x^\ast)|+|\ell_s(x^\ast)|\leq C_sLr_s+A_0
\end{equation}
as $\ell_s(x^\ast)=f_K(x^\ast)$.
Therefore \eqref{Eq.Prf.Aprx04} holds if
\begin{equation}\label{Eq.Prf.Aprx07} 
	\lambda r_s^2+C_sL|x-x_s|+C_sLr_s+2A_0<\lambda|x-x_s|^2\,.
\end{equation}
Note that $|x-x_s|\geq r_s+\sigma_s$. Let us consider the function
\begin{equation}\label{Eq.Prf.Aprx08}
	g(t)=\lambda t^2-\lambda r^2_s-C_sLt-C_sLr_s-2A_0\,.
\end{equation}
If we can find conditions for $\lambda$ such that $g(r_s+\sigma_s)>0$ and $g^\prime(t)>0$
when $t\geq r_s+\sigma_s$, then \eqref{Eq.Prf.Aprx07} holds and \eqref{Eq.Prf.Aprx04} will be satisfied.

\medskip
\noindent We see that $g(r_s+\sigma_s)>0$ is equivalent to
\begin{equation}\label{Eq.Prf.Aprx09}
	\lambda[(r_s+\sigma_s)^2-r_s^2]>C_sL(2r_s+\sigma_s)+2A_0\,.
\end{equation}
This last inequality is equivalent to \eqref{Thm.AprxRC.2}. Thus \eqref{Eq.Prf.Aprx04} holds and thus 
$g(r_s+\sigma_s)>0$.

\medskip
\noindent Next we have $g^\prime(t)=2\lambda t-C_sL$. Since $g^\prime(t)$ itself is an increasing function,
we only need to show that $g^\prime(r_s+\sigma_s)>0$, which is equivalent to
\begin{equation}\label{Eq.Prf.Aprx10}
	\lambda>\frac{C_sL}{2(r_s+\sigma_s)}\,,
\end{equation}
which follows from  \eqref{Thm.AprxRC.2}. 
This completes the proof for  case \textit{(i)}.\\

\medskip
\noindent \textit{(ii)}: Let $x\in\mathbb{R}^n\setminus K$, hence $-f_K^{-M}(x)=M$.
We need to prove that
\begin{equation}\label{Eq.Prf.Aprx11}
	\lambda r_s^2-\ell_s(x)< \lambda|x-x_s|^2+M\,.
\end{equation}
Again we have
\begin{equation}\label{Eq.Prf.Aprx12}
	\lambda r_s^2-\ell_s(x)\leq  \lambda r_s^2+C_sL|x-x_s|+C_sLr_s+A_0\,.
\end{equation}
Therefore we prove \textit{(ii)} if 
\begin{equation}\label{Eq.Prf.Aprx13}
	\lambda r_s^2+C_sL|x-x_s|+C_sLr_s+A_0<\lambda|x-x_s|^2+M\,.
\end{equation}
Since \eqref{Thm.AprxRC.2} is satisfied, then by inspection it is easy to verify that \eqref{Eq.Prf.Aprx13}
holds for all non-negative numbers $|x-x_s|\geq 0$, which completes the proof.
\end{proof}

%%%%%%%%%%%%%%%%%%%%%%%%%%%%%%%%%%%%%%%%%%%%%%%%%%%%%%%%%%%%%%%%%%%%%%%%%%%%%%%%%%%%%%%%%%%%%%
%%%%%%%%%%%%%%%%%%%%%%%%%%%%%%%%%%%%%%%%%%%%%%%%%%%%%%%%%%%%%%%%%%%%%%%%%%%%%%%%%%%%%%%%%%%%%%

\begin{proof}(\cref{Lem.AfAprx})
\textit{(i):}
We see that both $p_+$ and $p_-$ are well-defined functions in $D$ and clearly 
$p_-(x)\leq v\leq p_+(x)$ for every $(x,v)\in \co[\Gamma_s]$. 
It is also easy to see that the two different expressions for 
$p_+(x)$ and respectively for $p_-(x)$ are equal.\\

\medskip
\noindent \textit{(ii):} Since $\co[\Gamma_s]$ is a convex polytope, we have, for any $x_1,\, x_2\in D$
and for every $0<t<1$, that
\begin{equation*}
	t(x_1,p_+(x_1))+(1-t)(x_2,p_+(x_2))=(tx_1+(1-t)x_2, tp_+(x_1)+(1-t)p_+(x_2))\in \co(\Gamma_s)
\end{equation*}
as both $D$ and $\co[\Gamma_s]$ are convex.
Furthermore,  by definition of $p_+$, $tp_+(x_1)+(1-t)p_+(x_2)\leq p_+(tx_1+(1-t)x_2)$.
Thus $p_+$ is concave in $D$, hence is continuous in $D$. Similarly we can show that $p_-$ is convex, 
hence continuous in $D$. 
Also $p_+$ and $p_-$ are both piecewise affine functions.
In fact, since $\co[\Gamma_s]$ is a convex polytope, $\co[\Gamma_s]$ has finitely many
closed $n$-dimensional faces. We may write $\partial \co[\Gamma_s]=\Gamma_+\cup\Gamma_-\cup\Gamma_0$,
where $\Gamma_+=\cup_{k=1}^m F^+_k$, $\Gamma_-=\cup_{j=1}^l F^-_j$
and $\Gamma_0=\cup_{r=1}^s F^0_r$ with $F^+_k$, $F_j^-$ and $F_r^0$ $n$-faces of
$\co[\Gamma_s]$. For $F^+_k$, there is an affine function $\ell^+_k:\mathbb{R}^n\to \mathbb{R}$ 
such that $\ell^+_k(x)=v$ if $(x,v)\in F^+_k$
and $\ell^+_k(x)>v$ if $(x,v)\in \co[\Gamma_s]\setminus (F^+_k)$. Similarly, for
$F^-_j$, there is an affine function $\ell^-_j:\mathbb{R}^n\to \mathbb{R}$ 
such that $\ell^-_j(x)=v$ if $(x,v)\in F^-_j$
and $\ell^-_k(x)<v$ if $(x,v)\in \co[\Gamma_s]\setminus (F^-_j)$. Every  $F^0_r$ is an
$n$-face whose normal vectors are in $\mathbb{R}^n\times\{0\}\subset \mathbb{R}^n\times\mathbb{R}$,
that is, $F^0_r$ is perpendicular to $D\times\{0\}$. Since the vertices of each $F^+_k$ 
are extreme points of $\co[\Gamma_s]$ and every point $x\in S$ is an extreme point of $\co[S]$ 
we see that for every extreme point $(x,v)$ of $\co[\Gamma_s]$, $x$ is an extreme point of $D$.
Let $D^+_k=P_{\mathbb{R}^n}(F^+_k)$ be the orthogonal projection from $F^+_k$ to $\mathbb{R}^n$, 
then $D^+_k$ is a convex polytope contained in $D$ whose vertices are all in $S$.
The projection $P_{\mathbb{R}^n}$ also maps relative boundary of $F^+_k$ 
to boundary of $D^+_k$, and the relative interior $F^+_k$ to interior of $D^+_k$. 
Also on $D^+_k$, $p_+(x)=\ell^+_k(x)$. Thus $p_+(\cdot)$ is affine on $D^+_k$.

\medskip
\noindent Similarly, for each $F^-_j$, we define  $D^-_j=P_{\mathbb{R}^n}(F^-_j)$. Then
the vertices of $D^-_j$ belong to $S$ and $p_-(x):=\ell^-_j(x)$ is affine on $D^-_j$.\\

%%%%%%%%%%%%%%%%%%%%%%%%%%%%%%%%%%%%%%%%%%%%%%%%%%%%%%%%%%%%%%%%%%%%%%

\medskip
\noindent \textit{(iii):} It is easy to see that 
$\mathring{D}^+_k\cap \mathring{D}^+_j=\emptyset$ and
$\mathring{D}^-_k\cap \mathring{D}^-_j=\emptyset$ for $k\neq j$. 
Next we show that $D=\cup^m_{k=1}D_k^+=\cup^l_{j=1}D_j^-$.

\medskip
\noindent If $\cup^m_{k=1}D_k^+\neq D$, there is an interior point
$x\in D\setminus \cup^m_{k=1}D_k^+$. By definition $(x,p_+(x))\in \partial \co[\Gamma_s]$
and we may assume that $(x,p_+(x))$ lies in the relative interior of an 
$n$-face $F\subset\partial \co[\Gamma_s]$. If $F$ is one of the $F^-_j$'s, 
this implies $p_+(x)=p_-(x)$. This cannot happen inside $D$. If $F$ is one 
of the $F^0_r$'s, then $D^0_r:=P_{\mathbb{R}^n}(F^0_r)$ is an $n-1$-dimensional polytope. 
If $E$ is the $(n-1)$-dimensional plane in $\mathbb{R}^n$ containing $D^0_r$,
then $D$ must lie on one side of $D^0_r$. Therefore $D^0_r\subset \partial D$,
hence $x$ is a boundary point of $D$. This contradicts our assumption that $x$ 
is an interior point of $D$. Thus $D=\cup^m_{k=1}D_k^+$. Similarly, we can show that
$D=\cup^l_{j=1}D_j^-$.

\medskip
\noindent The other conclusions also follow from the above arguments.
\end{proof}

%%%%%%%%%%%%%%%%%%%%%%%%%%%%%%%%%%%%%%%%%%%%%%%%%%%%%%%%%%%%%%%%%%%%%%%%%%%%%%%%%%%%%%%%%%%%%%
%%%%%%%%%%%%%%%%%%%%%%%%%%%%%%%%%%%%%%%%%%%%%%%%%%%%%%%%%%%%%%%%%%%%%%%%%%%%%%%%%%%%%%%%%%%%%%

\begin{proof}(\cref{Thm.AfAprx})
Since $\co[S]=\cup^m_{k=1}D^+_k$ and on each $D^+_k$, 
there is an affine function $\ell_k^+:\mathbb{R}^n\to \mathbb{R}$ such that
$\ell_k^+(x)=p^+_k(x)$ for $x\in D^+_k$ and $\ell_k^+(x)> f_K(x)$ for $x\in S^+_k$, where
$S^+_k$ is the set of extreme points of $D^+_k$ given by Lemma \ref{Lem.AfAprx} which is a subset of $S$.
Let $C_k^+>0$ be the constant given by Lemma \ref{Lem.BndAf} so that $|D\ell^+_k(x)|<C_k^+L\leq C_sL$.

\medskip
\noindent If we can show that $\co[\lambda|(\cdot)-x_s|^2-f^{-M}_K]=\lambda r_s^2-\ell^+_k(x)$ for $x\in D_k^+$,
the proof is finished. As in the proof of Theorem \ref{Thm.AprxRC}, we have to consider different cases.
If $x\in \mathbb{R}^n$ or $x\in \mathbb{R}^n\setminus K$ or $x\in K\setminus S$, 
the proof for the inequality
$\lambda r_s^2-l^+_k(x)\leq \lambda|x-x_s|^2-f_K^{-M}(x)$ is the same as that in the proof of Theorem \ref{Thm.AprxRC}.
The only new case we have to consider is for $x\in S\setminus S^+_k$.

\medskip
\noindent But for $x\in S\setminus S^+_k$, the above inequality is
\begin{equation}\label{Eq.Prf.AfAprx01} 
	\lambda r_s^2-\ell^+_k(x)\leq \lambda|x-x_s|^2-f_K(x)=\lambda r_s^2-f_K(x)\,,
\end{equation}
which is equivalent to $\ell^+_k(x)\geq f_K(x)$ as $S\subset \partial B(x_s;\,r_s)$. We also know from
Lemma \ref{Lem.AfAprx} that $\ell^+_k(x)> f_K(x)$ for $x\in S\setminus S^+_k$. Therefore on each $D^+_k$,
\eqref{Eq.Prf.AfAprx01} holds as $p_+(x)=\ell^+_k(x)$ on $D^+_k$. The proof for the lower transform is similar.
The proof is finished.
\end{proof}

%%%%%%%%%%%%%%%%%%%%%%%%%%%%%%%%%%%%%%%%%%%%%%%%%%%%%%%%%%%%%%%%%%%%%%%%%%%%%%%%%%%%%%%%%%%%%%
%%%%%%%%%%%%%%%%%%%%%%%%%%%%%%%%%%%%%%%%%%%%%%%%%%%%%%%%%%%%%%%%%%%%%%%%%%%%%%%%%%%%%%%%%%%%%%

\begin{proof}(\cref{Cor.Inp})
For the proof of this result, we first follow the proof of 
Theorem \ref{Thm.AprxUF.Cmpct} so that the points
$x^i$'s for the convex envelope are in $\bar\Omega$. Then we follow the proof of \cite[Theorem 3.7]{ZCO16}
to show that $x^i$'s can only be in $K$. The rest of the proof then follows from that of Theorem \ref{Thm.AprxUF.Cmpct}.
\end{proof}

%%%%%%%%%%%%%%%%%%%%%%%%%%%%%%%%%%%%%%%%%%%%%%%%%%%%%%%%%%%%%%%%%%%%%%%%%%%%%%%%%%%%%%%%%%%%%%
%%%%%%%%%%%%%%%%%%%%%%%%%%%%%%%%%%%%%%%%%%%%%%%%%%%%%%%%%%%%%%%%%%%%%%%%%%%%%%%%%%%%%%%%%%%%%%

\section*{Acknowledgements}
The authors are extremely grateful to the anonymous referees, whose constructive comments 
on earlier versions of the manuscript have contributed to produce a better version of the paper.
The authors would also like to thank the Isaac Newton Institute for Mathematical Sciences for
support and hospitality during the programme `Variational methods and effective algorithms 
for imaging and vision'  when part of the work on this paper was undertaken. 
This work was partially supported by EPSRC  Grant Number EP/K032208/1.
KZ wishes then to thank The University of Nottingham for its support, whereas
EC is grateful for the financial support of the College of Science, Swansea University,
and
AO acknowledges the partial financial support of the Argentinian Research Council (CONICET)
through the project PIP 11220170100100CO,
the National University of Tucum\'{a}n through the project PIUNT CX-E625 and the FonCyT 
through the project PICT 2016 201-0105 Prestamo Bid. The authors finally would like to thank 
Simone Parisotto for delight discussions on the numerical simulations of surface reconstructions 
from real data.
%%%%%%%%%%%%%%%%%%%%%%%%%%%%%%%%%%%%%%%%%%%%%%%%%%%%%%%%%%%%%%%%%%%%%%%%%%%%%%%%%%%%%%%%%%%%%%
%%%%%%%%%%%%%%%%%%%%%%%%%%%%%%%%%%%%%%%%%%%%%%%%%%%%%%%%%%%%%%%%%%%%%%%%%%%%%%%%%%%%%%%%%%%%%%
%%%%%%%%%%%%%%%%%%%%%%%%%%%%%%%%%%%%%%%%%%%%%%%%%%%%%%%%%%%%%%%%%%%%%%%%%%%%%%%%%%%%%%%%%%%%%%

\bibliographystyle{siamplain}
\bibliography{references}

%%%%%%%%%%%%%%%%%%%%%%%%%%%%%%%%%%%%%%%%%%%%%%%%%%%%%%%%%%%%%%%%%%%%%%%%%%%%%%%%%%%%%%%%%%%%%%
%%%%%%%%%%%%%%%%%%%%%%%%%%%%%%%%%%%%%%%%%%%%%%%%%%%%%%%%%%%%%%%%%%%%%%%%%%%%%%%%%%%%%%%%%%%%%%
%%%%%%%%%%%%%%%%%%%%%%%%%%%%%%%%%%%%%%%%%%%%%%%%%%%%%%%%%%%%%%%%%%%%%%%%%%%%%%%%%%%%%%%%%%%%%%

\end{document}